\documentclass[12pt, twoside
]{article}

\usepackage{graphpap}
\usepackage{epsfig}

\usepackage{amscd}

\usepackage{amsmath}
\usepackage{amssymb}
\usepackage{amsthm}

\makeatletter
\def\@seccntformat#1{\csname the#1\endcsname{.}\hskip .5em}

\renewcommand{\section}{\@startsection {section}{1}{\z@}%
                                {-4.5ex \@plus -1ex \@minus-.2ex}%
                                   {2.3ex \@plus.2ex}%
                               {\reset@font\Large\scshape\centering}}

\makeatother

\normalbaselines

\newcommand{\C}{\mathbb{C}}

\newcommand{\e}{\varepsilon}

\newcommand{\dist}{\operatorname{dist}}
\newcommand{\s}{\sigma}
\renewcommand{\le}{\leqslant}
\renewcommand{\ge}{\geqslant}

\theoremstyle{definition}

\theoremstyle{remark}

\newtheorem*{rem*}{Remark}

\numberwithin{equation}{section}

\newcommand{\ci}[1]{_{ {}_{\scriptstyle #1}}}
\newcommand{\cci}[1]{_{ {}_{\scriptscriptstyle #1}}}

\newcommand{\f}{\varphi}

\makeatletter

\def\multilimits@{\bgroup
  \Let@
  \restore@math@cr
  \default@tag
 \baselineskip\fontdimen10 \scriptfont\tw@
 \advance\baselineskip\fontdimen12 \scriptfont\tw@
 \lineskip\thr@@\fontdimen8 \scriptfont\thr@@
 \lineskiplimit\lineskip
 \vbox\bgroup\ialign\bgroup\hfil$\m@th\scriptstyle{##}$\hfil\crcr}
\def\Sb{_\multilimits@}
\def\Sp{^\multilimits@}
\def\endSb{\crcr\egroup\egroup\egroup}

\makeatother

\newcommand{\Dl}{\Delta}
\newcommand{\D}{\mathcal{D}}
\newcommand{\al}{\alpha}
\newcommand{\cl}{\operatorname{clos}}
\newcommand{\supp}{\operatorname{supp}}

\newcommand{\la}{\lambda}

\newcommand{\diam}{\operatorname{diam}}

\newcommand{\Om}{\Omega}
\newcommand{\om}{\omega}

\newcommand{\E}{\mathbb{E}}

\newcommand{\wt}{\widetilde}

{\end{list}}

\begin{document}
\title{The $Tb$-theorem on non-homogeneous spaces that proves a conjecture of Vitushkin}

\author{F. Nazarov, S. Treil
\ and A. Volberg\thanks{All  authors are partially supported by the NSF grant
DMS 9970395}}

\date{}
\maketitle

\tableofcontents

\bigskip
 
\centerline{F.\,Nazarov,\, S.\, Treil,\, A.\, Volberg}

\begin{abstract}
This article was written in 1999, and was posted as a preprint in CRM (Barcelona) preprint series $n^0\, 519$ in 2000.
However, recently CRM (Barcelona) erased all preprints dated before 2006 from its site, and this paper became inacessible. It has certain importance though, as the reader shall see. Meanwhile this paper in bits and pieces appeared in several book formats, namely in Volberg's lecture notes [Vo], in Doudziak's book [Du], and in Tolsa's book [To].

Formally this paper is a proof of the (qualitative version of the) Vitushkin
conjecture. The last section is concerned with the quantitative version. This quantitative version turns out to be very important.
It allowed Xavier Tolsa to close the subject concerning Vtushkin's conjectures: namely, using the quantitative nonhomogeneous 
$Tb$ theorem proved in the present paper, he proved the semiadditivity of analytic capacity.
Another ``theorem'', which is implicitly contained in this paper,
is the statement that any non-vanishing
$L^2$-function is accretive in the sense that
if one has a finite measure $\mu$ on the complex plane $\C$ that is
Ahlfors at almost every point (i.e. for $\mu$-almost every $x\in\C$ there
exists a constant $M>0$ such that $\mu(B(x,r))\le Mr$ for every $r>0$) then
any one-dimensional antisymmetric
Calder\'on-Zygmund operator $K$ (i.e. a Cauchy integral type
operator) satisfies the following ``all-or-nothing'' princple:
if there exists at least one function $\f\in L^2(\mu)$ such that $\f(x)\ne 0$
for $\mu$-almost every $x\in\C$ and such that {\it the maximal singular operator} $K^*\f\in L^2(\mu)$, then there
exists an everywhere positive weight $w(x)$, such that $K$ acts from
$L^2(\mu)$ to $L^2(w d\mu)$. In particular, there exists a a set $E$ of positive $\mu$-measure, $\mu(E)>0$, such that operator $K$ is a bounded operator from $L^2(E, \mu)$ to itself. Moereover,  a concrete estimate can be given for the bound of its norm and the portion $\mu(E)/\|\mu\|$ if we have quantitative information on how non-zero is $\f$ and haow small is $\|K^*\f\|$.

\end{abstract}

\noindent{\bf Table of contents}

\medskip

\def\tt#1#2{\leftline{\bf #1\hfill #2}}

\tt{0. What  this is all about}{}
\tt{I. Suppressed operators $K\ci\Phi$}{}
\tt{II. Perfect random dyadic lattices and good functions}{}
\tt{III. Perfect hair}{}                                            
\tt{IV. Truncated mathematical expectation}{}                       
\tt{V. How to use perfect hair}{}                                   
\tt{VI. Lyric deviation: Hausdorff measure and Analytic capacity}{}
\tt{VII. Cauchy integral representation}{}                          
\tt{VIII. The Ahlfors radius $\mathcal{R}(x)$}{}                            
\tt{IX. The exceptional set $H$}{}                                      
\tt{X. Localization}{}                                              
\tt{XI. Construction of perfect hair}{}                             
\tt{XII. Projections $\Lambda$ and $\Dl_Q$}{}                           
\tt{XIII. Functions $\Phi\ci\D$}{}                                  
\tt{XIV. Action on good functions}{}                                
\tt{XV. Estimation of $\sigma_2$}{}                                 
\tt{XVI. Estimation of $\s_3$}{}                                    
\tt{XVII. Estimation of $\s_3^{term}$}{}                           
\tt{XVIII. Whitney decomposition}{}                                  
\tt{XIX. Estimation of $\s_3^{tr}$}{}                                
\tt{XX. Estimation of $\s_1$}{}                                     
\tt{XXI. Negligible contours}{}                                     
\tt{XXII. Estimation of probability}{}                              
\tt{XXIII. Quantitative pulling ourselves up by the hair}{}         
\tt{XXIV. The quantitative version of Vitushkin's conjecture}{} 
\tt{XXV. Cotlar's inequality for non-uniformly Ahlfors measures}{}

\bigskip

\bigskip

{\bf 0. What this is all about}

\medskip

Let us be a little bit more specific. The analytic capacity of a compact set on the plane was defined by Ahlfors in 1947 as

$$
\gamma(E) =\sup_f\lim_{z\rightarrow \infty} |z\,f(z)|,
$$
where the supremum is taken over all analytic functions in the complement of $E$ such that $|f(z)|\leq 1$ 
and $f(\infty) =0$. Ahlfors showed that $\gamma(E)=0$ if and only if $E$ is removable for bounded analytic
functions. It was very interesting to find a geometric characterization. This is often called the Painlev\'e problem since
Painlev\'e started to study it more than 100 years ago.

\bigskip

Vituskin's conjecture (1967): for sets $E$ such that $\mathcal{H}^1(E) <\infty$, $\gamma(E) =0$ if and only
if 
$\mathcal{H}^1(E\cap \Gamma)=0$ for every rectifiable curve $\Gamma$.

\bigskip

Alberto Calder\'on and Guy David found the geometric characterization of sets of positive analytic 
capacity and finite length (= finite $\mathcal{H}^1$-measure), thus proving one half of Vitushkin's conjecture each. 

\bigskip

{\bf Theorem:}
Let $E$ be a compact on the plane with $\mathcal{H}^1(E) <\infty$. Then $\gamma(E) =0$ if and only if 
$\mathcal{H}^1(E\cap \Gamma)=0$ for every rectifiable curve $\Gamma$.

\bigskip

Here $\mathcal{H}^1$ is $1$-dimensional Hausdorff measure. The sets of finite $1$-dimensional Hausdorff measure with the 
latter condition satisfied are called purely unrectifiable according to Federer. Besicovitch studied them and multidimensional
analogs in the 1920's and 1930's and proved many very difficult and beautiful results about such sets. He called them irregular. 

The ``only if" part of the theorem has been proved by Calder\'on in 1977. It amounts to establishing that the Cauchy integral 
operator on Lipschitz curves is bounded on $L^2$ (Calder\'on's problem, which he solved in 1977 for small Lipschitz constants:
this turned out to be sufficient for the ``only if" part). The ``if" part was considered to be super difficult. Finally it was
proved by Guy David in 1997 [D1] using also [DM]. But actually this was only the ``analytic part" of the proof. The ``geometric
part" was fortunately known because of the fantastic idea of Melnikov and Verdera [MV] and a geometric theorem due to David and
L\'eger [L].

\bigskip

Here we give another (probably simpler and more streamlined, more conceptual) proof of the ``if" part in the theorem, actually of 
the ``analytic" part.

\bigskip

To explain the approach we need the notion of the Cauchy integral operator. So let $E$ in the plane have finite $\mathcal{H}^1(E)$.
Call $\mu=\mathcal{H}^1|E$. 
The Cauchy singular integral operator $C_{\mu}$ is

$$
C_{\mu}g(z) =\lim_{\delta\rightarrow 0} \int_{E\setminus B(z,\delta)}\frac{g(\zeta)}{\zeta-z}\,d\mu(\zeta)\,.
$$

Actually, if $z\in E$, it is not clear when the limit exists (while outsude of $E$ the definition is always fine).
So we introduce the maximal Cauchy singular integral operator $C^*_{\mu}$:

$$
C^*_{\mu}g(z) =\sup_{\delta > 0}| \int_{E\setminus B(z,\delta)}\frac{g(\zeta)}{\zeta-z}\,d\mu(\zeta)|
$$
and the ``cut-off" Cauchy singular integral operator $C_{\mu}^{\delta}$:

$$
C_{\mu}^{\delta}g(z) = \int_{E\setminus B(z,\delta)}\frac{g(\zeta)}{\zeta-z}\,d\mu(\zeta)\,.
$$

\bigskip 

Suppose $\gamma(E) >0$. One should find a rectifiable $\Gamma$ such that $H^1(E\cap \Gamma)>0$.

\bigskip

The analytic part here will end by constructing a positive (this is very important, let us say this again, positive) $\phi$ such
that
$$
 |C^*_{\mu}\phi(z)|\leq 1 \quad \forall z\in \C\,.
$$

Setting $\nu =\phi\,d\mu$ and applying to this positive measure the  permutation idea from [MV] one gets

$$
c^2(\nu):= \int\int\int c(x,y,z)^2d\nu(x)d\nu(y)d\nu(z) <\infty
$$
where $c(x,y,z)$ is the reciprocal of the radius of the circle passing through $x,y,z$. The quantity $c(\nu)$ is  
called the Menger curvature of the  measure $\nu$.

\bigskip

The following theorem is from the abovementioned ``geometric part" of the proof. It is due to David and L\'eger [L].

\bigskip

{\bf Theorem:}
If $\nu=\phi\,dH^1|E, \phi \ge 0, \phi\in L^{\infty}(E), \,\mathcal{H}^1(E) <\infty$ and $c^2(\nu) <\infty $, then there are 
rectifiable curves $\Gamma_i$ such that $\nu(\C\setminus \cup_{i=1}^{\infty} \Gamma_i) =0$.

\bigskip

Now we see that after constructing a {\bf positive} $\phi$ such that $|C^*_{\mu}\phi(z)|\leq 1, \quad \forall z\in \C$, 
one refers to the geometric papers [MV] and [L] to finish the proof of Vitushkin's conjecture.

How to find such a positive $\phi$?  We have only the information that $\gamma(E) > 0$ and $\mathcal{H}^1(E) <\infty$. 
The first condition means that there is a nonconstant bounded analytic function $f$ in $\C\setminus E$ vanishing at infinity. The
second condition quite easily shows that  this $f$ is represented as a Cauchy integral of $\phi \,d\,H^1|E=\phi \,d\,\mu$: $ f(z)
= C_{\mu}\phi(z),\, \forall z\in \C\setminus E$. We do not explain this now, but it is easy to assume that our
$\mu:=\mathcal{H}^1|E$ satisfies $\mu(B(z,r)) \leq C\,r$ for all $z\in \C$ and all $r>0$. Then not only $C_{\mu}\phi(z)$ is bounded
on $\C\setminus E$, but one can prove that there exists a finite constant $C$ such that 

$$
 |C^*_{\mu}\phi(z)|\leq C <\infty \quad\quad\quad \forall z\in \C\,.
$$

But this is not at all what we need---even though it seems precisely what we wanted. The main problem is that $\phi$ is 
complex valued function! It is impossible to prove that it is positive. (Actually positivity will generically never happen.)

\bigskip

Here is the main result to the proof of which the rest of the paper is devoted:

\bigskip

${}$

\bigskip

{\bf Main Theorem:}
Let $\mu$ denote $\mathcal{H}^1|E$ for a set $E$ of finite $1$-dimensional Hausdorff measure.
If there is a nonzero $\phi\in L^{\infty}(E)$ (this $L^{\infty}$ part can be weakened) such that 
$\sup_{z\in\C}|C^*_{\mu}\phi(z)|\leq Const <\infty$, then there exists a nonnegative bounded function $\psi$, which is strictly 
positive on the set of positive measure $\mu$, such that $\sup_{z\in\C}|C^*_{\mu}\psi(z)|\leq Const <\infty$.

\bigskip

Actually the fact that we work with $\mathcal{H}^1$ is not important. Another way of expressing the essence of the previous 
theorem is to formulate its analog, which is as follows:

\bigskip

{\bf Theorem (on bounded Cauchy transforms of measures):}
Let $\nu$ denote a nonzero complex measure with compact support on the plane. Let its Cauchy transform $C_{\nu}$ be uniformly 
bounded: 
$\sup_{z\in \C\setminus \supp(\nu)}|C_{\nu}(z)|\leq Const <\infty$. Suppose that the area of $\supp\nu$ is zero. Then there 
exists a  positive  measure $\mu$,  absolutely continuous with respect to $\nu$, such that its Cauchy transform is uniformly
bounded too: $\sup_{z\in \C\setminus \supp(\mu)}|C_{\mu}(z)|\leq Const <\infty$.

\bigskip

We are grateful to V. Lomonosov and   N.K. Nikolski  who pointed out to us that this result has the following interpretation as 
a result about normal operators.

\bigskip

{\bf Theorem (on resolvents of normal operators):}
Let $N$ be a normal operator whose spectrum $\sigma(N)$ has zero area. Let $R_{\lambda}, \lambda\in \C\setminus \sigma(N)$, 
denote its resolvent. If there are two vectors  $f,g$ such that $g$ belongs to the closed linear span of $\{N^k f\}_{k\ge 0}$,
$g\neq 0$, and such that $(R_{\lambda}f,g)$ is a bounded function on $\C\setminus \sigma(N)$, then there exists a nonzero vector
$h$ in the closed linear span of $\{N^k f\}_{k\ge 0}$, such that $(R_{\lambda}h,h)$ is a bounded function on $\C\setminus
\sigma(N)$.

\bigskip

In other words, if a compact set supports a complex measure with bounded nonzero Cauchy transform, 
then this compact set supports a positive measure with bounded (and also automatically nonzero) 
Cauchy transform. Also if the
resolvent of a normal operator is uniformly bounded on a pair of vectors $f,g$, ($g\neq 0$ 
being in the invariant subspace
generated by $f$) then it is uniformly bounded on certain $h,h$, $h\neq 0$.

\bigskip

So this is what we will be proving  using the ``perfect hair" approach in what now follows.

Few words about methods used in the proof. 

The probabilistic argument is a very important thing here. 
It is used to compensate for the roughness of our underlying
measure. The other people have used before the arguments involving many dyadic lattices at once. We
mean  a paper by Garnett and Jones called "BMO from diadic BMO" [GJ].

We use dyadic
martingale decomposition in our proof. We want to mention that looking at dyadic
martingale decomposition   is also a variation of an old theme, initiated, at least
in the context of the Cauchy integral, by Coifman, Jones and Semmes in their paper [CJS].
There they proved a
$T(b)$ theorem for the Cauchy integral using a Haar basis adapted to $b$. 
The main strategy of our proof is looking at dyadic
martingale decomposition, but a random one! 

\bigskip

{\bf Going further.}

Let us recall the definitions of the Cauchy capacities.
The first is {\it the complex Cauchy capacity} (not a very good name because it is a non-negative set function). We define it for
$\nu\in M_c(K)$:= complex measures supported on $K$.
$$
\gamma_c(K):=\{\sup|\nu(K)|: |C^{\nu}(z)|\leq 1 \forall z\in\C\setminus K,\,\nu\in M_c(K)\}\,.
$$
The second is {\it the positive Cauchy capacity} or just {\it the Cauchy capacity}:
$$
\gamma_c(K):=\{\sup \mu(K): |C^{\mu}(z)|\leq 1 \forall z\in\C\setminus K,\,\nu\in M_+(K)\}\,.
$$
Here $M_+(K)$ is a set of all positive measures supported on $K$.
Obviously,
$$
\gamma_+(K)\leq \gamma_c(K)\leq \gamma(K)\,.
$$
We actually  prove in this paper the following theorem (a sort of inverse to the previous left inequality).

\bigskip

{\bf Theorem.}
Let $K$ be a compact set of zero area. Then
$$
\gamma_+(K) \geq A \Bigl(1+\Bigl(\frac{\diam K}{\gamma_c(K)}\Bigr)^2\Bigl(\frac{\|\nu\|}{\gamma_c(K)}\Bigr)^{42}\Bigr)^{-1/2}
\gamma_c(K)\,,\qquad\qquad\qquad\qquad (INV)
$$
where $\nu$ is a measure that (almost) gives the supremum in the definition of $\gamma_c$. 
Its total variation in (INV) hinders us
from proving that 
$$
\gamma_c \geq A \gamma_+\,.
$$

\bigskip

\noindent Recently Xavier Tolsa [XT3] used (INV) and a very clever ``induction on scales" that appeared  in the
preprint  by J. Mateu, X. Tolsa and J. Verdera [MTV],
in which it is shown  that the condition conjectured by Mattila characterizes the 
Cantor sets of vanishing analytic capacity, to prove:
. 
$$
\gamma_c \geq A \gamma_+\,.
$$

\noindent This solves an old open problem. Actually, this implies the positive answer to Vitushkin's question whether the analytic
capacity is semi-additive (with absolute constant). In fact, it is relatively easy to prove that $\gamma_+$ is semi-additive (see
[NTV2], [NTV3]). The uniform comparability of
$\gamma_c$ and
$\gamma_+$ implies uniform comparability of $\gamma$ and
$\gamma_+$ (indeed, this is just an easy approximation argument using the fact that for any compact set which is
a finite union of
rectifiable curves,
$\gamma_c$ coincides with
$\gamma_+$).

\bigskip

{\bf Acknowledgements.} We are grateful to Michael Frazier and Joan Verdera for many helpful remarks.

\bigskip

{\bf I. Suppressed operators $K\ci\Phi$}

\bigskip

Let $\Phi$ be a nonnegative Lipschitz function, i.e., $\Phi(x)\ge 0$ for
every $x\in\C$ and
$$
|\Phi(x)-\Phi(y)|\le |x-y|\text{\qquad for every } x,y\in\C.
$$
Define
$$
k\ci\Phi(x,y)=\frac{\overline{x-y}}{|x-y|^2+\Phi(x)\Phi(y)}\,.
$$
{\bf Lemma:}
The kernel $k\ci\Phi$ is an antisymmetric Calderon-Zygmund kernel.
It is also really well suppressed at the points where $\Phi(x)>0$
or $\Phi(y)>0$. Namely,
$$
|k\ci\Phi(x,y)|\le \frac{1}{\max\{\Phi(x),\Phi(y)\}}\qquad\text{ for all }
x,y\in\C.
$$

{\bf Proof:}
Clearly,
$$
|k\ci\Phi(x,y)|\le\frac{1}{|x-y|} \text{\qquad and \qquad}
k\ci\Phi(x,y)=-k\ci\Phi(y,x).
$$

Since $k\ci\Phi$ is antisymmetric, to prove the second claim of the lemma,
it is enough to show that $|k\ci\Phi(x,y)|<\frac{1}{\Phi(x)}$
for all $x,y\in\C$.
We have
$\Phi(y)\ge\Phi(x)-|x-y|$. Therefore
$$
|k\ci\Phi(x,y)|\le\frac{|x-y|}{|x-y|^2+\Phi(x)(\Phi(x)-|x-y|)}
=\frac{|x-y|}{|x-y|^2+\Phi(x)^2-\Phi(x)|x-y|}
$$
$$
=\frac{|x-y|}{\Phi(x)|x-y|+(\Phi(x)-|x-y|)^2}\le \frac{1}{\Phi(x)},
$$
and we are done.

To prove the first claim of the lemma, let us show that
$$
|\nabla_x k\ci\Phi(x,y)|\le \frac{4}{|x-y|^2}.
$$
Indeed,
$$
|\nabla_x k\ci\Phi(x,y)|\le \frac{1}{|x-y|^2+\Phi(x)\Phi(y)}+
\frac{2|x-y|^2+|x-y|\Phi(y)}{[|x-y|^2+\Phi(x)\Phi(y)]^2}
$$
$$
\le
\frac{3}{|x-y|^2}+\frac{|x-y|\Phi(y)}{[|x-y|^2+\Phi(x)\Phi(y)]^2}
=
\frac{3}{|x-y|^2}+\frac{\Phi(y)}{[|x-y|^2+\Phi(x)\Phi(y)]}|k\ci\Phi(x,y)|
$$
$$
\le
\frac{3}{|x-y|^2}+\frac{\Phi(y)}{|x-y|^2}\frac{1}{\Phi(y)}=\frac{4}{|x-y|^2},
$$
finishing the proof of the lemma.

\bigskip

From now on, we will denote by $K\ci\Phi$ the operator with kernel
$k\ci\Phi$.

Pick some 
very small number $\delta>0$. It will stay fixed throughout the rest of the paper
and will be used in many formulae without any special comment. The reader may
think that $\delta$ is just an abbreviation for $45^{-239}$.

\medskip

{\bf II. Perfect random dyadic lattices and good functions}

\bigskip

Let $\mu$ be a measure   on the complex plane $\C$
satisfying $0<\mu(\C)<+\infty$.

Assume that $\D$ is a random dyadic lattice
(this phrase means that
we have a {\it family} of dyadic lattices
endowed with some probability $P$, and we use
the letter
$\D$ to denote {\it an element} in the family),
and let $\Lambda$, $\{\Delta\ci Q\}\ci{Q\in\D}$ be the (random) family of
projections associated with $\D$.
As usual, this means that
$$
\Lambda,\Delta\ci Q: L^2(\mu)\to L^2(\mu), \quad
\Delta\ci Q\Lambda=\Lambda\Delta\ci Q =0\text{ for all }Q\in\D, \qquad
\Delta\ci Q\Delta\ci R=0
\text{ when }Q\ne R,
$$
and for every
function $\f\in L^2(\mu)$, one has
$$
\f=\Lambda\f+\sum_{Q\in\D}\Delta\ci Q\f,
$$
where the series converges at least in $L^2(\mu)$.
Assume also that for every $\f\in L^2(\mu)$,
$$
2^{-1}||\f||^2\ci{L^2(\mu)}
\le
||\Lambda\f||^2\ci{L^2(\mu)}+
\sum_{Q\in\D}||\Delta\ci Q\f||^2\ci{L^2(\mu)}\le 2
||\f||^2\ci{L^2(\mu)}
.
$$

{\bf Remark:}

Let us make a couple of useful observations
 about such families of projections.

First of all, note that for every sequence of complex numbers
$\{c\ci Q\}\ci{Q\in\D}$ that is finite in the sense that only finitely many
$c\ci Q$ do not vanish, we have
$$
2^{-1}
\sum_{Q\in\D}|c\ci Q|^{2}||\Delta\ci Q\f||^2\ci{L^2(\mu)}
\le
\Bigl\|
\sum_{Q\in\D} c\ci Q\Delta\ci Q\f
\Bigr\|^2\ci{L^2(\mu)}
\le 2
\sum_{Q\in\D}|c\ci Q|^2||\Delta\ci Q\f||^2\ci{L^2(\mu)}.
$$

Indeed, consider the function $\wt\f:=\sum_{Q\in\D} c\ci Q\Delta\ci Q\f$
and note that $\Lambda\wt\f=0$, $\Dl\ci Q\wt\f=c\ci Q\Delta\ci Q\f$.
Now it remains only to apply our assumption to the function $\wt\f$
instead of $\f$ itself.

Now take any function $\psi\in L^2(\mu)$.
We have
$$
\Bigl|\sum_{Q\in\D} c\ci Q\langle \Delta\ci Q\f, \psi\rangle
\Bigr|=
\Bigl|
\Bigl<\sum_{Q\in\D} c\ci Q\Delta\ci Q\f,\psi\Bigr>\Bigr|
$$
$$
\le
\Bigl\|
\sum_{Q\in\D} c\ci Q\Delta\ci Q\f
\Bigr\|\ci{L^2(\mu)}\|\psi\|\ci{L^2(\mu)}
\le \sqrt2 \|\psi\|\ci{L^2(\mu)}
\Bigl[\sum_{Q\in\D}|c\ci Q|^2||\Delta\ci Q\f||^2\ci{L^2(\mu)}\Bigr]^{\frac12}.
$$
In particular, this means that if $\mathcal{F}\subset\D$ is some family of dyadic
squares, then
$$
\sum_{Q\in\mathcal{F}} |\langle \Delta\ci Q\f, \psi\rangle|
\le \sqrt2 \|\psi\|\ci{L^2(\mu)}
\Bigl[\sum_{Q\in\mathcal{F}}||\Delta\ci Q\f||^2\ci{L^2(\mu)}\Bigr]^{\frac12}
$$
(just take $c\ci Q=0$ for $Q\notin\mathcal{F}$ and choose $c\ci Q$ for $Q\in\mathcal{F}$
in such a way that $|c\ci Q|=1$ and $c\ci Q\langle \Delta\ci Q\f,
\psi\rangle=|\langle \Delta\ci Q\f, \psi\rangle|$; if the family $\mathcal{F}$ is
infinite, do it for all its finite subfamilies and then pass to the supremum).

Also, let us take any finite family $\mathcal{F}\subset\D$ such that
$\|\Dl\ci Q\f\|\ci{L^2(\mu)}>0$ for every $Q\in \mathcal{F}$.
Take $c\ci Q=0$ for $Q\notin\mathcal{F}$ and choose $c\ci Q$ for $Q\in\mathcal{F}$
in such a way that $|c\ci Q|=\dfrac
{|\langle \Delta\ci Q\f, \psi\rangle|}{||\Delta\ci Q\f||^2\ci{L^2(\mu)}}$
and $c\ci Q\langle \Delta\ci Q\f,
\psi\rangle=\dfrac{|\langle \Delta\ci Q\f, \psi\rangle|^2}
{||\Delta\ci Q\f||^2\ci{L^2(\mu)}}$.

Then we get
$$
\sum_{Q\in\mathcal{F}} \frac{|\langle \Delta\ci Q\f, \psi\rangle|^2}
{||\Delta\ci Q\f||^2\ci{L^2(\mu)}}
\le \sqrt2 \|\psi\|\ci{L^2(\mu)}
\Bigl[\sum_{Q\in\mathcal{F}}
\frac{|\langle \Delta\ci Q\f, \psi\rangle|^2}
{||\Delta\ci Q\f||^2\ci{L^2(\mu)} }\Bigr]^{\frac12},
$$
or, which is the same,
$$
\sum_{Q\in\mathcal{F}} \frac{|\langle \Delta\ci Q\f, \psi\rangle|^2}
{||\Delta\ci Q\f||^2\ci{L^2(\mu)}}
\le 2\|\psi\|^2\ci{L^2(\mu)}.
$$
Now, of course, the summation on the left can be extended to all
squares $Q$ for which $\|\Dl\ci Q\f\|\ci{L^2(\mu)}> 0$.

We will not need anything beyond this, so we are not going to say the magic
words that the projections $\Lambda$ and $\{\Dl\ci Q\}\ci{Q\in\D}$ generate
a Riesz basis of subspaces in $L^2(\mu)$ to a reader who does not want to
hear them.

\medskip

 Let $\D_1$ and $\D_2$ be two independent copies of the random
 dyadic lattice $\D$. Suppose that
there is some rule which allows one
to tell, for every square $Q_1\in\D_1$, whether it is  ``bad'' or
``good'' with respect to the lattice $\D_2$. Of course,
since $\D_1$ and $\D_2$ are  copies of {\it the same} random
dyadic lattice, we can
use the same rule to define bad squares in $\D_2$ with respect to $\D_1$.

Our next assumption is that bad squares are very rare.
Namely, we suppose that for every fixed $\D_1$ and for every $Q_1\in\D_1$,
the probability
$$
P\ci{\D_2}\{Q_1 \text{ is bad}\}\le \delta
$$
(and vice versa, of course).

If all the above assumptions are satisfied, we will say that $\D$
is a {\it perfect random dyadic lattice}.

Let again $\D_1$ and $\D_2$ be two independent copies of a random
dyadic lattice $\D$.

A function $\f_1\in L^2(\mu)$ is called
good (the full name should be $\D_1$-good
with respect to the lattice $\D_2$, or something like that)
if for every bad square $Q_1\in\D_1$, we have
$$
\Delta\ci{Q_1}\f_1=0.
$$
Even if a function $\f_1\in L^2(\mu)$
is not good, we still can write the decomposition
$$
\f_1=
\bigl[\,\Lambda_1\f_1+\sum\Sb Q_1\in\D_1, \\ Q_1\text{ is good} \endSb
\Delta\ci{Q_1}\f_1\,\bigr]
+
\sum\Sb Q_1\in\D_1, \\ Q_1\text{ is bad} \endSb\Delta\ci{Q_1}\f_1
=: (\f_1)_{\text{good}}+(\f_1)_{\text{bad}}.
$$
Note that this decomposition depends on both $\D_1$ and $\D_2$,
and therefore
$(\f_1)_{\text{good}}$ and $(\f_1)_{\text{bad}}$ are random functions
even if $\f_1$ is a sure function.
If the dyadic lattice $\D$ is perfect, it is easy to show  that
always
$$
\|(\f_1)_{\text{good}}\|\ci{L^2(\mu)},
\|(\f_1)_{\text{bad}}\|\ci{L^2(\mu)}
\le
2
\|\f_1\|\ci{L^2(\mu)}.
$$
What is more, if $\f_1$ does not depend on $\D_2$, then
for every fixed $\D_1$,
$$
\E\ci{\D_2}
\|(\f_1)_{\text{bad}}\|\ci{L^2(\mu)}^2\le 4 \delta
\|\f_1\|\ci{L^2(\mu)}^2.
$$
Indeed, we have
$$
\E\ci{\D_2}
\|(\f_1)_{\text{bad}}\|\ci{L^2(\mu)}^2=
\E\ci{\D_2}
\Bigl\|\sum\Sb Q_1\in\D_1, \\ Q_1\text{ is bad}\endSb
\Dl\ci{Q_1}\f_1
\Bigr\|\ci{L^2(\mu)}^2\le
2
\E\ci{\D_2}
\sum\Sb Q_1\in\D_1, \\ Q_1\text{ is bad}\endSb
\|\Dl\ci{Q_1}\f_1\|\ci{L^2(\mu)}^2
$$
$$
=
2\sum_{Q_1\in\D_1}P\ci{\D_2}\{Q_1\text{ is bad}\}
\|\Dl\ci{Q_1}\f_1\|\ci{L^2(\mu)}^2\le
2\delta
\sum_{Q_1\in\D_1}
\|\Dl\ci{Q_1}\f_1\|\ci{L^2(\mu)}^2\le
4\delta
\|\f_1\|\ci{L^2(\mu)}^2.
$$

Hence for all sure functions $\f_1$, we have
$$
\E
\|(\f_1)_{\text{bad}}\|\ci{L^2(\mu)}^2\le 4\delta
\|\f_1\|\ci{L^2(\mu)}^2.
$$

\bigskip

{\bf III. Perfect hair}

\bigskip

Let again $\mu$ be a measure   on the complex plane $\C$
satisfying $0<\mu(\C)<+\infty$.

Assume that we have a perfect random dyadic lattice
$\D$ (i.e., a family of dyadic lattices
endowed with some probability so that the assumptions of the previous section
are satisfied) and suppose that with every dyadic lattice $\D$
in that family a
nonnegative Lipschitz function $\Phi\ci\D$ is associated in such a way
that the following properties hold:
\medskip
\it
1) $\mu\{x\in\C\,:\,\Phi\ci\D(x)>0\}\le \delta\mu(\C)$\ \ for every $\D$;
\smallskip
2) For every two dyadic lattices $\D_1,\D_2$, for every Lipschitz
function $\Theta$ satisfying $\inf_{\C}\Theta>0$,
$\Theta\ge\delta\max(\Phi\ci{\D_1},\Phi\ci{\D_2})$,  and for
any two good functions $\f_1$ and $\f_2$ (with respect to the lattices
$\D_1$ and $\D_2$, correspondingly), we have
$$
\bigl|\langle \f_1,K\ci\Theta\f_2\rangle\bigl|\le
N
\|\f_1\|\ci{L^2(\mu)}\|\f_2\|\ci{L^2(\mu)},
$$
where $N$ is some (large) positive constant, not depending on $\f_1$, $\f_2$
or $\Theta$.

\smallskip\rm

(The assumption $\inf_{\C}\Theta>0$ is purely technical, of course: it just
allows us to avoid tiresome discussions concerning the definition
of $K\ci\Theta\f_2$: the kernel is uniformly bounded,
the measure is finite, so everything
makes sense.)
\medskip
Then we will say that we have  {\it``perfect hair"}.

Our first aim is to show that every perfect hair generates a bounded
(in $L^2(\mu)$) operator, which coincides with the Cauchy integral operator
everywhere outside an exceptional set of small $\mu$-measure.

\bigskip

{\bf IV. Truncated mathematical expectation}

\bigskip

Let $\xi$ be a nonnegative random variable and let $0<\beta<1$.
Define
$$
\E_\beta\xi:=\inf\Bigl\{
\int_{\Omega\setminus\Omega_1}\xi\, dP\,:\, P\{\Omega_1\}\le\beta
\Bigr\}
$$
Note that
\medskip\it
A) If $P\{\xi>0\}\le\beta$, then $\E_\beta\xi=0$;
\smallskip
B) $P\{\xi\ge\beta^{-1}\E_\beta\xi\}\le 2\beta$;
\smallskip
C) If $\Phi_\omega(x)$ ($x\in\C$) is a random nonnegative
Lipschitz function, then
$\E_\beta\Phi_\omega(x)$ is a certain nonnegative Lipschitz
function.
\medskip\rm

\bigskip

{\bf V. How to use perfect hair}

\bigskip

{\bf Theorem:}

Assume that we have a perfect hair. Let $\beta=\sqrt\delta$.
Let $\Phi:=\E_\beta\Phi\ci\D$.

Then

1) $\mu\{x\in\C\,:\,\Phi(x)>0\}\le\sqrt\delta\mu(\C)$;
\smallskip
2) The operator $K\ci\Phi$ acts in $L^2(\mu)$ in the sense that
$\sup_{\la>0}||K_{\Phi+\la}||\ci{L^2(\mu)\to L^2(\mu)}<+\infty$.

%
%
%

{\bf Proof:}

The first claim is easy: note that
$$
\E\mu\{x\in\C\,:\, \Phi\ci\D(x)>0\}\le\delta\mu(\C),
$$
and thereby for the set
$$
E:=
\bigl\{x\in\C\,:\,
P\{
 \Phi\ci\D(x)>0
\}
\ge \beta=\sqrt\delta
\bigr\},
$$
we have $\mu(E)\le \sqrt\delta\mu(\C)$.
It remains only to recall that, according to property (A) of the truncated
mathematical expectation,
$\Phi=\E_\beta\Phi\ci\D\equiv 0$ outside $E$.

\medskip

Now we will prove even a little bit more than the second claim.
Namely, we will show that
$$
\sup
\{||K\ci{\Theta}||\,:\,\Theta\text{ is Lipschitz }, \Theta\ge\Phi
\}<+\infty
$$
(in the same sense as above; see the exact formulation below).

Fix $\la>0$ and let
$$
N_\la=
\sup
\{||K\ci{\Theta}||\,:\,\Theta\text{ is Lipschitz }, \Theta\ge\Phi+\la
\}.
$$
Clearly, for every $\la>0$, we have
$N_\la\le \frac{\mu(\C)}{\la}<+\infty$. We are going to prove
that $N_\la$ is bounded by some constant {\it independent} of $\la$.

{\bf ``Space" and ``frequency'' reductions}

Choose $\Theta\ge\Phi+\la$ and functions
$\f_1,\f_2\in L^2(\mu)$ with $||\f_1||\ci{L^2(\mu)}=
||\f_2||\ci{L^2(\mu)}=1$ such that
$$
\bigl|\langle\f_1, K\ci\Theta\f_2\rangle\bigr|
\ge \tfrac{9}{10} N_\la.
$$
Consider two independent copies $\D_1$ and $\D_2$ of a
perfect random dyadic lattice $\D$.
Let
$$
S:=\{x\in\C\,:\,\max_{j=1,2}\Phi\ci{\D_j}(x)\ge\beta^{-1}\Phi(x)\} .
$$
Put
$$
\f'_j:=\f_j\chi\ci{S},\qquad \wt\f_j:=\f_j\chi\ci{\C \setminus S}
=\f_j-\f'_j\qquad j=1,2
$$
(``space'' reduction) and, at last,
$$
\psi_j:=(\wt\f_j)_{\text{good}}
=\wt\f_j-(\wt\f_j)_{\text{bad}}
, \qquad j=1,2
$$
(``frequency'' reduction).

We expect both reductions to be just ``minor corrections''. Soon we will
show that this really is the case, namely, that
$$
||\f'_j||\ci{L^2(\mu)},
||(\wt\f_j)_{\text{bad}} ||\ci{L^2(\mu)}\le\frac{1}{10}\qquad\qquad (*)
$$
with probability close to $1$. Now let us demonstrate that these reductions
really make sense.

Pick a pair of dyadic lattices $\D_1$ and $\D_2$, for which $(*)$ holds.
We have
$$
\langle
\wt\f_1, K\ci\Theta\wt\f_2
\rangle              =
\langle
\f_1, K\ci\Theta\f_2
\rangle-
\langle
\f'_1, K\ci\Theta\f_2
\rangle-
\langle
\wt\f_1, K\ci\Theta\f'_2
\rangle
$$
and thereby
$$
\bigl|\langle
\wt\f_1, K\ci\Theta\wt\f_2
\rangle\bigr|\ge \tfrac{9}{10}N_\la-\tfrac{2}{10}\|K\ci\Theta\|\ge
\tfrac{7}{10}N_\la
$$
(here we used the obvious estimate $||\wt\f_1||\ci{L^2(\mu)}\le
||\f_1 ||\ci{L^2(\mu)}=1$ together with $(*)$ to get the first
inequality).

The key observation about the space reduction is that
$$
\langle\wt\f_1, K\ci\Theta\wt\f_2\rangle=\iint k\ci\Theta(x_1,x_2)
\wt\f_1(x_1)\wt\f_2(x_2)\,d\mu(x_1)d\mu(x_2)=
$$
$$
\iint k\ci{\Theta'}(x_1,x_2)
\wt\f_1(x_1)\wt\f_2(x_2)\,d\mu(x_1)d\mu(x_2)=
\langle\wt\f_1, K\ci{\Theta'}\wt\f_2\rangle,
$$
where
$$
\Theta':=\max\{\Theta, \beta\Phi\ci{\D_1}, \beta\Phi\ci{\D_2} \}.
$$
We still have $\Theta'$ Lipschitz and satisfying $\Theta'\ge\Phi+\la$, but
now also $\Theta'\ge\delta\max\{\Phi\ci{\D_1},\Phi\ci{\D_2}\}$, and therefore
we only need to make the functions $\wt\f_j$ good to apply property $(2)$ of
perfect hair and to finish the story. This is exactly what the frequency
reduction does.
Like above, we can write
$$
\langle
\psi_1, K\ci{\Theta'}\psi_2
\rangle              =
\langle
\wt\f_1, K\ci{\Theta'}\wt\f_2
\rangle-
\langle
(\wt\f_1)_{\text{bad}}, K\ci{\Theta'}\wt\f_2
\rangle-
\langle
\psi_1, K\ci{\Theta'}(\wt\f_2)_{\text{bad}}
\rangle
$$
and thereby
$$
\bigl|\langle
\psi_1, K\ci{\Theta'}\psi_2
\rangle\bigr|\ge \tfrac{7}{10}N_\la-\tfrac{3}{10}\|K\ci{\Theta'}\|\ge
\tfrac{4}{10}N_\la
$$
(here we used the estimate $||\psi_j||\ci{L^2(\mu)}\le
||\wt\f_j ||\ci{L^2(\mu)}+||(\wt\f_j)_{\text{bad}}
||\ci{L^2(\mu)}<2$ together with $(*)$ to get the first
inequality).

Now, according to property $(2)$ of perfect hair, the left hand part does not
exceed $4N$ and we get $N_\la\le 10N$.
It remains only to prove that $(*)$ holds with probability close to $1$.

Note that for any given point $x\in\C$, we have
$
P\{x\in S\}\le 4\beta
$, and therefore,
$$
\E||\f'_j||\ci{L^2(\mu)}^2\le 4\beta,\qquad j=1,2.
$$
Hence,
$$
P\bigl\{||\f'_j||\ci{L^2(\mu)}\ge \beta^{\frac13}\bigr\}\le
4\beta^{\frac13}\qquad j=1,2.
$$
Now we would like to say that the norms of the functions
$(\wt\f_j)_{\text{bad}}$ are small as well.
Unfortunately, as constructed, each of them depends on both
$\D_1$ and $\D_2$.
So it seems that we can only apply the obvious estimate
$
||(\wt\f_j)_{\text{bad}} ||\ci{L^2(\mu)}\le
2||(\wt\f_j)||\ci{L^2(\mu)}\le 2
$,
which is clearly useless.

Note, nevertheless, that
$$
(\wt\f_j)_{\text{bad}}=(\f_j)_{\text{bad}}-(\f'_j)_{\text{bad}}.
$$
The norm of  $(\f'_j)_{\text{bad}}$ does not exceed
$2\|\f'_j\|\ci{L^2(\mu)}$.
As to
$(\f_j)_{\text{bad}}$,
we can apply the estimate for sure functions to $\f_j$,
which yields
$$
\E||(\f_j)_{\text{bad}}||^2\ci{L^2(\mu)}\le 4\delta.
$$
So finally we conclude that with probability at least
$
1-8\beta^{\frac13}-8\delta^{\frac13}>\frac{9}{10}
$
all the norms in the left hand part of $(*)$
are bounded by $2\beta^{\frac13}+\delta^{\frac13}<\frac{1}{10}$.

\bigskip

{\bf VI. Lyric deviation: Hausdorff measure and
 analytic capacity}

\bigskip

We will start with a couple of definitions.

{\bf The 1-dimensional Hausdorff measure}

Let $\e>0$. For every set $E\subset\C$ define
$$
\mathcal{H}_\e(E):=\inf\Bigl\{
\sum_j r_j\,:\,E\subset\bigcup_j B(x_j,r_j),\quad x_j\in\C,\, r_j\le\e
\Bigr\}
$$
(the infimum is taken over all (countable) coverings of $E$ by open disks
$B(x_j,r_j)$ with radii $r_j\le\e$).

It is clear that $\mathcal{H}_\e$ is an outer measure and that if $\e'\le
\e''$, then $\mathcal{H}_{\e'}(E)\ge\mathcal{H}_{\e''}(E)$ for every $E\subset\C$. Since every
monotone function has a limit (maybe, infinite), we can define
$$
\mathcal{H}(E):=\lim_{\e\to 0}\mathcal{H}_\e(E)=\sup_{\e> 0}\mathcal{H}_\e(E).
$$
It is a trivial exercise to show that $\mathcal{H}$ is an outer measure. However, it
is much better than just that, namely, $\mathcal{H}$ is {\it a Borel measure}. The
proof of this remarkable theorem can be found in any (good) textbook on
measure theory. We can only regret that it is not included in the Leningrad
(or Michigan)
State University analysis course.

{\bf Analytic capacity}

Let $F\subset\C$ be a compact set. We will say that $F$ has positive analytic
capacity if there exists a bounded analytic function
$f:\C\setminus F\to \C$, which is not identically $0$ and such that
$f(x)\to 0$ as $x\to\infty$.
\medskip
Assume now that we have a compact set $F$ of positive analytic capacity
and such that $\mathcal{H}(F)<+\infty$. Let $f$ be the corresponding bounded analytic
function.

\bigskip

{\bf VII. Cauchy integral representation}

\bigskip

We devote this section to a well-known representation of bounded analytic function 
outside of a compact of finite length.  
Since $F$ is compact, we can consider only finite coverings in the definition
of $\mathcal{H}_\e(F)$. Now for every positive integer $n$, construct a covering
$$
\bigcup_{j=1}^{N(n)} B(x_j^{(n)},r_j^{(n)})\supset F
$$
such that all $r_j^{(n)}\le \frac1n$, $\sum_j r_j^{(n)}\le \mathcal{H}(F)+\frac1n$ and
$B(x_j^{(n)},r_j^{(n)})\cap F\ne\emptyset$ for every $j$.
Let $\Omega_n=\C\setminus\cl\bigcup_j B(x_j^{(n)},r_j^{(n)})$
and let $\Gamma_n:=\partial \Omega_n$.
Note that $\Gamma_n$ is a good contour (consisting of finitely many arcs) and
that $\Gamma_n\subset\C\setminus F$. Therefore for every point
$x\in\Omega_n$,
we can write the standard Cauchy formula
$$
f(x)=-\frac{1}{2\pi
i}\oint_{\Gamma_n}\frac{f(y)\,dy}{x-y}=\int_\C\frac{d\nu_n(y)}{x-y},
$$
where $\nu_n$ is a complex-valued measure defined (on Borel sets, say)
by
$$
\nu_n(E)=-\frac{1}{2\pi i}\oint_{\Gamma_n\cap E}f(y)\,dy.
$$
Note that the variations of the complex-valued measures $\nu_n$ are uniformly
bounded (by $||f||\ci{L^\infty}(\mathcal{H}(F)+1)$, say), and therefore (passing to a
subsequence, if needed) we may assume that $\nu_n$ weakly converge to a
complex-valued Borel measure $\nu$ (over the space
$C_0(\C)$  of all
compactly supported  complex-valued
continuous functions on $\C$).
Note now that $\Omega_n$ contains all points $x\in\C$ for which
$\dist(x,F)>\frac1n$.
Therefore for any
 $\eta\in C_0(\C)$ satisfying $\supp\f\subset\C\setminus F$,
we have
$$
\int_\C\eta\,d\nu=\lim_{n\to\infty}\int_\C\eta\,d\nu_n=0,
$$
i.e., $\supp\nu\subset F$.
Now passing to the limit in the Cauchy formula above, we see that for every
$x\in\C\setminus F$,
$$
f(x)=\int_\C\frac{d\nu(y)}{x-y}.
$$
Our next step will be to show that for every Borel set $E\subset\C$ we have
$$
|\nu|(E)\le \|f\|\ci{L^\infty} \mathcal{H}F(E),
$$
where
$\mathcal{H}F(E):=\mathcal{H}(E\cap F)$.

Recall that every finite Borel measure $\mu$ is regular in the sense that
for every Borel set $E$ one can find an open set $G\supset E$ such that
$\mu(G\setminus E)$ is as small as one wants. Therefore it is enough to
prove this inequality for open sets only.

Recall also that for an open set $G$,
$$
|\nu|(G)=\sup\Bigl\{
\Bigl|\int_\C\eta\,d\nu\Bigr|\,:\, \eta\in C_0(\C), \supp\eta\subset G,
||\eta||\ci{L^\infty}\le 1
\Bigr\}.
$$
Therefore we need only to show that for every such $\eta$,
$$
\Bigl|\int_\C\eta\,d\nu\Bigr|\le\|f\|\ci{L^\infty}\mathcal{H}F(G).
$$
But we have
$$
\int_\C\eta\,d\nu=\lim_{n\to\infty}\int_\C\eta\,d\nu_n,
$$
and for every $n$,
$$
\Bigl|\int_\C\eta\,d\nu\Bigr|\le\|f\|\ci{L^\infty}
\sum_{j:B(x_j^{(n)},r_j^{(n)})
\cap \supp\eta\ne\emptyset} r_j^{(n)}.
$$
Now notice that if $\frac1n<\dist(\supp\eta, \partial G)$, then the
disks $B(x_j^{(n)},r_j^{(n)})$ that intersect $\supp\eta$ cannot participate
in the covering of $F\setminus G$. Therefore,
$$
\sum_{j:B(x_j^{(n)},r_j^{(n)})   
\cap \supp\eta\ne\emptyset} r_j^{(n)}\le
\sum_{j=1}^{N(n)} r_j^{(n)}
-
\sum_{j:B(x_j^{(n)},r_j^{(n)})
\cap F\setminus G\ne\emptyset} r_j^{(n)}
$$
$$
\le \mathcal{H}(F)+\frac1n-\mathcal{H}_{\frac1n}(F\setminus G)\to \mathcal{H}(F)-\mathcal{H}(F\setminus G)=
\mathcal{H}(F\cap G)=\mathcal{H}F(G)
$$
as $n\to\infty$, proving the claim.

Applying the Radon-Nykodim theorem, we conclude that there exists a Borel
measurable function $h$ satisfying $\|h\|\ci{L^\infty}\le\|f\|\ci{L^\infty}$ and
such that
$$
f(x)=\int_{\C}\frac{h(y)}{x-y}d\mathcal{H}F(y)
$$
for every $\C\setminus F$ (there is no problem with convergence here,
because, as we remember, the integral is actually taken over $F$).
Note also that since $f(x)\ne 0$ for at least one $x\in \C$, we should have
$\mathcal{H}F\{y\in \C\,:\,h(y)\ne 0\}>0$. As a trivial consequence, we observe that
$\mathcal{H}(F)>0$.

\bigskip

{\bf VIII. The Ahlfors radius $\mathcal{R}(x)$}

\bigskip

Now take some large $M>1$. We will call a disk $B(x,r)$
($x\in \C, r>0$) {\it non-Ahlfors}, if
$$
\mathcal{H}F(B(x,r))>Mr.
$$
For every point $x\in\C$ define its Ahlfors radius $\mathcal{R}(x)$ by
$$
\mathcal{R}(x):=\sup\{r>0\,:\, B(x,r)\text{ is non-Ahlfors}\,
\}.
$$
Since $f$ is bounded on $\C\setminus F$, so is the Cauchy integral
$\int_{\C}\frac{h(y)}{x-y}d\mathcal{H}F(y)$. We are going to show that, in a sense,
this integral stays bounded on $F$ as well (where $f$, generally speaking,
 does not exist). Namely, for every $x\in\C$
$$
\sup_{\e>\mathcal{R}(x)}\Bigl|
\int_{\C\setminus B(x,\e)}\frac{h(y)}{x-y}d\mathcal{H}F(y)
\Bigr|\le 7M\|f\|\ci{L^\infty}.
$$

{\bf Proof:}

Note first of all, that the condition $\mathcal{H}(F)<+\infty$ implies that the
2-dimensional Lebesgue measure $m(F)=0$. Indeed, for every covering
$\bigcup_jB(x_j,r_j)\supset F$, we have
$$
m(F)\le \pi\sum_j r_j^2\le \pi(\max_j r_j)\sum_j r_j.
$$
Therefore
$$
m(F)\le \pi\e\mathcal{H}_\e(F)\le \pi\e\mathcal{H}(F)
$$
for every $\e>0$, and we are done.

Now compare
$
\int_{\C\setminus B(x,\e)}\frac{h(y)}{x-y}d\mathcal{H}F(y)
$
to
$$
\frac{4}{\pi\e^2}
\int_{B(x,\frac\e2)}f(z)\,dm(z)=
\frac{4}{\pi\e^2}
\int_{B(x,\frac\e2)\setminus F}f(z)\,dm(z),
$$
which is clearly bounded by $\|f\|\ci{L^\infty}\le M\|f\|\ci{L^\infty}$.
Using the Cauchy integral representation for $f$, we see that the difference
equals
$$
-\frac{4}{\pi\e^2}
\int_{B(x,\frac\e2)}\Bigl(
\int_{B(x,\e)}\frac{h(y)}{z-y}d\mathcal{H}F(y)
\Bigr)dm(z)+
$$
$$
\int_{\C\setminus B(x,\e)}h(y)
\Bigl(\frac{1}{x-y}
-
\frac{4}{\pi\e^2}
\int_{B(x,\frac\e2)}\frac{dm(z)}{z-y}\Bigr)d\mathcal{H}F(y)=: I'+I''.
$$
The integral $I'$ allows the rough estimate
$$
|I'|\le
\frac{4}{\pi\e^2}\|f\|\ci{L^\infty}
\int_{B(x,\e)}\Bigl(
\int_{B(x,\frac\e2)}\frac{dm(z)}{|z-y|}\Bigr)d\mathcal{H}F(y).
$$
Since the inner integral does not exceed $\pi\e$ for every $y\in\C$, we get
$$
|I'|\le 4\e^{-1}\|f\|\ci{L^\infty}\mathcal{H}F(B(x,\e))\le 4M\|f\|\ci{L^\infty},
$$
provided that $\e\ge \mathcal{R}(x)$.

To estimate $I''$, note that 
$$
\Bigl|\frac{1}{x-y}
-
\frac{4}{\pi\e^2}
\int_{B(x,\frac\e2)}\frac{dm(z)}{z-y}\Bigr|
=
\frac{4}{\pi\e^2}
\Bigl|
\int_{B(x,\frac\e2)}\frac{z-x}{(x-y)(z-y)}dm(z)\Bigr|\le\frac{\e}{|x-y|^2}
$$
because $|z-x|\le\frac{\e}{2}$ and $2|z-y|\ge |x-y|$
for every $y\in\C\setminus B(x,\e)$,
$z\in B(x,\frac{\e}{2})$.

Thus
$$
|I''|\le \|f\|\ci{L^\infty}\int_{\C\setminus
B(x,\e)}\frac{\e}{|x-y|^2}d\mathcal{H}F(y).
$$
To estimate the last integral, we need the following obvious
lemma, which we will
frequently use in the future.

{\bf Comparison Lemma:}
Let $S\subset \C$. Assume that we have a measure $\mu$ satisfying
$$
\mu\{x\in\C\,:\, \dist(x,S)<r\}\le Mr \text{\qquad for every }r\ge R_0
$$
and a nonnegative continuous decreasing function $U(t)$ ($t>0$).

Then for every $R\ge R_0$
$$
\int_{\{x:\dist(x,s)\ge R\}} U(\dist(y,S))d\mu(y)\le
M\Bigl(RU(R)+\int_R^{+\infty}U(t)dt\Bigr).
$$
Note also that the quantity in parentheses can be viewed as the
integral over the whole ray $[0,\infty)$ of $\min\{U(t),U(R)\}$ and therefore
is a decreasing function in $R$. So, we can replace $R$ on the right hand
side by any lesser number if we want to.

\medskip
The Comparison Lemma (with $S=\{x\}$, $R_0=\mathcal{R}(x)$, $R=\e$ and
$U(t)=\frac{\e}{t^2}$) yields
$$
\int_{\C\setminus
B(x,\e)}\frac{\e}{|x-y|^2}d\mathcal{H}F(y)\le
M\Bigl(1+\int_\e^{+\infty}\frac{\e}{t^2}dt\Bigr)=2M,
$$
and thereby $|I''|\le 2M\|f\|\ci{L^\infty}$.
It remains only to add the estimates to get the desired inequality.

The additional assumption $\e\ge \mathcal{R}(x)$ in the formulation of the last
statement seems quite unpleasant. We would prefer to have a result
that is valid for
every $\e>0$. This can be achieved if we replace the kernel $\frac{1}{x-y}$
by the suppressed kernel $k\ci\Phi$ with a Lipschitz function $\Phi$
satisfying $\Phi(x)\ge\delta \mathcal{R}(x)$ for every $x\in\C$.

{\bf Lemma:}

For every $x\in\C$
$$
\sup_{\e>0}\Bigl|
\int_{\C\setminus B(x,\e)}k\ci\Phi(x,y)h(y)d\mathcal{H}F(y)
\Bigr|\le  (11+\delta^{-1})M\|f\|\ci{L^\infty}.
$$

{\bf Proof:}
Recall that the kernel $k\ci\Phi$ is given by
$$
k\ci\Phi(x,y)=\frac{\overline{x-y}}{|x-y|^2+\Phi(x)\Phi(y)}.
$$
Put
$r:=\Phi(x)$, $R:=\delta^{-1}\Phi(x)\,(\ge \mathcal{R}(x))$ and, at last,
$R':=\max\{\e,R\}$.
Write
$$                                               
\int_{\C\setminus B(x,\e)}k\ci\Phi(x,y)h(y)d\mathcal{H}F(y)=
\int_{\C\setminus B(x,R')}\dots
+\int_{B(x,R)\setminus B(x,\e)}\dots =: I'+I''
$$
($R$ in the second integral is not a misprint: we need this second term
only for $R>\e$ when $R'=R$).

Recall that $|k\ci\Phi(x,y)|\le \frac{1}{\Phi(x)}= r^{-1}$
for all $y\in\C$.

Thus
$$
|I''|\le\|f\|\ci{L^\infty}r^{-1}\mathcal{H}F(B(x,R))\le \|f\|\ci{L^\infty}r^{-1}MR
= \delta^{-1}M \|f\|\ci{L^\infty}.
$$
As to $I'$, let us compare it to $\wt I:=\int_{\C\setminus
B(x,R')}\frac{h(y)}{x-y}d\mathcal{H}F(y)$, which is bounded by $7M\|f\|\ci{L^\infty}$,
because $R'\ge \mathcal{R}(x)$.
The difference does not exceed
$$
\|f\|\ci{L^\infty}\int_{\C\setminus B(x,R')}
\Bigl|
     \frac{1}{x-y}-\frac{\overline{x-y}}{|x-y|^2+\Phi(x)\Phi(y)}
\Bigr| d\mathcal{H}F(y).
$$
Representing $\frac{1}{x-y}$ as $\frac{\overline{x-y}}{|x-y|^2}$ and
observing that for every two numbers $t,s>0$, one has
$$
\frac{1}{t}-\frac{1}{t+s}=\frac{s}{t(t+s)}\le \frac{s}{t^2},
$$
we get
$$
\Bigl|
     \frac{1}{x-y}-\frac{\overline{x-y}}{|x-y|^2+\Phi(x)\Phi(y)}
\Bigr|\le
\frac{\Phi(x)\Phi(y)}{|x-y|^3}\le \frac{r(r+|x-y|)}{|x-y|^3}.
$$
Applying the Comparison Lemma again, we see that
$$
|I'-\wt I|\le M\|f\|\ci{L^\infty}\Bigl[
r\frac{r(r+r)}{r^3}+\int_r^{+\infty}
\frac{r(r+t)}{t^3}dt\Bigr]=
\tfrac72
M\|f\|\ci{L^\infty}\le 4
M\|f\|\ci{L^\infty}
$$
(here, in order to simplify calculations, we used the possibility
to replace $R'$ by the lesser number $r$).

Now it remains only to bring all the estimates together to get the conclusion
of the lemma.

\bigskip

{\bf IX. The exceptional set $H$}

\bigskip

The demand $\Phi(x)\ge\delta \mathcal{R}(x)$ is much less restrictive than it seems at
first glance. Let us show that if $M$ is sufficiently large, then
$\mathcal{R}(x)=0$ for most $x$. Indeed, for every non-Ahlfors point $x\in \C$, one can
find a disk $B(x,r)$ such that $\mathcal{H}F(B(x,r))>Mr$. Using the Vitali covering
theorem, we can construct a countable family of {\it pairwise disjoint}
non-Ahlfors disks $B(x_j,r_j)$ such that every non-Ahlfors disk $B(x,r)$ is
contained in the union
$$
H:=\bigcup_j B(x_j,5r_j).
$$
Note that $r_j<\frac{\mathcal{H}F(B(x_j,r_j))}{M}$ and therefore
$$
\sum_{j} r_j<\frac{\mathcal{H}F(\C)}{M}.
$$
Observing that every term in the sum is not greater than the whole sum,
we conclude that
$$
\mathcal{H}\cci{\frac{5\mathcal{H}F(C)}{M}}(H)\le\frac{5\mathcal{H}F(\C)}{M},
$$
and thereby,
$$
\mathcal{H}F(\C\setminus H)=\mathcal{H}(F\setminus H)\ge
\mathcal{H}\cci{\frac{5\mathcal{H}F(\C)}{M}}(F\setminus H)\ge
\mathcal{H}\cci{\frac{5\mathcal{H}F(\C)}{M}}(F)-
\mathcal{H}\cci{\frac{5\mathcal{H}F(\C)}{M}}(H)\geq
$$
$$
\mathcal{H}\cci{\frac{5\mathcal{H}F(\C)}{M}}(F)-
\frac{5\mathcal{H}F(\C)}{M}\to \mathcal{H}(F)=\mathcal{H}F(\C)
$$
as $M\to +\infty$.
Thus $\mathcal{H}F(H)=\mathcal{H}F(\C)-\mathcal{H}F(\C\setminus H)\to 0$ as $M\to+\infty$, proving the
claim.

Now define
$$
\wt{\mathcal{R}}(x):=\dist(x,\C\setminus H).
$$
Clearly $\wt{\mathcal{R}}$ is a Lipschitz function. Since every non-Ahlfors disk is
contained in $H$, we have $\wt{\mathcal{R}}(x)\ge \mathcal{R}(x)$. At last
$\mathcal{H}F\{x\in\C\,:\,\wt{\mathcal{R}}(x)>0\}=\mathcal{H}F(H)$ can be made as small as one wants by
choosing the constant $M$ large enough.

\bigskip

{\bf X. Localization}

\bigskip

Let $x_0$ be any $L^2$-Lebesgue point of $h$ with respect to the measure
$\mathcal{H}F$ satisfying $h(x_0)\ne 0$. Recall that it means
$$
\mathcal{H}F(B(x_0,r))>0\qquad \text{ for every } r>0;
$$
$$
\frac{1}{\mathcal{H}F(B(x_0,r))}\int_{B(x_0,r)}|h(x)-h(x_0)|^2\,d \mathcal{H}F(x)\to 0
\quad\text{ as }r\to 0.
$$
Since the measure $\mathcal{H}F$ is finite, $\mathcal{H}F$-almost every point $x\in\C$ is a
Lebesgue point of a bounded function $h$
(actually this statement is true for {\it any}
$L^2(\mathcal{H}F)$-function
$h$).
On the other hand, as we have seen above,
$\mathcal{H}F\{x\in\C\,:\,h(x)\ne 0\}>0$. So, the needed point $x_0$ really exists.

Now choose $0<\rho<\frac{1}{8}$ so small that
$$
\frac{1}{\mathcal{H}F(B(x_0,\rho))}\int_{B(x_0,\rho)}|h(x)-h(x_0)|^2\,d \mathcal{H}F(x)
<\delta^4|h(x_0)|^2.
$$
Choose $M>1$ so large that for the corresponding exceptional set $H$, we have
$$
\mathcal{H}F(H)\le\frac{\delta}{3}\mathcal{H}F(B(x_0,\rho)).
$$
Now let $\rho'<\rho$ be so close to $\rho$ that
$$
\mathcal{H}F(B(x_0,\rho)\setminus B(x_0,\rho'))<
\frac{\delta}{3}\mathcal{H}F(B(x_0,\rho)).
$$
Let
$$
\wt\Phi(x):=\max\{\wt{\mathcal{R}}(x), |x-x_0|-\rho'\}.
$$
Note that $\wt\Phi(x)$ is a nonnegative Lipschitz function majorizing
the Ahlfors radius $\mathcal{R}(x)$ and that
$$
\mathcal{H}F\{x\in B(x_0,\rho)\,:\, \wt\Phi(x)>0  \}\le\frac{2\delta}{3}\mathcal{H}F(B(x_0,\rho)).
$$
Define  the Borel measure $\mu$ by
$$
\mu(E):=\mathcal{H}F(E\cap B(x_0,\rho)).
$$
Note that for every Lipschitz function $\Theta\ge\delta\wt\Phi$ and for every
$x\in\C$ we have
$$
K\ci\Theta^{\sharp}h(x):=
\sup_{\e>0}
\Bigl|
\int_{\C\setminus B(x,\e)}k\ci\Theta(x,y)h(y)\,d\mu(y)
\Bigr|\le
$$
$$
\bigl[
 (11+\delta^{-1})M+\delta^{-1}(\rho-\rho')^{-1}\mathcal{H}F(\C)
\bigr]\|f\|\ci{L^\infty}=: B|h(x_0)|.
$$
Indeed, if we replace $d\mu(y)$ by $d\mathcal{H}F(y)$, we will have the bound
$(11+\delta^{-1})M\|f\|\ci{L^\infty}$ for the supremum.
The difference between the corresponding integrals does not exceed
$$
\int_{\C\setminus B(x_0,\rho)}|k\ci\Theta(x,y)|\cdot|h(y)|\,d\mathcal{H}F(y)\le
\|f\|\ci{L^\infty}
\int_{\C\setminus B(x_0,\rho)}\frac{d\mathcal{H}F(y)}{\Theta(y)}\le
\|f\|\ci{L^\infty}
\frac{\mathcal{H}F(\C)}{\delta(\rho-\rho')},
$$
and we are done.

\medskip

Now it is time to bring all the information together. Having started with a
compact set $F$ of finite Hausdorff measure and positive analytic capacity,
we have constructed a bounded Borel measurable function $h$, a point
$x_0\in\C$ for which $h(x_0)\ne 0$, a measure $\mu$ (which is just $\mathcal{H}F$
restricted to some small
disk centered at $x_0$),
a large constant $M$, an open set $H$, a Lipschitz
function $\wt\Phi$ and a (large) constant $B$ (they are listed in the order
one can choose them) such that the following properties hold:
\medskip\it

1) Every non-Ahlfors disk is contained in $H$, in particular,
$\mu(B(x,r))>Mr\Longrightarrow B(x,r)\subset H$
{\rm (recall that $\mu(B(x,r))\le\mathcal{H}F(B(x,r))$\,)}\/;
\smallskip

2) $h(x)=h(x_0)+g(x)$ with $\int_\C|g|^2d\mu\le\delta^4\mu(\C)$;
\smallskip

3) $\wt\Phi(x)\ge \dist(x,\C\setminus H),\qquad \mu\{x\in\C\,:\,\wt\Phi(x)>0
\}\le \frac{2\delta}{3}\mu(\C)$
\smallskip

4) For every Lipschitz function $\Theta\ge\delta\wt\Phi$ and for every
point $x\in\C$, one has $K\ci\Theta^{\sharp}h(x)\le B|h(x_0)|$.
\medskip\rm

We recommend the reader to reread this list of objects and their
properties attentively several times. They  are all completely natural,
but a little  too many to grasp at first glance.

\bigskip

{\bf XI. Construction of perfect hair}

\bigskip

Given $\delta,M,B,h,H$ and $\wt\Phi$ as above, let us construct perfect hair.
In order not to drag $x_0$ and $h(x_0)$ along all the time, assume that
$x_0=0$ and $h(x_0)=1$. Clearly, the problem is shift-invariant, and we
specially wrote all the above conditions in such a way that division
of $f$ and $h$ by the same constant would change  nothing in them.

First we should construct a perfect dyadic lattice $\D$. Our
construction will be surprisingly simple (compared to Guy David's
decomposition, say): we will just take the standard dyadic lattice and
consider its random shifts. 

Pick any point $\om\in[-\frac{1}{4},\frac14)^{2}$ and take the square
$Q^0(\om):=\om+[-\frac{1}{2},\frac12)^{2}$ as the ``starting'' square of the
dyadic lattice $\D=\D(\om)$. Recall that $\supp\mu\subset B(0,\frac{1}{8})$
and therefore $\supp\mu\subset Q^0(\om)$ for every such $\om$.

We are going to assign equal probability to every $\om$; so our probability
$P$ will be just $4$ times the Lebesgue measure restricted to
$[-\frac{1}{4},\frac14)^{2}$.

Once we have fixed the starting square $Q^0$, we have no choice of how to
{\it position} the smaller squares of $\D$: we just split $Q^0$ into four
equal subsquares (of the same kind $[a,b)\times[c,d)$\,), then split each new
square etc. Nevertheless, we still have the freedom of {\it how far down to
go} at every point. Now we are going to use this freedom.

We will call a square {\it terminal} in the following two cases:
\medskip\it
1) $Q\subset H$; or
\smallskip
2) $\int_Q|g|^2\,d\mu\ge\delta^2\mu(Q)$.
\medskip\rm
Note that in particular, $(2)$ holds if $\mu(Q)=0$.
If a square is not terminal, we will call it {\it transit}.

Now start the construction of $\D$ with the square $Q^0$, which is always
transit. It has size (side length) $l(Q^0)=2^{-0}=1$. Split it into four equal
subsquares. Some of them may be terminal and we will not touch those any
more. But we will further split each transit square of size $2^{-1}$ into
four subsquares of size $2^{-2}$ and so on.
 

\bigskip

{\bf XII. Projections $\Lambda$ and $\Dl_Q$}

\bigskip

Let $\D$ be one of the dyadic lattices constructed above. For a function
$\psi\in L^1(\mu)$ and for a square $Q\subset\C$, denote by
$\langle \psi\rangle\ci Q$ the average value of $\psi$ over $Q$ with respect
to the measure $\mu$, i.e.,
$$
\langle \psi\rangle\ci Q := \frac{1}{\mu(Q)}\int_Q\psi\,d\mu
$$
(of course, $\langle \psi\rangle\ci Q$ makes sense only for squares $Q$ with
$\mu(Q)>0$).

Put
$$
\Lambda\f :=\frac{\langle \f\rangle\ci {Q^0}}{\langle h\rangle\ci {Q^0}}h.
$$
Clearly, $\Lambda\f\in L^2(\mu)$ for all $\f\in L^2(\mu)$,
and $\Lambda^2=\Lambda$, i.e., $\Lambda$ is a projection. Note also, that actually $\Lambda$
does not depend on the lattice $\D$, because the average is taken over the
whole support of the measure $\mu$ regardless of the position of the square
$Q^0$.

From now on, we will always denote by $Q_j$ ($j=1,2,3,4$) the four subsquares
of a square $Q$ enumerated in some ``natural order''
 (to be chosen by the reader). In particular, that means that we will have to
give up our idea to denote the squares in
two copies $\D_1$ and $\D_2$ of the same random dyadic lattice $\D$
by $Q_1$ and $Q_2$, respectively.
This is okay, because while above it was important to emphasize the symmetry
between $\D_1$ and $\D_2$, below we will start almost every  claim with
``Assume (for definiteness) that $l(Q)\le l(R)$...''.

For every square $Q\in\D^{tr}$, define $\Dl\ci Q\f$ by
$$
\bigl.
\Dl\ci Q\f\bigr|\ci{\C\setminus Q}:=0,\qquad
\bigl.
\Dl\ci Q\f\bigr|\ci{Q_j}:=
\left\{
\aligned
\left[\frac{\langle \f\rangle\ci {Q_j}}{\langle h\rangle\ci {Q_j}}
 -\frac{\langle \f\rangle\ci {Q}}{\langle h\rangle\ci {Q}}\right]
h, &\text{\quad if $Q_j$ is transit;}
\\
\f-\frac{\langle \f\rangle\ci {Q}}{\langle h\rangle\ci {Q}}
h, &\text{\quad if $Q_j$ is terminal}
\endaligned\right.
$$
($j=1,2,3,4$).
Observe that for every transit square $Q$, we have $\mu(Q)>0$
and 
$$
\langle h\rangle\ci {Q}=1+\langle g\rangle\ci {Q}; \qquad\quad
|\langle g\rangle\ci {Q}|\le \sqrt{\langle |g|^2\rangle\ci {Q}}\le \delta,
$$
so our
definition makes sense: no zero can appear in the denominator.

{\bf Easy properties of $\Dl\ci Q\f$}

For every $\f\in L^2(\mu)$ and $Q\in\D^{tr}$,
\medskip\it
1) $\Dl\ci Q\f\in L^2(\mu)$ ;
\smallskip
2) $\int_\C\Dl\ci Q\f\,d\mu=0$;
\smallskip
3) $\Dl\ci Q$ is a projection, i.e., $\Dl\ci Q^2=\Dl\ci Q$;
\smallskip
4) $\Dl\ci Q\Lambda=\Lambda\Dl\ci Q=0$;
\smallskip
5) If $R\in\D^{tr}$ and $R\ne Q$, then $\Dl\ci Q\Dl\ci R=0$.
\medskip\rm
To check these properties is left to the reader as an exercise.

{\bf Lemma:}

For every $\f\in L^2(\mu)$ we have
$$
\f=\Lambda\f+\sum_{Q\in\D^{tr}}\Dl\ci Q\f,
$$
the series converges in $L^2(\mu)$ and, moreover,
$$
2^{-1}\|\f\|^2\ci{L^2(\mu)}\le
\|\Lambda\f\|^2\ci{L^2(\mu)}+\sum_{Q\in\D^{tr}}\|\Dl\ci Q\f\|^2\ci{L^2(\mu)}  \le
2\|\f\|^2\ci{L^2(\mu)}.
$$

{\bf Proof:}

Note first of all that if one understands the sum $\sum_{Q\in\D^{tr}}$
as $\lim_{n\to\infty}\sum_{Q\in\D^{tr}:l(Q)> 2^{-n}}$, then for
$\mu$-almost every $x\in\C$, one has
$$
\f(x)=\Lambda\f(x)+\sum_{Q\in\D^{tr}}\Dl\ci Q\f(x).
$$
Indeed, the claim is obvious if the point $x$ lies in some terminal square.
Suppose now that it is not the case. Observe that
$$
\Lambda\f(x)+\sum_{Q\in\D^{tr}:l(Q)> 2^{-n}}\Dl\ci Q\f(x)=
\frac{\langle \f\rangle\ci{Q^{n}} }{\langle h\rangle\ci{Q^{n}} }h(x),
$$
where $Q^{n}$ is the dyadic square of size $2^{-n}$, containing $x$.
Therefore, the claim is true if
$$
\langle \f\rangle\ci{Q^{n}}\to \f(x) \qquad\text{ and }\qquad
\langle h\rangle\ci{Q^{n}}\to h(x)\qquad \text{ as }n\to\infty
$$
(since for every transit square $Q$ the average
$\langle h\rangle\ci{Q}$ is close to $1$, we surely have $h(x)\ne 0$ for such
$x$).  But the exceptional set for any of these conditions has $\mu$-measure
$0$.

Now let us compare $\Lambda\f$ and $\Dl\ci Q\f$ to the corresponding terms in the
standard martingale decomposition, i.e., to
$$
\wt\Lambda\f :={\langle \f\rangle\ci {Q^0}}
$$
and
$$
\bigl.
\wt\Dl\ci Q\f\bigr|\ci{\C\setminus Q}:=0,\qquad
\bigl.
\wt\Dl\ci Q\f\bigr|\ci{Q_j}:=
\left\{
\aligned
\langle \f\rangle\ci {Q_j}
 -\langle \f\rangle\ci {Q}
, &\text{\quad if $Q_j$ is transit;}
\\
\f-\langle \f\rangle\ci {Q}
, &\text{\quad if $Q_j$ is terminal}
\endaligned\right.
$$
($j=1,2,3,4$).
It is well-known (and easy to prove) that
$$
\|\wt\Lambda\f\|^2\ci{L^2(\mu)}+\sum_{Q\in\D^{tr}}\|\wt\Dl\ci Q\f\|^2\ci{L^2(\mu)}
=\|\f\|^2\ci{L^2(\mu)}.
$$
A direct computation yields
$$
\|\wt\Lambda\f\|^2\ci{L^2(\mu)}=|\langle \f\rangle\ci{Q^0}|^2\mu(Q^0),\qquad
\|\Lambda\f\|^2\ci{L^2(\mu)}=
\frac{\langle |h|^2\rangle\ci {Q^0}}{|\langle h\rangle\ci{Q^0}|^2}
|\langle \f\rangle\ci {Q^0}|^2\mu(Q^0),
$$
i.e.,
$$
\|\Lambda\f\|^2\ci{L^2(\mu)}=
\frac{\langle |h|^2\rangle\ci {Q^0}}{|\langle h\rangle\ci{Q^0}|^2}
\|\wt\Lambda\f\|^2\ci{L^2(\mu)}.
$$
We are going to show that the ratio
$
\dfrac{\langle |h|^2\rangle\ci {Q^0}}{|\langle h\rangle\ci{Q^0}|^2}
$
is close to 1. Indeed, we can write
$$
\frac{\langle |h|^2\rangle\ci{Q^0}}{|\langle h\rangle\ci{Q^0}|^2}-1=
\frac{\langle |h|^2\rangle\ci{Q^0}
-|\langle h\rangle\ci{Q^0}|^2}
{|\langle h\rangle\ci{Q^0}|^2}=
\frac{\langle |g|^2\rangle\ci{Q^0}
-|\langle g\rangle\ci{Q^0}|^2}
{|\langle h\rangle\ci{Q^0}|^2}.
$$
Now note that
$$
|\langle h\rangle\ci {Q^0}|\ge
1-\langle |g|\rangle\ci {Q^0}
\ge 1-\sqrt{\langle |g|^2\rangle\ci {Q^0}}\ge 1-\delta,
$$
while the numerator is not less than $0$ (Cauchy inequality) and not greater
than $\langle |g|^2\rangle\ci {Q^0}\le\delta^2$. Therefore the whole ratio lies
between $0$ and $\delta^2(1-\delta)^{-2}\le \delta$.
So, we finally get
$$
\|\wt\Lambda\f\|^2\ci{L^2(\mu)}\le
\|\Lambda\f\|^2\ci{L^2(\mu)}
\le(1+\delta)\|\wt\Lambda\f\|^2\ci{L^2(\mu)}.
$$
As to the terms $\Dl\ci Q\f$, we will represent each of them as
the difference $\Dl'\ci Q\f-\dfrac{\langle \f\rangle\ci Q}
{\langle h\rangle\ci Q}h\ci Q$, where
$$
\bigl.
\Dl'\ci Q\f\bigr|\ci{\C\setminus Q}:=0,\qquad
\bigl.
\Dl'\ci Q\f\bigr|\ci{Q_j}:=
\left\{
\aligned
\frac{\langle \f\rangle\ci {Q_j}
 -\langle \f\rangle\ci {Q}}{\langle h\rangle\ci {Q_j}}
h, &\text{\quad if $Q_j$ is transit;}
\\
\f-\langle \f\rangle\ci {Q}
, &\text{\quad if $Q_j$ is terminal,}
\endaligned\right.
$$
and
$$
\bigl.
h\ci Q\bigr|\ci{\C\setminus Q}:=0,\qquad
\bigl.
h\ci Q\bigr|\ci{Q_j}:=
\left\{
\aligned
\frac{\langle h\rangle\ci {Q_j}
 -\langle h\rangle\ci {Q}}{\langle h\rangle\ci {Q_j}}
h, &\text{\quad if $Q_j$ is transit;}
\\
h-\langle h\rangle\ci {Q}
, &\text{\quad if $Q_j$ is terminal}
\endaligned\right.
$$
($j=1,2,3,4$).
Note that $\Dl'\f\equiv \wt\Dl\f$ on $\C\setminus Q$ and on every terminal
square $Q_j$. Also, if $Q_j$ is a transit subsquare of $Q$, then
$$
\int_{Q_j}|\wt\Dl\ci Q\f|^2d\mu\le
\int_{Q_j}|\Dl'\ci Q\f|^2d\mu
\le(1+\delta)
\int_{Q_j}|\wt\Dl\ci Q\f|^2d\mu
$$
(the reasoning is exactly the same as we had for $\Lambda\f$ and $\wt\Lambda\f$).
Using the elementary inequality
$$
\frac23|a|^2-2|b|^2\le |a-b|^2\le \frac32|a|^2+3|b|^2\qquad (a,b\in\C),
$$
we get
$$
\tfrac23\|\f\|^2\ci{L^2(\mu)}-2\sigma\le
\|\Lambda\f\|^2\ci{L^2(\mu)}+\sum_{Q\in\D^{tr}}\|\Dl\ci Q\f\|^2\ci{L^2(\mu)}  \le
\tfrac32(1+\delta)\|\f\|^2\ci{L^2(\mu)}+3\sigma,
$$
where
$$
\sigma:=\sum_{Q\in\D^{tr}}\frac{|\langle \f\rangle\ci {Q}|^2}
{|\langle h\rangle\ci {Q}|^2}\|h\ci Q\|^2\ci{L^2(\mu)}\le
\frac{1+\delta}{(1-\delta)^2}
\sum_{Q\in\D^{tr}}{|\langle \f\rangle\ci {Q}|^2}
\|\wt\Dl\ci Q g\|^2\ci{L^2(\mu)}\leq
$$
$$
2\sum_{Q\in\D^{tr}}{|\langle \f\rangle\ci {Q}|^2}
\|\wt\Dl\ci Q g\|^2\ci{L^2(\mu)},
$$
because  $|\langle h\rangle\ci {Q}|\ge 1-\delta$; the same reasoning as we
used when comparing $\Dl'\ci Q\f$ to $\wt\Dl\ci Q\f$, allows us to conclude
that $\|h\ci Q\|^2\ci{L^2(\mu)}\le (1+\delta)\|\wt\Dl\ci Q h\|^2\ci{L^2(\mu)}$,
and, at last, $\wt\Dl\ci Q h=\wt\Dl\ci Q g$.

Now let us remind the reader of the celebrated

{\bf Dyadic Carleson Imbedding Theorem}

Assume that we have a dyadic lattice $\D$ as above
and a family of nonnegative numbers
$\{a\ci Q\}\ci{Q\in\D}$. Suppose also that for every square $R\in\D$, we have
$$
\sum_{Q\in\D:Q\subset R}a\ci Q\le A\mu(R).
$$
Then for every function $\f\in L^2(\mu)$ we have
$$
\sum_{Q\in\D:\mu(Q)\ne 0}a\ci Q{|\langle \f\rangle\ci {Q}|^2}\le
4A\|\f\|^2\ci{L^2(\mu)}.
$$
Now observe that for every transit square $R\in\D$, we have
$$
\!\!\sum_{Q\in\D^{tr}:Q\subset R}
\|\wt\Dl\ci Q g\|^2\ci{L^2(\mu)}\!=
\sum_{Q\in\D^{tr}:Q\subset R}
\|\wt\Dl\ci Q (g\chi\ci R)\|^2\ci{L^2(\mu)}\le
\| g\chi\ci R\|^2\ci{L^2(\mu)}=\int_{R}|g|^2 d\mu\le \delta^2\mu(R).
$$
Thus, applying the Dyadic Carleson Imbedding Theorem to
$a\ci Q=\|\wt\Dl\ci Q g\|^2\ci{L^2(\mu)}$, if $Q$ is transit, and
$a\ci Q=0$, if $Q$ is terminal, we get
$$
\sigma\le 8\delta^2\|f\|\ci{L^2(\mu)}.
$$
To finish the proof of the lemma, it remains only to note that
$$
\frac{2}{3}-16\delta^2\ge\frac12\qquad\text{ and }\qquad
\frac{3}{2}(1+\delta)+24\delta^2\le 2.
$$

\bigskip

{\bf XIII. Functions $\Phi\ci\D$}

\bigskip

Recall that we already have the Lipschitz function $\wt\Phi$ and that
$\wt\Phi(x)\ge\dist(x,\C\setminus H)$. In particular it follows that
$$
\wt\Phi(x)\ge \dist(x,\partial Q)\qquad\text{for all }x\in Q,
$$
if $Q\in\D^{term}$ and $Q\subset H$.

We would like to extend this property to {\it all} terminal squares in $\D$.
So, let us define
$$
\Phi\ci\D(x)=\sup\{\wt\Phi(x), \dist(x,\C\setminus Q)\,:\,
Q\in\D^{term},\,\int_Q|g|^2d\mu\ge\delta^2\mu(Q)\}.
$$
Clearly, $\Phi\ci\D$ is Lipschitz, $\Phi\ci\D\ge \wt\Phi$, and
$
\Phi\ci\D(x)\ge \dist(x,\partial Q)
$
whenever $x\in Q\in \D^{term}$.

Now note that
$$
\mu\{x\in\C\,:\,\Phi\ci\D(x)>0\}\le \mu\{x\in\C\,:\,\wt\Phi(x)>0\}+
\sum_{Q\in\D^{term}, \int_Q|g|^2\ge\delta^2\mu(Q)}\mu(Q).
$$
The list in the end of
Section X shows
$\mu\{x\in\C\,:\,\wt\Phi(x)>0\}\le\frac{2\delta}{3}\mu(\C)$. On the other hand, the squares in $\D^{term}$ are
pairwise disjoint. Therefore the second sum does not exceed
$\delta^{-2}\int_\C|g|^2d\mu\le\delta^2\mu(\C)$, and we finally get
$$
\mu\{x\in\C\,:\,\Phi\ci\D(x)>0\}\le\bigl(\tfrac{2\delta}{3}+\delta^2
\bigr)\mu(\C)\le\delta\mu(\C).
$$

\bigskip

{\bf XIV. Action on good functions}

\bigskip

Now let $\D_1$ and $\D_2$ be two dyadic lattices of the above kind. We need
to show that for every two good functions $\f,\psi\in L^2(\mu)$ (they play
the roles of the functions $\f_1$ and $\f_2$
in the definition of perfect hair, respectively) and for every Lipschitz
function $\Theta\ge\delta\max\{\Phi\ci{\D_1},\Phi\ci{\D_2}\}$ satisfying
$\inf_{\C}\Theta>0$, one has
$$
|\langle\f,K\ci\Theta\psi\rangle|\le N\|\f\|\ci{L^2(\mu)}
\|\psi\|\ci{L^2(\mu)}.
$$
The reader may be surprised by the fact that we are talking about good
functions without defining the bad squares first. Actually, to tell the
truth, the bad squares are those with which we do not know what to do.
Almost all the statements below are very hard or even impossible to prove
directly for {\it arbitrary} squares $Q\in\D_1$, $R\in\D_2$. But they become
next to trivial, if we introduce some additional assumptions. All we need to
do is to show that all our auxiliary assumptions hold with probability close
to $1$, and this can be postponed to the very end.

Note first of all, that it is enough to prove the desired inequality
for functions $\f$ and $\psi$ such that $\Lambda\f=\Lambda\psi=0$.

Indeed, for any $\f\in L^2(\mu)$, we have
$$
\|K\ci\Theta\Lambda\f\|\ci{L^2(\mu)}=
\frac{|\langle\f\rangle\ci{Q^0}|}
{|\langle h\rangle\ci{Q^0}|}
\|K\ci\Theta h\|\ci{L^2(\mu)}\leq
$$
$$
\frac{1}{1-\delta}
|\langle\f\rangle\ci{Q^0}|\cdot
\|K\ci\Theta h\|\ci{L^\infty(\mu)}
\sqrt{\mu(Q^0)}\le
2B|\langle\f\rangle\ci{Q^0}|
\sqrt{\mu(Q^0)}\le 2 B
\|\f\|\ci{L^2(\mu)}.
$$
Taking into account that $\langle\f,K\ci\Theta\psi\rangle=-
\langle K\ci\Theta\f,\psi\rangle$ for all $\f,\psi\in L^2(\mu)$, we get
$$
\langle\f,K\ci\Theta\psi\rangle=
-\langle K\ci\Theta\Lambda\f,\psi\rangle+
\langle\f-\Lambda\f,K\ci\Theta\Lambda\psi\rangle+
\langle\f-\Lambda\f,K\ci\Theta(\psi-\Lambda\psi)\rangle.
$$
The first two terms do not exceed $2B\|\f\|\ci{L^2(\mu)}\|\psi\|\ci{L^2(\mu)}$
and $4B\|\f\|\ci{L^2(\mu)}\|\psi\|\ci{L^2(\mu)}$, correspondingly (because
$\|\f-\Lambda\f\|\ci{L^2(\mu)}\le 2\|\f\|\ci{L^2(\mu)}$). Meanwhile, the functions
$\f'=\f-\Lambda\f$ and $\psi'=\psi-\Lambda\psi$ clearly satisfy
the condition $\Lambda\f=\Lambda\psi=0$
and their $L^2(\mu)$-norms are bounded by $2\|\f\|\ci{L^2(\mu)}$ and
$2\|\f\|\ci{L^2(\mu)}$, respectively.
So, if we prove the desired inequality for all good
 $\f$ and $\psi$ satisfying
$\Lambda\f=\Lambda\psi=0$ with some constant $N_1$, then we will get it for two arbitrary
good functions with the constant $N=4N_1+6B$.

We would like to write
$$
\langle\f,K\ci\Theta\psi\rangle
=
\sum_{Q\in\D_1^{tr},\, R\in\D_2^{tr}}
\langle\Dl\ci Q\f,K\ci\Theta\Dl\ci R\psi\rangle.
$$
The question arises of why this series converges in any reasonable sense. But
let us observe that, since $\inf_\C\Theta>0$, the operator $K\ci\Theta$ is bounded
in $L^2(\mu)$ and therefore we can restrict ourselves to the good functions
$\f$ and $\psi$ that have only finitely many non-zero terms in their
decompositions (clearly, if $\f$ is good, then any partial sum of the series
$\Lambda\f+\sum_{Q\in\D_1}\Dl\ci Q\f$ is good as well). This not only removes any
questions about the convergence, but also allows us to rearrange and to group
the terms in the sum in any way we want.

Due to this observation and due to the (anti)symmetry, it is enough to
estimate the sum over $Q\in\D_1^{tr}$ and $R\in\D_2^{tr}$, for which
$l(Q)\leq l(R)$. For the sake of notational simplicity, everywhere below
instead of
$$
\sum_{Q\in\D_1^{tr},\,Q\text{ is good},
\,R\in\D_2^{tr},\,l(Q)\leq l(R),\text{ other conditions}}\ ,
$$
we will write
$$
\sum_{Q,R:\text{ other conditions}}\ .
$$
Also we will always reduce\ \ 
$\sum_{Q\in\D_1^{tr}:\, Q\text{ is good},\,\text{other conditions}}$  to $\sum_{Q\,:\,\text{other conditions}}$ 

\noindent  and  $\sum_{R\in\D_2^{tr}\,:\,
\text{other conditions}}$ to $\sum_{R\,:\,\text{other conditions}}$.

Note, that (unless otherwise specified) we will always think that the
summation over $Q$ goes only over {\it good} squares $Q\in\D_1^{tr}$,
while the summation over $R$ goes over {\it all} $R\in\D_2^{tr}$.

Of course, formally it doesn't matter, because, since the functions
$\f$ and $\psi$ are good, it is merely a business of  adding or omitting
several zeros. But it will allow us (and the reader) to see
clearly where and  what property  is used. As the reader might
have  already
guessed, for
the sum over pairs $Q,R$ with $l(Q)>l(R)$, this point of view
should be changed to the opposite.

Pick some large positive integer $m$ and write
$$
\sum_{Q,R}\langle\Dl\ci Q\f,K\ci\Theta\Dl\ci R\psi\rangle=
\sum_{Q,R:l(Q)\ge 2^{-m}l(R)}+\sum_{Q,R:l(Q)<2^{-m}l(R)}=
$$
$$
\sum\Sb Q,R:l(Q)\ge 2^{-m}l(R),
\\      \dist(Q,R)\le l(R)\endSb
+\Biggl[\,
\sum\Sb Q,R:l(Q)\ge 2^{-m}l(R),
\\      \dist(Q,R)> l(R)\endSb
+
\sum\Sb Q,R:l(Q)< 2^{-m}l(R),
\\      Q\cap R=\emptyset\endSb\,
\Biggr]
+
\sum\Sb Q,R:l(Q)< 2^{-m}l(R),
\\      Q\cap R\ne\emptyset\endSb
$$
$$
=:\sigma_1+\sigma_2+\sigma_3.
$$
Recall that the kernel $k\ci\Theta$ satisfies the estimates
$$
|k\ci\Theta(x,y)|\le\frac{1}{\max\{|x-y|,\Theta(x),\Theta(y)\}}
\qquad\text{and}\qquad|\nabla_x k\ci\Theta(x,y)\le\frac{4}{|x-y|^2}.
$$
The second inequality implies that
$$
|k\ci\Theta(x,y)-k\ci\Theta(x',y)|\le \frac{16|x-x'|}{|x-y|^2},
$$
provided that $|x-x'|\le\tfrac12|x-y|$.
Actually, we do not need the kernel to be that smooth. We will see that
the estimate
$$
|k\ci\Theta(x,y)-k\ci\Theta(x',y)|\le \frac{A|x-x'|^\e}{|x-y|^{1+\e}}
$$
with some (fixed) $0<\e\le 1$ and $0<A<+\infty$
is sufficient for all our tricks. The reader may ask:
``Why introduce a special notation for the parameter, which is actually
equal to 1; isn't it merely a generalization for the sake of
generalization?'' Well, first of all, we want to show that there is nothing
very special about the Cauchy kernel $\frac{1}{x-y}$; it can be replaced by
any other (antisymmetric) Calderon-Zygmund kernel. And secondly, it will
allow the reader to check that our proof works not because of some ``magic''
numerical identities like $\frac13-\frac12+\frac16=0$, but because we really
have found a good way to go around the main drawback of the Haar system:
the impossibility to make good estimates near jumps. And once this main drawback
is removed, the old-fashioned  Haar system becomes more
elegant and powerful than any ultramodern and superfamous wavelets.

\bigskip

{\bf XV. Estimation of $\sigma_2$}

\bigskip

Recall that the sum $\sigma_2$ is taken over pairs $Q,R$ such that
$Q\cap R=\emptyset$. If $l(Q)\ge 2^{-m}l(R)$, then the squares not only do not
intersect, but are well-separated: $\dist(Q,R)\ge l(R)$. We would like to
extend this property onto the case $l(Q)< 2^{-m}l(R)$. Though we
cannot achieve exactly the same separation by the length of the larger
square, we can get as close to it as we want. Namely, for any $\al>0$ and for
any $Q\in\D_1$, the probability
$$
P\ci{\D_2 }\{
\text{there exists }R\in\D_2\text{ : }l(R)> 2^ml(Q),
R\cap Q=\emptyset\, \text{ and }\,\dist(Q,R) < l(Q)^{\al}l(R)^{1-\al}
\}
$$
allows an estimate that does not depend on $Q$ and tends to $0$
as $m\to\infty$.

We shall need this result for $\al=\frac{\e}{2(1+\e)}$ ($\frac14$ in the case
of the Cauchy kernel). We will postpone the proof of this claim to the end of
the paper, as we said before; and now let us observe that if we declare the
corresponding squares $Q$ bad and if $\f$ is good, then
for {\it every} pair $Q,R$, participating in $\s_2$, we have
$\dist(Q,R)\ge l(Q)^{\al}l(R)^{1-\al}$.

Define the {\it long distance} $D(Q,R)$ between the squares $Q$ and $R$ by
$$
D(Q,R)=l(Q)+l(R)+\dist(Q,R).
$$
{\bf Far Interaction Lemma:}

Suppose that $Q$ and $R$ are two squares on the complex plane $\C$, such that
$l(Q) \leq l(R)$. Let $\f\ci
Q,\psi\ci R\in L^2(\mu)$. Assume that $\f\ci Q$ vanishes outside $Q$,
$\psi\ci R$
vanishes outside $R$; $\int_\C\f\ci Q=0$ and, at last,
$\dist(Q,\supp\psi\ci R)\ge l(Q)^{\al}l(R)^{1-\al}$.

Then
$$
|\langle \f\ci Q, K\ci\Theta \psi\ci R \rangle|\le
3^{1+\e}A\frac{l(Q)^{\frac\e2}l(R)^{\frac\e2}}{D(Q,R)^{1+\e}}
\sqrt{\mu(Q)}\sqrt{\mu(R)}\|\f\ci Q\|\ci{L^2(\mu)}\|\psi\ci R\|\ci{L^2(\mu)}.
$$

{\bf Remark}

Note, that we require only that the support of the function $\psi$ lies far
from $Q$; the squares $Q$ and $R$ themselves may intersect!
We will really have such a situation when estimating $\s_3$.

{\bf Proof:}

Let $x\ci Q$ be the center of the square $Q$. Note that for all $x\in Q$,
$y\in \supp\psi\ci R$, we have
$$
|x\ci Q-y|\ge \frac{l(Q)}{2}+\dist(Q,\supp\psi\ci R)\ge
\frac{3l(Q)}{2}\ge\sqrt2l(Q)\ge 2|x-x\ci Q|.
$$
Therefore,
$$
|\langle \f\ci Q, K\ci\Theta \psi\ci R \rangle|=
\Bigl|\iint k\ci\Theta(x,y)\f\ci Q(x)\psi\ci R(y)\, d\mu(x)\, d\mu(y)\Bigr|=
$$
$$
\Bigl|\iint [k\ci\Theta(x,y)-k\ci\Theta(x\ci Q,y)]
\f\ci Q(x)\psi\ci R(y)\, d\mu(x)\, d\mu(y)\Bigr|\leq
$$
$$
A\frac{l(Q)^\e}{\dist(Q,\supp\psi\ci R)^{1+\e}}
\|\f\ci Q\|\ci{L^1(\mu)}\|\psi\ci R\|\ci{L^1(\mu)}.
$$
There are two possible cases:

{\bf Case 1: $\dist(Q,\supp\psi\ci R)\ge l(R)$}

Then
$$
D(Q,R)=l(Q)+l(R)+\dist(Q,R)\le
3\dist(Q,\supp\psi\ci R)
$$
and therefore
$$
\frac{l(Q)^\e}{\dist(Q,\supp\psi\ci R)^{1+\e}}\le 3^{1+\e}
\frac{l(Q)^\e}{\D(Q,R)^{1+\e}}\le
3^{1+\e}\frac{l(Q)^{\frac\e2}l(R)^{\frac\e2}}{D(Q,R)^{1+\e}}.
$$

{\bf Case 2: $l(Q)^{\al}l(R)^{1-\al}\le \dist(Q,\supp\psi\ci R)\le l(R)$}

Then $D(Q,R)\le 3l(R)$ and we get
$$
\frac{l(Q)^\e}{\dist(Q,\supp\psi\ci R)^{1+\e}}\le
\frac{l(Q)^\e}{[l(Q)^{\al}l(R)^{1-\al}]^{1+\e}}=
\frac{l(Q)^{\frac\e2}l(R)^{\frac\e2}}{l(R)^{1+\e}}\le
3^{1+\e}\frac{l(Q)^{\frac\e2}l(R)^{\frac\e2}}{D(Q,R)^{1+\e}}.
$$
Now, to finish the proof of the lemma, it remains only to note that
$$
\|\f\ci Q\|\ci{L^1(\mu)}\le
\sqrt{\mu(Q)}\|\f\ci Q\|\ci{L^2(\mu)}
\text{\qquad and\qquad }
\|\psi\ci R\|\ci{L^1(\mu)}\le
\sqrt{\mu(R)}\|\psi\ci R\|\ci{L^2(\mu)}.
$$


Applying this lemma to $\f\ci Q=\Dl\ci Q\f$ and $\psi\ci R=\Dl\ci R\psi$, we
obtain
$$
|\s_2|\le
3^{1+\e}A
\sum_{Q,R}
\frac{l(Q)^{\frac\e2}l(R)^{\frac\e2}}{D(Q,R)^{1+\e}}
\sqrt{\mu(Q)}\sqrt{\mu(R)}\|\Dl\ci Q\f\|\ci{L^2(\mu)}
\|\Dl\ci R\psi\|\ci{L^2(\mu)} \qquad (**)
$$

We are going to show that the matrix $T\ci{Q,R}$ defined by
$$
T\ci{Q,R}:=
\frac{l(Q)^{\frac\e2}l(R)^{\frac\e2}}{D(Q,R)^{1+\e}}
\sqrt{\mu(Q)}\sqrt{\mu(R)}\qquad (Q\in\D_1^{tr},\,
R\in\D_2^{tr},\,l(Q)l(R)\,)
$$
generates a bounded operator in $l^2$.

{\bf Lemma:}

For any two
``sequences''  $\{a\ci Q\}\ci{Q\in\D_1^{tr}}$ and $\{b\ci R\}\ci{R\in\D_2^{tr}}$
of nonnegative numbers, one has
$$
\sum_{Q,R}T\ci{Q,R}a\ci Q b\ci R\le
\frac{3^{1+\e}(3+\e^{-1})M}{ 1-2^{-\frac\e2} }
\Bigl[\sum_{Q}a\ci Q^2\Bigr]^{\frac12}
\Bigl[\sum_{R}b\ci R^2\Bigr]^{\frac12}.
$$
{\bf Remark:}

Note that $T\ci{Q,R}$ are defined for all $Q,R$ with $l(Q)\le l(R)$
and that
the condition $\dist(Q,R)\ge l(Q)^{\al}l(R)^{1-\al}$ (or even the condition
$Q\cap R=\emptyset$)
no longer appears in the summation!

{\bf Proof:}

Let us ``slice'' the matrix $T\ci{Q,R}$ according to the ratio
$\frac{l(Q)}{l(R)}$. Namely, let
$$
T^{(n)}_{Q,R}=\left\{
\aligned
T\ci{Q,R},&\qquad\text{if } l(Q)=2^{-n}l(R);
\\
0,&\qquad\text{otherwise}
\endaligned
\right.
$$
($n=0,1,2,\dots$).
To prove the lemma, it is enough to show that for every $n\ge 0$,
$$
\sum_{Q,R}T^{(n)}\ci{Q,R}a\ci Q b\ci R\le
2^{-\frac\e2 n }3^{1+\e}(3+\e^{-1})M
\Bigl[\sum_{Q}a\ci Q^2\Bigr]^{\frac12}
\Bigl[\sum_{R}b\ci R^2\Bigr]^{\frac12}.
$$
The matrix $\{T^{(n)}_{Q,R}\}$ has a ``block'' structure: the variables
$b\ci R$ corresponding to the squares $R\in\D_2^{tr}$,
for which $l(R)=2^{j}$, can interact only with variables
$a\ci Q$ corresponding to the squares $Q\in\D_1^{tr}$,
for which $l(Q)=2^{j-n}$. Thus, to get the desired inequality, it is
enough to estimate each block separately, i.e., to demonstrate that
$$
\sum_{Q,R\,:\,l(Q)=2^{j-n},l(R)=2^j}T^{(n)}\ci{Q,R}a\ci Q b\ci R\leq
$$
$$
2^{-\frac\e2 n }3^{1+\e}
(3+\e^{-1})M
\Bigl[\sum_{Q\,:\,l(Q)=2^{j-n}}a\ci Q^2\Bigr]^{\frac12}
\Bigl[\sum_{R\,:\,l(R)=2^j}b\ci R^2\Bigr]^{\frac12}.
$$
Let us introduce the functions
$$
F:=\sum_{Q\,:\,l(Q)=2^{j-n}}\frac{a\ci Q}{ \sqrt{\mu(Q)} }
\chi\ci Q
\qquad
\text{and}
\qquad
G:=\sum_{R\,:\,l(R)=2^{j}}\frac{b\ci R}{ \sqrt{\mu(R)} }
\chi\ci R.
$$
Note that the squares of a given size in one dyadic lattice do not intersect,
and therefore at each point $x\in\C$, at most one term in the sum can be
non-zero. Also observe that
$$
\|F\|\ci{L^2(\mu)}=
\Bigl[\sum_{Q\,:\,l(Q)=2^{j-n}}a\ci Q^2\Bigr]^{\frac12}
\qquad
\text{and}
\qquad
\|G\|\ci{L^2(\mu)}=
\Bigl[\sum_{R\,:\,l(R)=2^{j}}b\ci R^2\Bigr]^{\frac12}.
$$
Then the estimate we need can be rewritten as
$$
\iint k_{j,n}(x,y)F(x)G(y)\, d\mu(x)\,d\mu(y) \le
2^{-\frac\e2 n }3^{1+\e}
(3+\e^{-1})M
\|F\|\ci{L^2(\mu)}\|G\|\ci{L^2(\mu)},
$$
where
$$
k_{j,n}(x,y)=\sum_{Q,R\,:\,l(Q)=2^{j-n}, l(R)=2^j}
\frac{l(Q)^{\frac\e2}l(R)^{\frac\e2}}{D(Q,R)^{1+\e}}\chi\ci Q(x)\chi\ci
R(y).
$$
Again, for every pair of points $x,y\in \C$, only one term in the sum can be
nonzero.
Since $|x-y|+l(R)\le 3D(Q,R)$ for any $x\in Q$, $y\in R$, we obtain
$$
k_{j,n}(x,y)=
2^{-\frac\e2 n }\frac{ l(R)^\e }{ D(Q,R)^{1+\e} }\le
2^{-\frac\e2 n }3^{1+\e}
\frac{2^{j\e}}{[2^j+|x-y|]^{1+\e}}
=:2^{-\frac\e2 n }3^{1+\e}k_j(x,y).
$$
So, it is enough to check that
$$
\iint k_{j}(x,y)F(x)G(y)\, d\mu(x)\,d\mu(y)
\le
(3+\e^{-1})M
\|F\|\ci{L^2(\mu)}\|G\|\ci{L^2(\mu)}.
$$
According to the Schur test, it would suffice to prove that
for every $y\in \C$, one has the estimate $\int_{\C}k_j(x,y)\,d\mu(x)\le
(3+\e^{-1})M$ and vice versa (i.e.,
for every $x\in \C$, one has $\int_{\C}k_j(x,y)\,d\mu(y)\le
(3+\e^{-1})M$). Then the norm of the integral operator with kernel $k_j$
in $L^2(\mu)$ would be bounded by the same constant
$(3+\e^{-1})M$, and the story would be over.

If we assumed a priori that $\mathcal{R}(y)\le 2^{j+1}$, then the needed estimate would
be next to trivial: we could write
$$
\int_{\C}k_j(x,y)\,d\mu(x)=
\int_{B(y,2^{j+1})}k_j(x,y)\,d\mu(x)+
\int_{\C\setminus B(y,2^{j+1})}k_j(x,y)\,d\mu(x)\leq
$$
$$
2^{-j}\mu(B(y,2^{j+1}))+
\int_{\C\setminus B(y,2^{j+1})}\frac{2^{j\e}}{|x-y|^{1+\e}}\,d\mu(x)\leq
$$
$$
M\Bigl(2+1+\int_{2^j}^{+\infty}\frac{2^{j\e}}{t^{1+\e}}dt\Bigr)
=(3+\e^{-1})M
$$
(we applied Comparison Lemma to estimate the integral over $\C\setminus
B(y,2^{j+1})$, and again we used the possibility to switch from the radius
$2^{j+1}$ to the smaller number $2^j$)

The problem is that we cannot guarantee that $\mathcal{R}(y)\le 2^{j+1}$ for
{\it every} $y\in\C$. So, generally speaking, we are unable to show that the
integral operator with kernel $k_j(x,y)$ acts in $L^2(\mu)$. But we
{\it do not
need\/} that much! We only need to check that the corresponding bilinear form
 is bounded on two {\it given} functions $F$ and $G$. So, we are not
interested in the points $y\in\C$ for which $G(y)=0$ (or in the points
$x\in\C$, for which $F(x)=0$). But, by  definition,
$G$ can be non-zero on transit squares in $\D_2$ of size $2^j$ only.
Now let us notice that if $R\in\D_2^{tr}$, then
$\mathcal{R}(y)\le 2l(R)$ for every ${y\in R}$. Indeed, otherwise there exists a
non-Ahlfors disk $B(y,r)$ of radius $r>2l(R)$. But then
$R\subset B(y,r)\subset H$, which is impossible for a transit square!

The same reasoning shows that $\mathcal{R}(x)\le 2^{j-n+1}\le 2^{j+1}$ whenever
$F(x)\ne 0$, and we are done with $|\s_2|$.

\medskip

Now, we hope, the reader will agree that the decision to declare the squares
contained in $H$ terminal was a good one: not only does  the fact that
the measure $\mu$ is not Ahlfors not put us in any real trouble, but we just
hardly have a chance to notice this fact at all. Also, it is clear
why the squares with large average of $|g|^2$ have been declared terminal: this
allowed us to treat $h$ like an accretive function all the time.

But it still remains unexplained why we were so eager to suppress the Cauchy
kernel on every terminal square. The answer is in the next two sections.

\bigskip

{\bf XVI. Estimation of $\s_3$}

\bigskip

Recall that the sum $\s_3$ is taken over the pairs $Q,R$, for which
$l(Q)<2^{-m}l(R)$ and $Q\cap R\ne\emptyset$. We would like to improve
this condition to the demand that $Q$ lie ``deep inside'' one of the four
subsquares $R_j$ ($j=1,2,3,4$).

Define the {\it skeleton} $\sc R$ of the square $R$ by
$$
\text{sk} R:=\bigcup_{j=1}^4\partial R_j.
$$
We will declare a square $Q\in\D_1$ bad if there exists a square $R\in\D_2$
such that $l(R)>2^ml(Q)$ and $\dist(Q,\text{sk} R)\le
8l(Q)^{\al}l(R)^{1-\al}$. Note that any square bad in the sense of the
previous section is bad in this new sense as well.

Now, for every good square $Q\in\D_1$, the conditions $l(Q)<2^{-m}l(R)$
and
 $Q\cap R\ne \emptyset$
together imply
 that $Q$ lies inside one of the four subsquares $R_j$. We will denote
this subsquare by $R\ci Q$. The sum $\s_3$ can now be split into
$$
\s_3^{term}:=\sum\Sb Q,R\,:\,Q\subset R,\, l(Q)<2^{-m}l(R),
                    \\ \mathcal{R}_Q\text{ is terminal} \endSb
\langle\Dl\ci Q\f, K\ci\Theta\Dl\ci R\psi\rangle
$$
and
$$
\s_3^{tr}:=\sum\Sb Q,R\,:\,Q\subset R,\, l(Q)<2^{-m}l(R),
                    \\  \mathcal{R}_Q\text{ is transit} \endSb
\langle\Dl\ci Q\f, K\ci\Theta\Dl\ci R\psi\rangle.
$$

\bigskip

{\bf XVII. Estimation of $\s_3^{term}$}

\bigskip

First of all, write (recall that $R_j$ denote the children of $R$):
$$
\s_3^{term}=\sum_{j=1}^4\,\,
\sum\Sb Q,R\,:\,l(Q)<2^{-m}l(R),
\\
Q\subset R_j\in\D_2^{term}\endSb
\langle\Dl\ci Q\f, K\ci\Theta\Dl\ci R\psi\rangle.
$$
Clearly, it is enough to estimate the inner sum for every fixed $j$. Let us
do it for $j=1$. We have
$$
\sum\Sb Q,R\,:\,l(Q)<2^{-m}l(R), \\ Q\subset
R_1\in\D_2^{term}\endSb
\langle\Dl\ci Q\f, K\ci\Theta\Dl\ci R\psi\rangle=
\sum_{R:R_1\in\D_2^{term}}\,\, \sum\Sb Q:\,l(Q)<2^{-m}l(R),
\\ 
Q\subset R_1\endSb
\langle\Dl\ci Q\f, K\ci\Theta\Dl\ci R\psi\rangle.
$$
Roughly speaking, our main idea here is the following. If
$R_1\in\D_2^{term}$, then for all $x\in R_1$, one has
$$
\Theta(x)\ge\delta\Phi\ci{\D_2}(x)\ge\delta\dist(x,\partial R_1).
$$
For the points $x$ that lie in the ``central part'' of $R_1$, the right hand
side is at least $\frac{\delta(R)}{8}$. Assume that it is so for {\it every}
point $x\in R_1$. Then
$$
k\ci\Theta(x,y)\le\frac{1}{\Theta(x)}\le\frac{8}{\delta(R)}\qquad\text{for all }
x\in R_1, y\in\C.
$$
Hence
$$
|K\ci\Theta\Dl\ci R\psi(x)|\le\frac{8\|\Dl\ci R\psi \|\ci{L^1(\mu)}}{\delta(R)}
\qquad\text{ for all }x\in R_1,
$$
and therefore
$$
\|\chi\ci{R_1}\cdot K\ci\Theta\Dl\ci R\psi \|\ci{L^2(\mu)}\le
8\|\Dl\ci R\psi \|\ci{L^1(\mu)}\frac{\sqrt{\mu(R_1)}}{\delta(R)}\le
\frac{8{\mu(R)}}{\delta(R)}\|\Dl\ci R\psi \|\ci{L^2(\mu)}\le
\frac{8M}{\delta}\|\Dl\ci R\psi \|\ci{L^2(\mu)},
$$
because $\mu(R_1)\le \mu(R)$, $\|\Dl\ci R\psi \|\ci{L^1(\mu)}
\le \sqrt{\mu(R)} \|\Dl\ci R\psi \|\ci{L^2(\mu)}$, and $\mu(R)\le Ml(R)$
(otherwise the disk of radius $l(R)$, centered at the same point as $R$,
would be non-Ahlfors, and we would have $R\subset H$, which is impossible).

Now, recalling the remark from Section II, and taking into account that
$\Dl\ci Q\f\equiv 0$ outside $Q$, we get
$$
\sum_{Q:\,Q\subset
R_1}
|\langle\Dl\ci Q\f, K\ci\Theta\Dl\ci R\psi\rangle|=
\sum_{Q:\,Q\subset
R_1}
|\langle\Dl\ci Q\f, \chi\ci{R_1}\cdot K\ci\Theta\Dl\ci R\psi\rangle|\leq
$$
$$
\sqrt2\|\chi\ci{R_1}\cdot K\ci\Theta\Dl\ci R\psi \|\ci{L^2(\mu)}
\Bigl[\sum_{Q:\,Q\subset
R_1}
\|\Dl\ci Q\f
\|^2\ci{L^2(\mu)}\Bigr]^{\frac{1}{2}}\leq
$$
$$
\frac{16M}{\delta}\|\Dl\ci R\psi \|\ci{L^2(\mu)}
\Bigl[\sum_{Q:\,Q\subset
R_1}
\|\Dl\ci Q\f
\|^2\ci{L^2(\mu)}\Bigr]^{\frac{1}{2}}.
$$
So, we obtain
$$
\sum_{R:\,R_1\in\D_2^{term}}\,
\sum_{Q:\,Q\subset
R_1}
|\langle\Dl\ci Q\f, K\ci\Theta\Dl\ci R\psi\rangle|
\le
$$
$$
\frac{16M}{\delta}
\sum_{R:\,R_1\in\D_2^{term}}
\|\Dl\ci R\psi \|\ci{L^2(\mu)}
\Bigl[\sum_{Q:\,Q\subset
R_1}
\|\Dl\ci Q\f
\|^2\ci{L^2(\mu)}\Bigr]^{\frac{1}{2}}\le
$$
$$
\frac{16M}{\delta}
\Bigl[\sum_{R:\,R_1\in\D_2^{term}}
\|\Dl\ci R\psi \|^2\ci{L^2(\mu)}\Bigr]^{\frac12}
\Bigl[
\sum_{R:\,R_1\in\D_2^{term}}\,\,
\sum_{Q:\,Q\subset
R_1}
\|\Dl\ci Q\f
\|^2\ci{L^2(\mu)}\Bigr]^{\frac{1}{2}}.
$$
But the terminal squares in $\D_2$ do not intersect! Therefore every
$\Dl\ci Q\f$ can appear at most once in the last double sum, and we get the
bound
$$
\sum_{R:\,R_1\in\D_2^{term}}
\sum_{Q:\,Q\subset
R_1}
|\langle\Dl\ci Q\f, K\ci\Theta\Dl\ci R\psi\rangle|
\le
$$
$$
\frac{16M}{\delta}
\Bigl[\sum_{R}
\|\Dl\ci R\psi \|^2\ci{L^2(\mu)}\Bigr]^{\frac12}
\Bigl[
\sum_{Q}
\|\Dl\ci Q\f\|^2\ci{L^2(\mu)}\Bigr]^{\frac{1}{2}}\le
\frac{32M}{\delta}
\|\f\|\ci{L^2(\mu)}\|\psi\|\ci{L^2(\mu)}.
$$

The problem is that we cannot guarantee the estimate $\Theta(x)\ge
\frac{\delta(R)}{8}$ for {\it every} point $x\in R_1$. So, the kernel
$k\ci\Theta$ can grow near the boundary. Nevertheless, due to our definition of
good squares, we need only to consider the squares $Q\subset R_1$, for which
$\dist(Q,\partial R_1)\ge 8l(Q)^{\al}l(R)^{1-\al}$. So, if such a square
$Q$ lies close to the boundary of $R_1$, the size $l(Q)$ has to be very
small and the corresponding function $\Dl\ci Q\f$ should oscillate very fast.
We may hope that this fast oscillation will compensate for the
growth of the kernel. To show that it is really the case, we need one more
standard technical tool.

\bigskip

{\bf XVIII. Whitney decomposition}

\bigskip

Let $S^0$ be an arbitrary square on the complex plane $\C$. Consider the
standard dyadic lattice starting with the square $S^0$, and denote by
$W(S^0)$ the family of all maximal subsquares $S$ in this lattice, for which
$\dist(S,\partial S^0)\ge l(S)$ (see Picture 2). The Whitney decomposition
$W(S^0)$ has the following remarkable properties:
\medskip\it
1) The squares $S\in W(S^0)$ are pairwise disjoint and cover the interior
of $S^0$;
\smallskip
2) $\dist(S,\partial S^0)=l(S)$ for every $S\in W(S^0)$;
\smallskip
3) The expanded squares $\wt S:=2S$ $(S\in W(S^0)\,)$ still lie ``deep 
inside'' $S$, namely, $\dist(\wt S,\partial S^0)=\frac{l(S)}{2}=\frac{l(\wt S)}{4}$,
and every point $x\in \C$ belongs to at most $6$ squares $\wt S$.
\medskip\rm
Denote again the center of a square $Q$ by $x\ci Q$. For $S\in W(R_1)$ put
$$
\psi\ci{R,S}:=\chi\ci{\wt S}\,\Dl\ci R\psi\qquad\text{ and }\qquad
\wt\psi\ci{R,S}:=\chi\ci{R\setminus \wt S}\,\Dl\ci R.
$$
We have
$$
\sum\Sb
Q:\,l(Q)<2^{-m}l(R),
\\ Q\subset R_1
\endSb
\langle \Dl\ci Q\f, K\ci\Theta \Dl\ci R\psi\rangle=
\sum_{S\in W(R_1)}\,\,
\sum\Sb
Q:\,l(Q)<2^{-m}l(R),
\\ Q\subset R_1,\, x\cci Q\in S
\endSb
\langle \Dl\ci Q\f, K\ci\Theta \Dl\ci R\psi\rangle
=
$$
$$
\sum_{S\in W(R_1)}\,\,
\sum\Sb
Q:\,l(Q)<2^{-m}l(R),
\\ Q\subset R_1,\, x\cci Q\in S
\endSb
\langle \Dl\ci Q\f, K\ci\Theta \psi\ci{R,S}\rangle+
\sum_{S\in W(R_1)}\,\,
\sum\Sb
Q:\,l(Q)<2^{-m}l(R),
\\ Q\subset R_1,\, x\cci Q\in S
\endSb
\langle \Dl\ci Q\f, K\ci\Theta \wt\psi\ci{R,S}\rangle.
$$
Note now that for every good $Q\subset R_1$ such that $x\ci Q\in S\in
W(R_1)$, one has
$$
8l(Q)\le 8l(Q)^\al l(R)^{1-\al}\le \dist(Q,\partial R_1)\le
\dist(x\ci Q,\partial R_1)\le 2l(S),
$$
and therefore
$$
\dist(Q,\supp \wt\psi\ci{R,S})\ge \dist(Q,\partial \wt S)\ge
\frac{l(S)-l(Q)}{2}\ge \frac{l(S)}{4}\ge l(Q)^\al l(R)^{1-\al}.
$$
Now the Far Interaction Lemma yields
$$
|\langle \Dl\ci Q\f, K\ci\Theta \wt\psi\ci{R,S}\rangle|\le
3^{1+\e}A\frac{l(Q)^{\frac\e 2} l(R)^{\frac\e 2}}{D(Q,R)^{1+\e}}
\sqrt{\mu(Q)}\sqrt{\mu(R)}\|\Dl\ci Q\f\|\ci{L^2(\mu)}
\|\wt\psi\ci{R,S}\|\ci{L^2(\mu)}.
$$
Taking into account that $\|\wt\psi\ci{R,S}\|\ci{L^2(\mu)}\le
\|\Dl\ci R\psi\|\ci{L^2(\mu)}$ and summing over all $R\in\D_2^{tr}$, we
arrive at the same sum as in the long term interaction of Section XV (actually, we arrive
at the part of that sum which {\it has not been used yet}, but 
{\it has already been estimated} there).

So, it remains to find a good upper bound for
$$
\sum_{S\in W(R_1)}\,\,
\sum\Sb
Q:\,l(Q)<2^{-m}l(R),
\\ Q\subset R_1,\, x\cci Q\in S
\endSb
\langle \Dl\ci Q\f, K\ci\Theta \psi\ci{R,S}\rangle.
$$
Observe once more that the conditions $Q\in\D_1^{tr}$, $Q$ is good,
$l(Q)<2^{-m}l(R)$, $Q\subset R_1$ and $x\ci Q\in S$ together
imply $Q\subset \wt S$ (as we have seen above, they even imply
that $Q$ lies deep inside $\wt S$).
So, it is enough to estimate the sum
$$
\sum_{S\in W(R_1)}\,\,
\sum\Sb
Q:\,Q\subset \wt S,\, x\cci Q\in S
\endSb
|\langle \Dl\ci Q\f, K\ci\Theta \psi\ci{R,S}\rangle|.
$$

Note now that for {\it every} $x\in\wt S$, we have
$$
\Theta(x)\ge\delta\dist(x,\partial R_1)\ge \frac{\delta(\wt S)}{4}.
$$
Recall that the ``naive" reasoning from Section XVII could not be used for the whole $R_1$. 
But it can be used for $\wt S$.
Repeating our ``naive'' reasoning from Section XVII for the square $\wt S$ instead of the
whole $R_1$, we obtain
$$
\sum\Sb
Q:\,Q\subset \wt S,
\, x\cci Q\in S
\endSb
|\langle \Dl\ci Q\f, K\ci\Theta \psi\ci{R,S}\rangle|\le
\|\chi\ci{\wt S}\cdot K\ci\Theta\psi\ci{R,S}\|\ci{L^2(\mu)}
\Bigl[
\sum\Sb
Q:\,Q\subset \wt S,
\, x\cci Q\in S
\endSb
\|\Dl\ci Q\f\|^2\ci{L^2(\mu)}
\Bigr]^{\frac12}
\le
$$
$$
\frac{4\mu(\wt S)}{\delta(\wt S)}
\|\psi\ci{R,S}\|\ci{L^2(\mu)}
\Bigl[
\sum\Sb
Q:\,Q\subset \wt S,
\, x\cci Q\in S
\endSb
\|\Dl\ci Q\f\|^2\ci{L^2(\mu)}
\Bigr]^{\frac12}.
$$
We would like to say again that $\mu(\wt S)\le Ml(\wt S)$. If not, then,
of course, we can conclude that $\wt S\subset H$,
but this {\it does not} yield a contradiction immediately, because $\wt S$
{\it is not\/} a transit square in $\D_2$ (actually, it is not in $\D_2$ at
all!). Note, nevertheless, that if we have at least one good square
$Q\in\D_1^{tr}$ such that $Q\subset \wt S$ (otherwise the sum
is $0$, and we have nothing to worry about), then we can extend the above chain of
inclusions to $Q\subset\wt S\subset H$,
which {\it is} a contradiction! So, as before, despite the fact that we cannot use the
Ahlfors condition {\it whenever we want to}, we can use it
{\it whenever we need to}.

Thus, we get
$$
\!\!\sum_{S\in W(R_1)}\,\,
\sum\Sb
Q:\, Q\subset \wt S,\\ x\cci Q\in S
\endSb
|\langle \Dl\ci Q\f, K\ci\Theta \psi\ci{R,S}\rangle|\le
\frac{4M}{\delta}
\sum_{S\in W(R_1)}\|\psi\ci{R,S}\|\ci{L^2(\mu)}
\Bigl[
\sum\Sb
Q:\,Q\subset \wt S,\\
x\cci Q\in S
\endSb
\|\Dl\ci Q\f\|^2\ci{L^2(\mu)}
\Bigr]^{\frac12}
\le
$$
$$
\frac{4M}{\delta}
\Bigl[
\sum_{S\in W(R_1)}\|\psi\ci{R,S}\|^2\ci{L^2(\mu)}
\Bigr]^{\frac12}
\Bigl[
\sum_{S\in W(R_1)}\,\,
\sum\Sb
Q:\,Q\subset R_1,
\, x\cci Q\in S
\endSb
\|\Dl\ci Q\f\|^2\ci{L^2(\mu)}
\Bigr]^{\frac12}
$$
(we relaxed the condition $Q\subset \wt S$ in the last sum to $Q\subset R_1$;
it causes no harm now).
But
$$
\sum_{S\in W(R_1)}\|\psi\ci{R,S}\|^2\ci{L^2(\mu)}=
\sum_{S\in W(R_1)}\int_{\wt S}|\Dl\ci R\psi|^2\,d\mu\le
6\int_\C|\Dl\ci R\psi|^2\,d\mu=6\|\Dl\ci R\psi\|^2\ci{L^2(\mu)}
$$
(because every point lies in not more than $6$ squares $\wt S$).

Meanwhile,
$$
\sum_{S\in W(R_1)}\,\,
\sum\Sb
Q:\,Q\subset R_1,
\, x\cci Q\in S
\endSb
\|\Dl\ci Q\f\|^2\ci{L^2(\mu)}
=
\sum_{Q:\,Q\subset R_1}
\|\Dl\ci Q\f\|^2\ci{L^2(\mu)}.
$$
Hence, summing over all $R\in\D_2^{tr}$, for which $R_1\in \D_2^{term}$,
we get
$$
\sum_{R:\,R_1\in\D_2^{term}}\,
\sum_{S\in W(R_1)}\,\,
\sum\Sb
Q:\, Q\subset \wt S,\\ x\cci Q\in S
\endSb
|\langle \Dl\ci Q\f, K\ci\Theta \psi\ci{R,S}\rangle|
\le
$$
$$
\frac{4\sqrt6 M}{\delta}
\sum_{R:\,R_1\in\D_2^{term}}\,
\|\Dl\ci R\psi\|\ci{L^2(\mu)}
\Bigl[
\sum_{Q:\,Q\subset R_1}
\|\Dl\ci Q\f\|^2\ci{L^2(\mu)}
\Bigr]^{\frac12}
\le
$$
$$
\frac{10M}{\delta}
\Bigl[
\sum_{R:\,R_1\in\D_2^{term}}
\|\Dl\ci R\psi\|^2\ci{L^2(\mu)}
\Bigr]^{\frac12}
\Bigl[
\sum_{R:\,R_1\in\D_2^{term}}\,
\sum_{Q:\,Q\subset R_1}
\|\Dl\ci Q\f\|^2\ci{L^2(\mu)}
\Bigr]^{\frac12}
\le
$$
$$
\frac{10M}{\delta}
\Bigl[
\sum_{R}\,
\|\Dl\ci R\psi\|^2\ci{L^2(\mu)}
\Bigr]^{\frac12}
\Bigl[
\sum_{Q}
\|\Dl\ci Q\f\|\ci{L^2(\mu)}
\Bigr]^{\frac12}
\le
\frac{20M}{\delta}
\|\f\|\ci{L^2(\mu)}
\|\psi\|\ci{L^2(\mu)},
$$
finishing the story with $\s_3^{term}$.

\bigskip

{\bf XIX. Estimation of $\s_3^{tr}$}

\bigskip

Recall that
$$
\s_3^{tr}=\sum\Sb Q,R\,:\,Q\subset R,\, l(Q)<2^{-m}l(R),
                    \\  \mathcal{R}_Q\text{ is transit} \endSb
\langle\Dl\ci Q\f, K\ci\Theta\Dl\ci R\psi\rangle.
$$
Split every term in the sum as
$$
\langle\Dl\ci Q\f, K\ci\Theta\Dl\ci R\psi\rangle=
\langle\Dl\ci Q\f, K\ci\Theta(\chi\ci{R\cci Q}\Dl\ci R\psi)\rangle+
\langle\Dl\ci Q\f, K\ci\Theta(\chi\ci{R\setminus R\cci Q}\Dl\ci R\psi)\rangle.
$$
Observe that since $Q$ is good, $Q\subset R$ and $l(Q)<2^{-m}l(R)$, we
have
$$
\dist(Q,\supp \chi\ci{R\setminus R\cci Q}\Dl\ci R\psi)\ge
\dist(Q,\text{sk} R)\ge l(Q)^{\al}l(R)^{1-\al}.
$$
Using the Far Interaction Lemma and taking into account that the norm
$\|\chi\ci{R\setminus R\cci Q}\Dl\ci R\psi\|\ci{L^2(\mu)}$ does not exceed
$\|\Dl\ci R\psi\|\ci{L^2(\mu)}$, we conclude that the sum
$$
\sum\Sb Q,R\,:\,Q\subset R,\, l(Q)<2^{-m}l(R),
                    \\  \mathcal{R}_Q\text{ is transit} \endSb
|\langle\Dl\ci Q\f, K\ci\Theta(\chi\ci{R\setminus R\cci Q}\Dl\ci R\psi)\rangle|
$$
can be estimated by the sum $(**)$ from Section XV.

Thus, our task is to find a good bound for the sum
$$
\sum\Sb Q,R\,:\,Q\subset R,\, l(Q)<2^{-m}l(R),
                    \\  \mathcal{R}_Q\text{ is transit} \endSb
\langle\Dl\ci Q\f, K\ci\Theta(\chi\ci{R\cci Q}\Dl\ci R\psi)\rangle.
$$
Recalling the definition of $\Dl\ci R\psi$ and recalling that $R\ci Q$ is a
{\it transit\/} square, we get
$$
\chi\ci{R\cci Q}\Dl\ci R\psi=c\ci{R,Q}\chi\ci{R\cci Q}h,
$$
where
$$
c\ci{R,Q}=\frac{\langle \psi\rangle\ci{R\cci Q}}
{\langle h\rangle\ci{R\cci Q}}-
\frac{\langle \psi\rangle\ci{R}}{\langle h\rangle\ci{R}}
$$
is a {\it constant}.
So, our sum can be rewritten as
$$
\sum\Sb Q,R\,:\,Q\subset R,\, l(Q)<2^{-m}l(R),
                    \\  \mathcal{R}_Q\text{ is transit} \endSb
c\ci{R,Q}\langle\Dl\ci Q\f, K\ci\Theta(\chi\ci{R\cci Q}h)\rangle.
$$
Our next aim will be to extend the function $\chi\ci{R\cci Q}h$ to the whole
function $h$ in every term (which is exactly the opposite of the idea of the
previous section, where, in a similar situation, we tried to ``shrink'' the
function $\Dl\ci R\psi$ to $\psi\ci{R,S}$).

Let us observe that
$$
\langle\Dl\ci Q\f, K\ci\Theta(\chi\ci{\C\setminus R\cci Q}h)\rangle=
\int_{\C\setminus R\cci Q}k\ci\Theta(x,y)\Dl\ci Q\f(x) h(y)\,d\mu(x)\,d\mu(y)
=
$$
$$
\int_{\C\setminus R\cci Q}
[k\ci\Theta(x,y)-k\ci\Theta(x\ci Q,y)]\Dl\ci Q\f(x) h(y)\,d\mu(x)\,d\mu(y).
$$
Note again that for every $x\in Q$, $y\in \C\setminus R\ci Q$, we have
$$
|x\ci Q-y|\ge \frac{l(Q)}{2}+\dist(Q,\C\setminus R\ci Q)\ge
\frac{3l(Q)}{2}\ge\sqrt2l(Q)\ge 2|x-x\ci Q|.
$$
Therefore
$$
|k\ci\Theta(x,y)-k\ci\Theta(x\ci Q,y)|\le \frac{A|x-x\ci Q|^\e}{|x\ci Q-y|^{1+\e}}
\le \frac{Al(Q)^\e}{|x\ci Q-y|^{1+\e}},
$$
and
$$
|\langle\Dl\ci Q\f, K\ci\Theta(\chi\ci{\C\setminus R\cci Q}h)\rangle|\le
Al(Q)^\e\|\Dl\ci Q\f\|\ci{L^1(\mu)}
\int_{\C\setminus R\cci Q}
\frac{|h(y)|\,d\mu(y)}{|x\ci Q-y|^{1+\e}}.
$$
Now let us consider the sequence of squares $R^{(j)}\in\D_2$, beginning with
$R^{(0)}=R\ci Q$ and gradually ascending ($R^{(j)}\subset R^{(j+1)}$,
$l(R^{(j+1)})=2l(R^{(j)})$) to the starting square $R^0=R^{(N)}$ of the
lattice $\D_2$. Clearly, all the squares $R^{(j)}$ are transit.

We have
$$
\int_{\C\setminus R\cci Q}
\frac{|h(y)|\,d\mu(y)}{|x\ci Q-y|^{1+\e}}=
\int_{R^0\setminus R\cci Q}
\frac{|h(y)|\,d\mu(y)}{|x\ci Q-y|^{1+\e}}=
\sum_{j=1}^N
\int_{R^{(j)}\setminus R^{(j-1)}}
\frac{|h(y)|\,d\mu(y)}{|x\ci Q-y|^{1+\e}}=:
\sum_{j=1}^N I_j.
$$
Note now that, since $Q$ is good and $l(Q)<2^{-m}l(R)\le
2^{-m}l(R^{(j)})$ for all $j=1,\dots,N$, we have
$$
\dist(Q,R^{(j)}\setminus R^{(j-1)})\ge
\dist(Q,\sc R^{(j)})\ge
l(Q)^{\al}l(R^{(j)})^{1-\al}.
$$
Hence
$$
I_j\le
\frac{1}{[l(Q)^{\al}l(R^{(j)})^{1-\al}]^{1+\e}}\int_{R^{(j)}}|h|\,d\mu.
$$
Recalling that $\al=\frac{\e}{2(1+\e)}$, we see that the first factor equals
$\dfrac{1}{l(Q)^{\frac\e2} l(R^{(j)})^{1+\frac\e2} }$.

Since $R^{(j)}$ is transit, we have
$$
\int_{R^{(j)}}|h|\,d\mu\le\int_{R^{(j)}}(1+|g|)\,d\mu
\le(1+\delta)\mu(R^{(j)})\le
(1+\delta)Ml(R^{(j)}).
$$
Thus,
$$
I_j\le\frac{(1+\delta)M}{l(Q)^{\frac\e2}l(R^{(j)})^{\frac\e2}}=
2^{-(j-1)\frac{\e}{2}}\frac{(1+\delta)M}{l(Q)^{\frac\e2}l(R)^{\frac\e2}}.
$$
Summing over $j\ge 1$, we get
$$
\int_{\C\setminus R\cci Q}
\frac{|h(y)|\,d\mu(y)}{|x\ci Q-y|^{1+\e}}=
\sum_{j=1}^N I_j\le
\frac{(1+\delta)M}{1-2^{-\frac\e2}}
\frac{1}{l(Q)^{\frac\e2}l(R)^{\frac\e2}}.
$$
Now let us note that, since $R\ci Q\in \D_2^{tr}$, we have
$$
\|\Dl\ci R\psi\|^2\ci{L^2(\mu)}\ge
\int_{R\cci Q}|\Dl\ci R\psi\|^2\,d\mu=
|c\ci{Q,R}|^2
\int_{R\cci Q}|h|^2\,d\mu
\ge
$$
$$
|c\ci{Q,R}|^2
|\langle h\rangle|^2\ci{R\cci Q}\mu(R\ci Q)\ge (1-\delta)^{2}
|c\ci{Q,R}|^2\mu(R\ci Q).
$$
So,
$$
|c\ci{Q,R}|\le\frac{1}{1-\delta}
\frac{\|\Dl\ci R\psi\|\ci{L^2(\mu)} }{\sqrt{\mu(R\ci Q)}}.
$$
Combining this estimate with the Cauchy inequality
$\|\Dl\ci Q\f\|\ci{L^1(\mu)}\le\sqrt{\mu(Q)}\|\Dl\ci Q\f\|\ci{L^2(\mu)}$, we
finally obtain
$$
|\langle\Dl\ci Q\f, K\ci\Theta(\chi\ci{\C\setminus R\cci Q}h)\rangle|\le
\frac{(1+\delta)MA}{(1-\delta)(1-2^{-\frac\e2})}
\left[\frac{l(Q)}{l(R)}\right]^{\frac\e2}
\sqrt{\frac{\mu(Q)}{\mu(R\ci Q)}}
\|\Dl\ci Q\f\|\ci{L^2(\mu)}
\|\Dl\ci R\psi\|\ci{L^2(\mu)}
$$
and
$$
\sum\Sb Q,R\,:\,Q\subset R,\, l(Q)<2^{-m}l(R),
                    \\ \mathcal{R}_Q\text{ is transit} \endSb
|c\ci{R,Q}|\cdot|\langle\Dl\ci Q\f, K\ci\Theta(\chi\ci{\C\setminus
R\cci Q}h)\rangle|
\le
$$
$$
\frac{(1+\delta)MA}{(1-\delta)(1-2^{-\frac\e2})}
\sum_{j=1}^4\,\sum_{Q,R\,:\,Q\subset R_j}
\left[\frac{l(Q)}{l(R)}\right]^{\frac\e2}
\sqrt{\frac{\mu(Q)}{\mu(R_j)}}
\|\Dl\ci Q\f\|\ci{L^2(\mu)}
\|\Dl\ci R\psi\|\ci{L^2(\mu)}.
$$
So, it is enough to demonstrate that, say, the matrix $\{T\ci{Q,R}\}$ defined
by
$$
T\ci{Q,R}:=\left[\frac{l(Q)}{l(R)}\right]^{\frac\e2}
\sqrt{\frac{\mu(Q)}{\mu(R_1)}}\qquad\quad(Q\subset R_1),
$$
generates a bounded operator in $l^2$ in the sense that for every two
``sequences''  $\{a\ci Q\}\ci{Q\in\D_1^{tr}}$ and $\{b\ci R\}\ci{R\in\D_2^{tr}}$
of nonnegative numbers, one has
$$
\sum_{Q,R:Q\subset R_1}T\ci{Q,R}a\ci Q b\ci R\le
\frac{1}{ 1-2^{-\frac\e2} }
\Bigl[\sum_{Q}a\ci Q^2\Bigr]^{\frac12}
\Bigl[\sum_{R}b\ci R^2\Bigr]^{\frac12}.
$$
Again
let us ``slice'' the matrix $T\ci{Q,R}$ according to the ratio
$\frac{l(Q)}{l(R)}$. Namely, let
$$
T^{(n)}_{Q,R}=\left\{
\aligned
T\ci{Q,R},&\qquad\text{if } Q\subset R_1,\  l(Q)=2^{-n}l(R);
\\
0,&\qquad\text{otherwise}
\endaligned
\right.
$$
($n=1,2,\dots$).
It is enough to show that for every $n\ge 0$,
$$
\sum_{Q,R}T^{(n)}\ci{Q,R}a\ci Q b\ci R\le
2^{-\frac\e2 n }
\Bigl[\sum_{Q}a\ci Q^2\Bigr]^{\frac12}
\Bigl[\sum_{R}b\ci R^2\Bigr]^{\frac12}.
$$
The matrix $\{T^{(n)}_{Q,R}\}$ has a very good ``block'' structure:
every $a\ci Q$ can interact with {\it only one} variable
$b\ci R$.
So, it is enough to estimate each block separately, i.e., to show that
for every fixed $R\in\D_2^{tr}$,
$$
\sum_{Q:\,Q\subset R_1,\, l(Q)=2^{-n}l(R)}
2^{-\frac\e2 n }
\sqrt{\frac{\mu(Q)}{\mu(R_1)}}
a\ci Q b\ci R\le
2^{-\frac\e2 n }
\Bigl[\sum_{Q}a\ci Q^2\Bigr]^{\frac12}
b\ci R.
$$
But, reducing both parts by the non-essential factor
$2^{-\frac\e2 n }b\ci R$, we see that this estimate is equivalent to the trivial
estimate
$$
\!\!\sum_{Q:\,Q\subset R_1,\, l(Q)=2^{-n}l(R)}
\sqrt{\frac{\mu(Q)}{\mu(R_1)}}
a\ci Q \le
\Bigl[
\sum_{Q:\,Q\subset R_1,\, l(Q)=2^{-n}l(R)}
\frac{\mu(Q)}{\mu(R_1)}\Bigr]^{\frac12}
\Bigl[\sum_{Q}a\ci Q^2\Bigr]^{\frac12}\le
\Bigl[\sum_{Q}a\ci Q^2\Bigr]^{\frac12},
$$
(since squares $Q\in \D_1$  of fixed size do not intersect,
$
\sum_{Q:\,Q\subset R_1,\, l(Q)=2^{-n}l(R)}
\mu(Q)\le \mu(R_1)
$\,).

So, the extension of $\chi\ci{R\cci Q}h$ to the whole $h$ does not cause much
harm, and we get the sum
$$
\sum\Sb Q,R\,:\,Q\subset R,\, l(Q)<2^{-m}l(R),
                    \\  \mathcal{R}_Q\text{ is transit} \endSb
c\ci{R,Q}\langle\Dl\ci Q\f, K\ci\Theta h\rangle
$$
to estimate. Note that the inner product
$\langle\Dl\ci Q\f, K\ci\Theta h\rangle$ {\it does not depend\/} on $R$ at all,
so it seems to be a good idea to sum over $R$ for fixed $Q$ first.
Recalling that
$$
c\ci{R,Q}=\frac{\langle \psi\rangle\ci{R\cci Q}}
{\langle h\rangle\ci{R\cci Q}}-
\frac{\langle \psi\rangle\ci{R}}{\langle h\rangle\ci{R}}
$$
and that $\Lambda\psi=0\Longleftrightarrow \langle \psi\rangle\ci{R^0}=0$, we
conclude that for every $Q\in\D_1^{tr}$ that really appears in the above sum,
$$
\sum\Sb R\,:\,R\supset Q,\, l(R)>2^{m}l(Q),
                    \\  \mathcal{R}_Q\text{ is transit} \endSb
c\ci{R,Q}
= \frac{\langle \psi\rangle\ci{R(Q)}}
{\langle h\rangle\ci{R(Q)}} ,
$$
where $R(Q)$ is the smallest {\it transit} square $R\in \D_2$ containing $Q$
and such that $l(R)\ge 2^ml(Q)$.
So, we obtain the sum
$$
\sum_{Q:\,l(Q)<2^{-m}l(R)}\frac{\langle \psi\rangle\ci{R(Q)}}
{\langle h\rangle\ci{R(Q)}}\langle\Dl\ci Q\f, K\ci\Theta h\rangle
$$
to take care of.

Actually, the range of summation should be
$Q\in \D_1^{tr}$, $Q$ is good (default); there exists a square $R\in
\D_2^{tr}$ such that $l(Q)<2^{-m}l(R)$,
$Q\subset R$ and $R\ci Q$ is transit,
so the last sum we wrote  includes some extra terms compared to the
original one, namely, the terms corresponding to the squares $Q$, for which
$R(Q)=R^0$. But first, we remember that $\langle \psi\rangle\ci{R^0}=0$,
and second, now (but not before!) we are going to put the absolute value
bars around each term, so we may {\it add\/}
as many terms as we want; the point is
not to {\it lose} any of them. In this respect everything is obviously fine.

Clearly, the squares with $\|\Dl\ci Q\f\|\ci{L^2(\mu)}=0$ do not contribute
anything to the sum. Also, since $R(Q)$ is transit,
$|\langle h\rangle\ci{R(Q)}|\ge 1-\delta$. So, we can write
$$
\sum_{Q:\,l(Q)<2^{-m}l(R)}\left|\frac{\langle \psi\rangle\ci{R(Q)}}
{\langle h\rangle\ci{R(Q)}}\langle\Dl\ci Q\f, K\ci\Theta h\rangle\right|\le
$$
$$
\frac{1}{1-\delta}
\sum\Sb Q:\,l(Q)<2^{-m}l(R),
 \\    \|\Dl\cci Q\f\|\cci{L^2(\mu)}>0 \endSb
|\langle \psi\rangle\ci{R(Q)}|\frac
{|\langle\Dl\ci Q\f, K\ci\Theta h\rangle|}{\|\Dl\ci Q\f\|\ci{L^2(\mu)}}
\cdot \|\Dl\ci Q\f\|\ci{L^2(\mu)}\le
$$
$$
\frac{1}{1-\delta}
\Biggl[
\sum\Sb Q:\,l(Q)<2^{-m}l(R),
 \\    \|\Dl\cci Q\f\|\cci{L^2(\mu)}>0 \endSb
|\langle \psi\rangle\ci{R(Q)}|^2
\frac
{|\langle\Dl\ci Q\f, K\ci\Theta h\rangle|^2}{\|\Dl\ci Q\f\|^2\ci{L^2(\mu)}}
\Biggr]^{\frac{1}{2}}\,\,
\Bigl[
\sum_{Q}
\|\Dl\ci Q\f\|^2\ci{L^2(\mu)}
\Bigr]^{\frac{1}{2}}.
$$
The last factor does not exceed $\sqrt2
\|\f\|\ci{L^2(\mu)}$. So, it is sufficient to show that the middle factor
squared is bounded by some constant times $\|\psi\|^2\ci{L^2(\mu)}$.
Switching to the summation over $R$, we see that the middle factor squared
equals
$$
\sum_{R}
|\langle \psi\rangle\ci{R}|^2
\sum_{Q\in\mathcal{F}(R)}
\frac
{|\langle\Dl\ci Q\f, K\ci\Theta h\rangle|^2}{\|\Dl\ci Q\f\|^2\ci{L^2(\mu)}}
=:
\sum_{R}a\ci R
|\langle \psi\rangle\ci{R}|^2,
$$
where
$$
\mathcal{F}(R):=\{Q:\,R(Q)=R,\,\|\Dl\ci Q\f\|\cci{L^2(\mu)}>0\}\,.
$$
So, in order to finish the story with $\s_3^{tr}$, it is enough to show that
the numbers $a\ci R$ satisfy the Carleson condition. Note that for every
$Q\in\mathcal{F}(R)$, one has $Q\subset R$ and that the families $\mathcal{F}(R)$ are
pairwise disjoint (one could say much more, but these two trivial
observations are the only ones that will matter).
Now, for every $S\in\D_2$, we have
$$
\sum_{R:\,R\subset S}a\ci R\le \sum\Sb Q:\,Q\subset S,
\\    \|\Dl\cci Q\f\|\cci{L^2(\mu)}>0 \endSb
\frac
{|\langle\Dl\ci Q\f, K\ci\Theta h\rangle|^2}{\|\Dl\ci Q\f\|^2\ci{L^2(\mu)}}=
\sum\Sb Q:\,Q\subset S
\\    \|\Dl\cci Q\f\|\cci{L^2(\mu)}>0 \endSb
\frac
{|\langle\Dl\ci Q\f,\chi\ci S\cdot
K\ci\Theta h\rangle|^2}{\|\Dl\ci Q\f\|^2\ci{L^2(\mu)}}\le
$$
$$
2\|\chi\ci S\cdot
K\ci\Theta h\|^2\ci{L^2(\mu)}=\int_{S}|K\ci\Theta h|^2\,d\mu\le 2B^2\mu(S)
$$
(because $\Theta\ge\delta\Phi\ci{\D_2}\ge \delta\wt\Phi$),
and we are through.

\bigskip

{\bf XX. Estimation of $\s_1$}

\bigskip

Recall that
$$
\s_1=\sum\Sb Q,R:\,l(Q)\ge 2^{-m}l(R),\\ \dist(Q,R)l(R)\endSb
\langle\Dl\ci Q\f, K\ci\Theta\Dl\ci R\psi\rangle.
$$
We are going to put the absolute value signs around every term and to restore
the symmetry between $Q$ and $R$ (so, we will add the corresponding part from
the sum over pairs $Q,R$, for which $l(Q)\ge l(R)$).
Thus, we have to estimate the sum
$$
\s'_1=\sum\Sb Q\in\D_1^{tr}\,,\,R\in\D_2^{tr}:\,
\\ Q,R\text{ are good}, \\
2^{-m}\le \frac{l(Q)}{l(R)}\le 2^{m},
\\ \dist(Q,R)\le\max\{l(Q),l(R)\}\endSb
|\langle\Dl\ci Q\f, K\ci\Theta\Dl\ci R\psi\rangle|
$$
(now all the conditions for the range of summation are written explicitely).

The key observation about this sum is that every square $Q$ can interact with
not more than $2^{2m}(4\cdot 2^m+1)^2(2m+1)$ squares $R$ and vice versa (the
estimate is quite rough, of course, and is obtained as  follows: one has
$2m+1$ possible values for $l(R)$; once the size
$l(R)\in[2^{-m}l(Q),2^ml(Q)]$ is fixed, the corresponding squares
$R$ are contained in the square of size $(4\cdot 2^m+1)l(Q)$, centered at
the same point as $Q$, are pairwise disjoint, and the area of each of them
is not less than $2^{-2m}l(Q)^2$). Therefore, it is enough to show that
for some large constant $U>0$, not depending on $\f$, $\psi$ and $\Theta$, one
has
$$
|\langle\Dl\ci Q\f, K\ci\Theta\Dl\ci R\psi\rangle|
\le U \|\Dl\ci Q\f\|\ci{L^2(\mu)}
\|\Dl\ci R\psi\|\ci{L^2(\mu)},
$$
provided that  $Q\in\D_1^{tr}$, $R\in\D_2^{tr}$,
$Q,R$ are good,
$2^{-m}\le \frac{l(Q)}{l(R)}\le 2^{m}$ and
$\dist(Q,R)\le\max\{l(Q),l(R)\}$.

\bigskip

{\bf XXI. Negligible contours}

\bigskip

Let $G$  be a contour on the complex plane $\C$. Let $\wt M$ be some large
positive number. We will call $G$ {\it negligible} (the full name should be $\wt
M$-negligible with respect to the measure $\mu$), if for every $r>0$,
$$
\mu\{x\in\C\,:\,\dist(x,G)\le r\}\le\wt M r.
$$
{\bf Lemma:}
Let $G$ be a negligible contour splitting the complex plane $\C$ into two
(open) regions $\Om_1$ and $\Om_2$. Then for any two functions
$\eta_1,\eta_2\in L^2(\mu)$ such that $\eta_j$ vanishes outside
$\Om_j\cup G$, one has
$$
|\langle\eta_1,K\ci\Theta\eta_2\rangle|\le 4\wt M
\|\eta_1\|\ci{L^2(\mu)}
\|\eta_2\|\ci{L^2(\mu)}.
$$
{\bf Proof:}
Note that the condition that $G$ is negligible immediately implies that
$\mu(G)=0$. So, we may assume that $\eta_j$ vanishes outside $\Om_j$.
We have
$$
|\langle\eta_1,K\ci\Theta\eta_2\rangle|\le
\iint|k\ci\Theta(x_1,x_2)|\cdot|\eta_1(x_1)|\cdot|\eta_2(x_2)|\,
d\mu(x_1)\,d\mu(x_2).
$$
Clearly, the integrand can be non-zero only if $x_1\in\Om_1$ and
$x_2\in\Om_2$. According to the Schur test (full $L^2$-version), it is enough
to find a function $\la:\C\setminus G\rightarrow (0,+\infty)$, such that
$$
\int_{\Om_1}|k\ci\Theta(x_1,x_2)|\la(x_1)\,
d\mu(x_1)\le 4\wt M\la(x_2)\qquad\text{ for every }x_2\in\Om_2,
$$
and vice versa, i.e.,
$$
\int_{\Om_2}|k\ci\Theta(x_1,x_2)|\la(x_2)\,
d\mu(x_2)\le 4\wt M\la(x_1)\qquad\text{ for every }x_1\in\Om_1.
$$
We will check that these properties hold for
$$
\la(x)=\frac{1}{\sqrt{\dist(x,G)}}\ .
$$
Indeed, for $x_1\in\Om_1$ and
$x_2\in\Om_2$, we have
$$
|k\ci\Theta(x_1,x_2)|\le\frac{1}{|x_1-x_2|}\le
\frac{1}{\max\{\dist(x_1,G),\dist(x_2, G) \} }\ .
$$
Thus, according to the Comparison lemma,
$$
\int_{\Om_1}|k\ci\Theta(x_1,x_2)|\la(x_1)\,
d\mu(x_1)\le
\int_{\Om_1}
\frac{1}{\max\{\dist(x_1,G),\dist(x_2, G) \} }\,
\frac{1}{\sqrt{\dist(x_1,G)}}
d\mu(x_1)
\le
$$
$$
\wt M
\int_{0}^{+\infty}
\frac{1}{\max\{t,\dist(x_2, G) \} }\,
\frac{1}{\sqrt{t}}
dt=\frac{4\wt M}{\sqrt{\dist(x_2,G)}}=4\wt M\la(x_2).
$$
Now observe that
$$
\langle\Dl\ci Q\f, K\ci\Theta\Dl\ci R\psi\rangle
=\sum_{i,j=1}^4
\langle\f\ci Q^{(i)}, K\ci\Theta\psi\ci R^{(j)}\rangle,
$$
where $\f\ci Q^{(i)}:=\chi\ci{Q_i}\Dl\ci Q\f$, and
$\psi\ci R^{(j)}:=\chi\ci{R_j}\Dl\ci R\psi$.

Assume that the boundaries of all the subsquares $Q_i$ and $R_j$ are $\wt
M$-negligible contours.
Then it makes sense to write
$$
\langle\f\ci Q^{(i)}, K\ci\Theta\psi\ci R^{(j)}\rangle
=
$$
$$
\langle\chi\ci{Q_i\setminus R_j}\cdot\f\ci Q^{(i)}, K\ci\Theta\psi\ci R^{(j)}\rangle
+
\langle\chi\ci{Q_i\cap R_j}\cdot\f\ci Q^{(i)}, K\ci\Theta
(\chi\ci{R_j\setminus Q_i}\cdot\psi\ci R^{(j)})\rangle
+
$$
$$
\langle\chi\ci{Q_i\cap R_j}\cdot\f\ci Q^{(i)}, K\ci\Theta
(\chi\ci{Q_i\cap R_j}\cdot\psi\ci R^{(j)})\rangle.
$$
In the first two terms the supports of the functions are separated by
negligible contours ($\partial R_j$ and $\partial Q_i$, respectively). So,
the corresponding inner products are bounded by
$$
4\wt M\|\f\ci Q^{(i)}\|\ci{L^2(\mu)}
\|\psi\ci R^{(j)}\|\ci{L^2(\mu)}\le
4\wt M\|\Dl\ci Q\f\|\ci{L^2(\mu)}
\|\Dl\ci R\psi\|\ci{L^2(\mu)}\,.
$$

As to the inner product
$
\langle\chi\ci{Q_i\cap R_j}\cdot\f\ci Q^{(i)}, K\ci\Theta
(\chi\ci{Q_i\cap R_j}\cdot\psi\ci R^{(j)})\rangle
$, there are two possibilities:

\bigskip

{\bf Case 1: one of the squares (say, $Q_i$) is terminal}

Then we have the estimate
$$
|k\ci\Theta(x,y)|\le
\frac{1}{\delta\max\{\dist(x,\partial Q_i),\dist(y, \partial Q_i) \} }
$$
for all $x,y\in Q_i\cap R_j$ and, repeating our above reasoning with the Schur
test, we obtain
$$
|\langle\chi\ci{Q_i\cap R_j}\cdot\f\ci Q^{(i)}, K\ci\Theta
(\chi\ci{Q_i\cap R_j}\cdot\psi\ci R^{(j)})\rangle|
\le
$$
$$
\frac{4\wt M}{\delta}\|\f\ci Q^{(i)}\|\ci{L^2(\mu)}
\|\psi\ci R^{(j)}\|\ci{L^2(\mu)}\le
\frac{4\wt M}{\delta}\|\Dl\ci Q\f\|\ci{L^2(\mu)}
\|\Dl\ci R\psi\|\ci{L^2(\mu)}.
$$

{\bf Case 2: both squares $Q_i$ and $R_j$ are transit}

Then both  functions
$\chi\ci{Q_i\cap R_j}\cdot\f\ci Q^{(i)}$ and
$\chi\ci{Q_i\cap R_j}\cdot\psi\ci R^{(j)}$ are constant multiples of {\it the
same} function $\eta:=\chi\ci{Q_i\cap R_j}\cdot h$. But the kernel $k\ci\Theta$
is antisymmetric, and therefore $\langle\eta,K\ci\Theta\eta\rangle=0$.

\medskip

What if the boundary of some square $Q_i$ (or $R_j$) is not negligible? We do
not know how to get a good estimate in this case; instead, we will try to rule it
out by declaring the corresponding squares bad. But we should be very careful
here: the temptation to declare a square $Q\in\D_1^{tr}$ bad if $\partial Q$
is not negligible should be severely suppressed, because, as we remember,
``badness'' of the square $Q$ should depend rather on $\D_2$, than on $Q$
itself. So, we are going to use a little bit less straightforward definition.

Namely, we will call a transit square $Q\in\D_1$ bad if there exists a
transit square $R\in \D_2$ such that $2^{-m}l(Q)\le l(R)\le 2^ml(Q)$,
$\dist (R,Q)\le 2^ml(Q)$ and for at least one $j=1,2,3,4$, the boundary
$\partial R_j$ is not $\wt M$-negligible (we do not care about the terminal
squares, so let them all be ``good by the definition''). Then for every pair
of squares $Q,R$ appearing in the sum $\s'_1$, the assumption that $Q$ is
good allows to conclude that all the four subsquares $R_j$ of the square $R$
are negligible and vice versa! Now it remains only to show that we can choose
$m$ and $\wt M$ (in this order) so that $P\ci{\D_2}\{Q\text{ is bad}\}\le\delta$
for every $Q\in\D_1^{tr}$.

\bigskip

{\bf XXII. Estimation of probability}

\bigskip

Let $Q\in\D_1^{tr}$. Consider the ``extended lattice''
$$
\wt\D_2=\wt\D_2(\om)=\left\{\om+\left[\tfrac{j}{2^n},\,\tfrac{j+1}{2^n}
\right)\times \left[\tfrac{k}{2^n},\,\tfrac{k+1}{2^n}
\right)\,:\, j,k,n\in \Bbb Z, n\ge 1
\right\}.
$$
Clearly, $\wt D_2$ contains every square $R\in\D_2$ of size $\frac12$ or
less.
Note that when $\om$ runs over $\left[-\tfrac{1}{4},\,\tfrac{1}{4}
\right)\times \left[-\tfrac{1}{4},\,\tfrac{1}{4}
\right)$, the lattice $\wt D_2$ runs over its whole period.

Starting now, we will declare a square $Q\in\D_1^{tr}$ bad if either
\medskip\it

{\bf 1)} there exists a square $R\in\wt\D_2$ such that
 $\dist(Q,\partial R)\le 16l(Q)^{\al}l(R)^{1-\al}$ 
and $l(R)\ge
2^ml(Q)$,
\smallskip\rm
or
\smallskip\it
{\bf 2)} there exists a square $R\in\wt\D_2$ such that $R\subset (4\cdot
2^m+1)Q$, $l(R)\ge 2^{-(m+1)}l(Q)$ and $\partial R$ is not $\wt
M$-negligible.
\medskip\rm
We leave it to the reader to check that every square $Q$ bad in the sense of
Section XV, XVI or XXI is bad according to this new definition as well.

{\bf Choice of $m$}

Fix $k\ge m$.
Let us estimate the probability that there exists a square $R\in\wt\D_2$ of
size $l(R)=2^kl(Q)$  such that
$\dist(Q,\partial R)\le 16l(Q)^{\al}l(R)^{1-\al}$.
Since the lattice $\wt\D_2$ runs over its whole period, we can find this
probability exactly: it equals to the ratio of the area of the dashed rim
on Picture 3 to the area of the whole square with side $2^kl(Q)$ (just
look at where the center $x\ci Q$ should lie with respect to the lattice
$\D_2$). Observing that
$$
16l(Q)^{\al}l(R)^{1-\al}+\tfrac{l(Q)}{2}\le 17l(Q)^{\al}l(R)^{1-\al},
$$
we conclude that this ratio is less than $68\left[\frac{l(Q)}{l(R)}
\right]^{\al}=68\cdot 2^{-k\al}$.

Therefore the probability that the square $Q$ is bad according to the first
part of our definition does not exceed
$$
68\sum_{k=m}^\infty2^{-k\al}=\frac{68\cdot 2^{-m\al}}{1-2^{-\al}}\le
\frac{\delta}{3},
$$
provided that $m$ is taken large enough.

\bigskip

{\bf Choice of $\wt M$}

\bigskip

Now let us look at how large the probability that $Q$ is bad according to the
second part of our definition may be.
Recall that $\partial R$ is $\wt M$-negligible if
$\mu\{x\in\C\,:\,\dist(x,\partial R)\le r\}\le \wt M r$ for all $r>0$.
Note first of all, that we do not have any trouble with $r\ge l(Q)$. Indeed,
since $R\subset(4\cdot 2^m+1)Q$, we have
$$
\{x\in\C\,:\,\dist(x,\partial R)\le r\}\subset
B(x\ci Q,\,(4\cdot 2^m+1)l(Q)+r))\subset B(x\ci Q,(4\cdot 2^m+2)r).
$$
But $\mu(B(x\ci Q,(4\cdot 2^m+2)r)\le (4\cdot 2^m+2)Mr$, because
$Q$ is a transit square, $r\ge l(Q)$ and therefore
$\mathcal{R}(x\ci Q)\le l(Q)\le r<(4\cdot 2^m+2)r$.
So, everything is okay with such $r$, provided that $\wt M\ge (4\cdot
2^m+2)M$.

Now observe that for $r<l(Q)$ we have
$$
\{x\in\C\,:\,\dist(x,\partial R)\le r\}\subset
B(x\ci Q,(4\cdot 2^m+2)l(Q)\,).
$$
So, the part of the measure $\mu$ that lies outside the disk
$B(x\ci Q,(4\cdot 2^m+2)l(Q)\,)$ does not matter and we can replace the
whole measure $\mu$ by its restriction $\wt\mu$ to this disk, defined as
$$
\wt\mu(E):=\mu(E\cap B(x\ci Q,(4\cdot 2^m+2)l(Q)\,)\,).
$$
Though we do not know {\it much} about $\wt\mu$, there is one thing we can
say for certain:
$$
\wt\mu(\C)\le (4\cdot 2^m+2)Ml(Q);
$$
and this will be enough for us.

Consider the grid $\mathcal{L}=\mathcal{L}(\om)$ consisting of all vertical lines
serving as boundaries of squares in $\wt\D_2$ of size $2^{-(m+1)}l(Q)$.
We are going to show that if $\wt M$ is sufficiently large, then,
with probability $1-\frac{\delta}{3}$ or more, this entire
grid is $\frac{\wt M}{2}$-negligible with respect to the measure $\wt\mu$.
Of course (together with the same estimate for horizontal lines), this will imply that 
the probability that the square $Q$ is bad
according to the second part of our definition does not exceed
$\frac{2\delta}{3}$, finishing the story.

Note that the grid $\mathcal{L}$ runs (several times)
over its whole period when $\om$ runs over
$\left[-\tfrac{1}{4},\,\tfrac{1}{4}
\right)\times \left[-\tfrac{1}{4},\,\tfrac{1}{4}
\right)$. So, we can change the random parameter $\om$ to another
random parameter $\tau\in[0,2^{-(m+1)}l(Q)\,)$ (which is just
the real part of $\om\mod 2^{-(m+1)}l(Q)$, of course) and reformulate our problem as the
following: we should demonstrate that the one-dimensional Lebesgue measure of
such $\tau\in[0,2^{-(m+1)}l(Q)\,)$ that the grid $\mathcal{L}(\tau)$ consisting of
all vertical lines intersecting the real axis at the points of the kind
$\tau+\frac{k}{2^{m+1}}$, $k\in \Bbb Z$, is not $\frac{\wt M}{2}$-negligible
with respect to the measure $\wt\mu$, does not exceed $\frac{\delta}{3}
2^{-(m+1)}l(Q)$.

Consider the $2^{-(m+1)}l(Q)$-periodic sweeping $\nu$ of the measure $\wt\mu$,
i.e. the measure defined on Borel subsets $E$ of the real line $\Bbb R$ by
$$
\nu(E)=\wt\mu\Bigl(
\bigcup_{k\in\Bbb Z}\bigl(k\cdot2^{-(m+1)}l(Q)+E\bigr)\times\Bbb R
\Bigr).
$$
Note that $\mathcal{L}(\tau)$ is not $\frac{\wt M}{2}$-negligible if and only if
$\mathcal{M}\nu(\tau)>\frac{\wt M}{2}$, where
$$
\mathcal{M}\nu(\tau)=\sup_{r>0}\frac{\nu([\tau-r,\tau+r])}{2r}
$$
is the Hardy-Littlewood maximal function. But the standard estimate for the
maximal function of a periodic measure yields
$$
m_1\{\tau\in [0,2^{-(m+1)}l(Q))\,:\,\mathcal{M}\nu(\tau)>\frac{\wt M}{2} \}
\le \frac{4\nu([0,2^{-(m+1)}l(Q))\,)}{\wt M}
=
$$
$$
\frac{4\wt\mu(\C)}{\wt M}
\le\frac{4(4\cdot 2^m+2)M}{\wt M}l(Q).
$$
So, if $\wt M\ge 12\delta^{-1}(4\cdot 2^m+2)M$, we are okay.

\bigskip

{\bf XXIII. Quantitative pulling ourselves up by the hair}

\bigskip

We are going to present the succession of ``fancy" $Tb$ theorems in the nonhomogeneous setting.

The first one is the least ``fancy" because $b$ will be accretive in it, but it solves a problem of P. Mattila about an analytic 
characterization of Besicovitch rectifiable sets. 

The middle one is the theorem proved in the previous sections; it gives an alternative proof of the result of Guy David [D1], 
thus solving the analytic part of Vitushkin's conjecture. 

The last one---and the most difficult---gives a quantitative information in the solution of Vitushkin's conjecture. 
Namely, given a set $E$ of positive analytic capacity $\gamma$ and length $M$, this last theorem allows us to say  
quantitatively what
portion of the length is rectifiable, and ``how" rectifiable it is.

\bigskip

In what follows $\mu$ is a positive measure on $\C$ satisfying the following {\it non-uniform linear growth condition}:

$$
\limsup_{r\rightarrow 0}\frac{\mu(B(x,r))}{r} < \infty\,\,\,\text{for}\,\,\mu\,\,\text{a.e}\,\,x\,.
$$

The truncated Cauchy integral is

$$
(C_{\mu}^{\e}b)(\zeta)= \int_{|z-\zeta|>\e}\frac{b(z)\,d\,\mu(z)}{\zeta -z}\,.
$$

The maximal Cauchy integral is 

$$
(C_{\mu}^{*}b)(\zeta) = \sup_{\e>0}|(C^{\mu}_{\e}b)(\zeta)|\,.
$$

Recall that for any $1$-Lipschitz function $\Phi$ on $\C$ the following Calder\'on-Zygmund kernel was introduced

$$
k_{\Phi}(x,y) = \frac{\overline{x-y}}{|x-y|^2 + \Phi(x)\Phi(y)}
$$
and let $K_{\Phi}$ be the canonical Calder\'on-Zygmund operator with this antisymmetric kernel. 
Recall that
$$
k_{\Phi}(x,y) \leq \min [\frac1{\Phi(x)},\frac1{\Phi(y)}]\,.
$$

Consider another truncation of the Cauchy integral:

$$
(C_{\Phi}b)(\zeta) = \int_{|z-\zeta|\geq\Phi(\zeta)}\frac{b(z)\,d\,\mu(z)}{\zeta-z}\,.
$$

As usual $M_1$ denotes the following maximal function 

$$
(M_1 f)(\zeta) = \sup_{r>0}\frac1{r}\int_{B(\zeta,r)} |f(z)|\,d\,\mu(z)\,.
$$

\bigskip 

{\bf Lemma 1:}
$$|(K_{\Phi} f)(x) - (C_{\Phi} f)(x)| \le A (M_1 f) (x)\,.$$

{\bf Proof}
Fix $x$ and consider the absolute value of the difference of the kernels. For $y\in B(x, \Phi(x))$ it is at most $\frac1{\Phi(x)}$. 
For $y$ such that $|y-x|\ge \Phi(x)$ it is 
$$
\le \frac{(\Phi(x)\Phi(y)|x-y|}{[|x-y|^2 +\Phi(x)\Phi(y)]|x-y|^2}\le |k_{\Phi}|\frac{\Phi(x)\Phi(y)}{|x-y|^2}\le 
\frac{\Phi(x)}{|x-y|^2}\,.
$$
Splitting $\{y: |y-x|\ge \Phi(x)\}$ into annuli $\{y: 2^{k+1}\Phi(x)> |y-x|\ge 2^k\Phi(x)\}$ finishes the proof.

\bigskip

Recall that we have assumption of non-uniform linear growth on 
$\mu$.

Let us also normalize $\mu$ and think (if otherwise  not stated) that $\|\mu\|=1$. Recall that $M$-non-Ahlfors disc is a $B(x,R)$
such that

$$
\mu(B(x,R)) \ge MR,\,\, x\in\supp\mu\,.
$$
In this case the point $x$ is called  an $M$-non-Ahlfors point.

{\bf Lemma 2:}
There exists $\e = \e(M), \, \e\rightarrow 0$ if $M\rightarrow \infty$, such that the union of all $M$-non-Ahlfors discs has $\mu$-measure 
at most $\e$.

{\bf Proof.}
It follows from non-uniform linear growth condition that 

$$
\mu\{x\in\supp\mu: \sup_{r}\frac{\mu(B(x,r))}{r} \ge \sqrt{M}\} =\delta(M)\rightarrow 0,\,\,\text{when}\,\, M\rightarrow \infty\,.
$$

Denote this set by $G_{M}$. For $x \in \supp\mu \setminus G_M$ we choose the maximal $M$-non-Ahlfors disc centered at $x$ (if any). 
Their union will be called $O$. By Vitali's lemma, $O$ is covered by $\cup B(x_j, 5r_j)$, where $x_j\in \supp\mu\setminus G_M$, and $B(x_j,
r_j)$ are disjoint and $M$-non-Ahlfors. Thus,

$$
\Sigma \,r_i \le \frac1{M}\Sigma \,\mu(B(x_i,r_i)) \le\frac1{M}
$$

On the other hand, $\mu (B(x_j, 5r_j))\le 5\sqrt{M}\, r_j$. Thus,

$$
\mu(O) \le \Sigma\,\mu(B(x_j, 5r_j)) \le 5\sqrt{M} \Sigma \,r_j \le \frac5{\sqrt{M}}
$$
$M$-non-Ahlfors points can be only in $O\cup G_M$. So we see that

$$
\mu(\text{M-non-Ahlfors points}) \rightarrow 0,\,\,\text{when} \,\, M\rightarrow \infty
$$

But we want a bit more---the smallness of measure of the union of all $M$-non-Ahlfors discs. To get this, consider points in
$G_M$,  and consider  maximal $M$-non-Ahlfors disc centered at each of them. Call their union $G$. The set $G$ is covered by
$\cup_j B(y_j, 5R_j)$, where $B(y_j, R_j)$ are disjoint $M$-non-Ahlfors discs. Consider  $y\in B(y_j, 5R_j)$. Then
$
\frac{\mu(B(y, 10R_j))}{10R_j}\ge \frac{\mu(B(y_j, R_j))}{10R_j}\ge \frac{\sqrt{M}}{10}
$. In our notations this means that $y\in G_{M/100}$. Thus, $G\subset G_{M/100}$. So $\mu(G)$ is small if $M$ is large.
The lemma is proved.

{\bf Lemma 3:}
Let $\|\mu\|=1$, let $\mu$ be a positive measure with non-uniform linear growth, and let 
$H=H_M$ be the union of all $M$-non-Ahlfors discs. Let $\Phi$ be a $1$-Lipschitz function such that $\Phi(x)\ge \dist(x,\C\setminus H)$.
Then $K_{\Phi}$ and $C_{\Phi}$ are bounded or unbounded simultaneously on $L^2(\mu)$.

{\bf Proof.}
In Lemma 1 we saw that $|(K_{\Phi} - C_{\Phi}) (f)(x)| \le A (M_1 f)(x)$. Actually, the proof says more, namely

$$
|(K_{\Phi} - C_{\Phi}) (f)(x)| \le A (M_{1,\Phi} f)(x) := A\sup_{r\ge \Phi(x)} \frac1r\int_{B(x,r)}|f(y)|\,d\,\mu(y)\,.
$$
But for $r\ge \Phi(x)$ we have $\mu (B(x,r)) \le Mr$. Therefore,

$$
 (M_{1,\Phi} f)(x)  \le M\sup\frac1{\mu(B(x,r))} \int_{B(x,r)} |f(y)|\,d\,\mu(y)\,.
$$
It is well-known that this maximal operator is bounded in $L^2(\mu)$.
Lemma is proved.

\bigskip

Now we are ready to present several conditions for $K_{\Phi}$ ($\Phi$ is a $1$-Lipschitz function) to be bounded on $L^2(\mu)$. 
While doing that we are interested in such $\Phi$'s that $F_{\Phi}:=\{x\in\C: \Phi(x)=0\}$ has positive measure (or, if circumstances
permit, even measure close to $1$). This interest is easy to explain: for such $\Phi$ we have
$$
\quad\qquad\quad\qquad\qquad\qquad (K_{\Phi}\,f,\,g)=(Cf,g)\quad\qquad\quad\qquad\quad\qquad\quad\qquad\quad\qquad\quad\qquad (*)
$$
for $f,g$ supported on $F_{\Phi}$. And after all, we are interested in estimates of the Cauchy operator $C$.

\bigskip

{\bf Theorem 1:} 
Let $\mu$ be a measure with non-uniform linear growth, let $H_M$ be the union of all $M$-non-Ahlfors discs, 
$$
\Phi(x) \ge \dist(x,\C\setminus H_M)
$$
and let $\Phi$ be a $1$-Lipschitz function. Consider 
$(K_{\Phi,\e}f)(x) := \int_{|y-x|\ge\e}k_{\Phi}(x,y) f(y)\,d\,\mu(y)$.
If there exists a constant $B$ such that 
$$
(K^*_{\Phi}1)(x) := \sup_{\e>0} |(K_{\Phi,\e}1)(x)| \le B <\infty \,\,\,\text{for} \,\,\mu\,\,\text{a. e.} \,\,x
$$
then 
$$
\|K_{\Phi}\|_{L^2(\mu)\rightarrow L^2(\mu)} \le ABM\,.
$$

\bigskip

An assumption on {\it the maximal singular function} $(K^*_{\Phi}1)$ can be conveniently modified. Let us consider 
the following  assumption of {\it a.e.  finiteness of the maximal singular function}:

$$
(K^*_{\Phi}1)(x) <\infty  \,\,\,\text{for} \,\,\mu\,\,\text{a. e.} \,\,x\,.
$$

Fix a large $M$ and $L>100 M$. Fix a bounded measurable function $b$. 
If $x$ is such that $(K^*_{\Phi}b)(x)>L$, then there exists a maximal $\e_0(x)$ such that $|(K_{\Phi,\e_0}b)(x)|\ge L$ (the
function
$\e\rightarrow (K_{\Phi,\e_0}b)(x)$ is right continuous). Consider
$$
G_L(b):= \cup_{x\in\supp\mu} B(x, 2\e_0(x))\,.
$$

{\bf Lemma 4:}
Let us  assume the a.e. finiteness of the maximal singular function. Then $\mu(G_L(b)\setminus H_M)\rightarrow 0$ 
if $L\rightarrow
\infty$.

\bigskip

Remind that Calder\'on-Zygmund constants of kernel $k_{\Phi}$ are bounded by $C$.

\bigskip

{\bf Proof.}
Let $y\in B(x, 2\e_0(x)$. Let us prove first that
$$
K_{\Phi}^*1 (y) \geq L- ACM
$$
for an absolute constant $A$. In fact, let us consider two cases: a) $\Phi(x) 
\geq \frac15 \e_0(x)$, b) $\Phi(x) < \frac15 \e_0(x)$. In the
first case let $\e =20 \Phi(x)$. Then 
$$
|K_{\Phi}^{\e}(y) - K_{\Phi}^{\e_0(x)}(x)| \leq \int_{z:|z-y|\geq \e}|k_{\Phi}(y,z)-
k_{\Phi}(x,z)|\
,d\mu(z) + \int_{z: |z-y|\leq 20 \Phi(x)}|k_{\Phi}(x,z)|\,d\mu(z)\,.
$$ 
The first integral can be estimated as usual using the
Calder\'on-Zygmund property of the kernel $k_{\Phi}$ and the fact that all ``large" disks centered at $y$ are contained in 
discs centered at $x$ of ``almost" the same radii. These radii will be larger than $\Phi(x)$, and, hence, they will be $M$-Ahlfors.
The second integral is bounded by
$AM$ because  $k_{\Phi}(x,z)\leq \frac1{\Phi(x)}$ and $\mu (B(y, 20 \Phi(x))\leq \mu (B(x, 40 \Phi(x))\leq 40M \Phi(x)$ 
(the first
inequality holds because we are in the first case). 

Let us consider case b) now. Put  
$\e = 4 \e_0(x)$. Then 
$$
K_{\Phi}^{\e}(y) - K_{\Phi}^{\e_0(x)}(x)| 
\leq \int_{z:|z-y|\geq \e}|k_{\Phi}(y,z)-k_{\Phi}(x,z)|\,d\mu
(z) +
$$
$$
 \int_{z: |z-y|\leq 4 \e_0(x)}|k_{\Phi}(x,z)|\,d\mu(z)\,.
$$
The first integral can be estimated exactly as in the case a). The
second integral is bounded by $\frac{C}{\e_0(x)} \,\mu(B(x, 6\e_0(x)))$, where $C$ is the constant from 
Calder\'on-Zygmund properties
of our kernel. The disc $B(x, 6\e_0(x))$ is $M$-Ahlfors because we are in case b). Thus the second integral is also bounded by
$ACM$.

Now it is clear that the assumption

$$
K_{\Phi}^*1(x) <\infty \,\,\,\text{for}\,\,\mu\,\,\text{a.e}\,\,x
$$
implies that 

$$
\mu (G_L(1))\rightarrow 0\,\,\,\text{when}\,\, L\rightarrow \infty\,.
$$
Lemma 4 is proved. 

\bigskip

Recall that for a given $M$, $H_M$ denotes the union of all $M$-non-Ahlfors disks.

\bigskip

{\bf Theorem 1a:}
Let $\mu$ satisfy the non-uniform linear growth condition, let $\Phi$ be a $1$-Lipschitz function such that $\Phi(x) \ge \dist(x,
\C\setminus H_M)$, and assume the a.e. finiteness of the maximal singular  function $(K^*_{\Phi}1)(x)$. Let $\Psi$ 
be a
$1$-Lipschitz function such that
$
\Psi(x)\ge\max [\dist(x,
\C\setminus(G_L(1)), \Phi(x)]$. Then

\noindent 1) $(K^*_{\Psi}1)(x) \le AC(L+M)$ uniformly, and 

\noindent 2) $\|K_{\Psi}\|_{L^2(\mu)\rightarrow L^2(\mu)} \le AC(L+M)M$.

\bigskip

{\bf Remark.} In the first claim of Theorem 1a one can replace $1$ by any bounded function $b, \|b\|_{\infty}\leq 1$.

\bigskip

{\bf Proof.}
The second claim of the Theorem follows from the first claim and from Theorem 1. 
The first claim is a simple calculation using Lemma 1. Let us do it for the sake of completeness.

Let $C_{\Psi}^{\e} f(x):= \int_{y: |y-x| \ge \max [\e, \Psi(x)]} k(x,y)\,d\mu(y)$.

\noindent Step 1. For any $1$-Lipschitz $\Psi \ge \Phi$ and any $\e$,
$$
|K_{\Psi}^{\e}1(x) -C_{\Psi}^{\e}1(x)|\le A M_{1,\Psi}1(x) \le  A M_{1,\Psi}1(x)\le AM\,.
$$
In fact, if $\e\le \Psi(x)$, then 
$$
|K_{\Psi}^{\e}1(x) -C_{\Psi}^{\e}1(x)|\le \int_{\e\le |y-x|\le \Psi(x)}
\frac{d\mu(y)}{\Psi(x)} + \int_{y: |y-x| >\Psi(x)}\frac{\Psi(x)\Psi(y)|x-y|}{(|x-y|^2 + \Psi(x)\Psi(y))|x-y|^2}\,.
$$
The first term is bounded by $\frac{\mu(B(x, \Psi(x))}{\Psi(x)}\le M$. 

The second term can be estimated precisely as in Lemma 1 if we use
that $\Psi \ge \Phi$. So it is also bounded by $AM$.
If now $\e > \Psi(x)$, then only the second term will appear. We are done with the first step.

\noindent Step 2. Recall that $\e_0(x) = \max\{\e : |K_{\Phi}^{\e}1(x)|\le L\}$. Fix $x_0$ and let $\e \le \e_0(x_0)$. 
Then $\Psi(x_0) \ge
2\e_0(x_0) >\e$. Then $ |K_{\Psi}^{\e}1(x_0)|\le \int_{\e \le |y-x| \le 2\e_0}|k_{\Psi}(x_0, y)\,d\mu(y) + 
|K_{\Psi}^{2\e_0}1(x_0)|$.
The first term is bounded by $\frac{\mu(B(x_0, 2\e_0)}{\Psi(x_0)}\le \frac{\mu(B(x_0, \Psi(x_0))}{\Psi(x_0)} \le M$ 
since $\Psi\ge \Phi$.
Using step 1 we can reduce the estimate of the second term to the estimate of $|C_{\Psi}^{2\e_0}1(x_0)|$ (with the error 
bounded by $AM$).
Let us now use the fact that $\Psi(x_0) \ge
2\e_0(x_0),\, \Psi(x_0) \ge \Phi(x_0)$. This means that $C_{\Psi}^{2\e_0}1(x_0) = C_{\Phi}^{\Psi(x_0)}1(x_0)$. 
By another application of
step 1 we can see that the last quantity is within $AM$ of $K_{\Phi}^{\Psi(x_0)}1(x_0)$. The absolute value of this  
expression is bounded by
$L$ by the definition of $\e_0$ and the fact that $\Psi(x_0) \ge 2\e_0$. In particular, our second term is bounded by $L+AM$.

\noindent Step 3. $\e_0(x_0) < \e \le \Psi(x_0)$. Then $|K_{\Psi}^{\e}1(x_0) -K_{\Phi}^{\e}1(x_0)|\le \int_{\e \le |y-x_0| \le 
\Psi(x_0)}
|k_{\Psi} (x_0,y)|\,d\mu(y) + |\int_{\e \le |y-x_0| \le \Psi(x_0)} k_{\Phi}
(x_0,y)\,d\mu(y)| + \int_{y: |y-x_0| \ge \Psi(x_0)} |k_{\Phi}(x_0,y) - k_{\Psi}(x_0,y)|\,d\mu(y)$. The first term is bounded by
$\frac{\mu(B(x_0, \Psi(x_0))}{\Psi(x_0)} \le M$. The second term is bounded by $|K_{\Phi}^{\e}1(x_0)|+
|K_{\Phi}^{\Psi(x_0)}1(x_0)|\le 2L$, because $\e_0(x_0) < \e \le \Psi(x_0)$, just by the definition of $\e_0$. 
The third term
can be estimated precisely as in Lemma 1 if we notice that the integrand is bounded by $|k_{\Phi}(x_0,y) - k(x_0,y)| + |k(x_0,y) -
k_{\Psi}(x_0,y)|\le \frac{2\Psi(x_0)}{|x-x_0|^2}$. The integral then is bounded by $AM$.

\noindent Step 4. $\e > \Psi(x_0)$. We use the first step to write $|K_{\Psi}^{\e}1(x_0) - C_{\Psi}^{\e}(x_0)| \le AM$ and also 
$|K_{\Psi}^{\e}1(x_0) - C_{\Psi}^{\e}(x_0)| \le AM$. Therefore, we are left to estimate $|C_{\Phi}^{\e}1(x_0) -
C_{\Psi}^{\e}(x_0)|$. But this quantity vanishes because $\e > \Psi(x_0)\ge \Phi(x_0)$.

\noindent The first claim of Theorem 1a is completely proved. We have already made a remark that the second claim follows 
from the
first one combined with Theorem 1.

\bigskip

Before proving Theorem 1, we would like to present its beautiful application found by
Xavier Tolsa [XT2].

\bigskip

Recall that $R(x,y,z)$ denotes the radius of the circle passing through $x,y,z\in\C$. We will call a measure $\mu\ge 0$ 
on $\C$ Besicovitch-Melnikov-Verdera rectifiable if $\mu =\sum_{n=0}^{\infty} \mu|E_n$, $E_n, n=1,2,3,...$ are compact sets, and 

$$
c^2(\mu|E_n):= \iiint_{E_n^3} R^{-2}(x,y,z) \,d\,\mu(x) \,d\,\mu(y)\,d\,\mu(z)<\infty, \,\, n=1,2,3,...; \mu(E_0) =0\,.
$$
The curvature $c^2(\mu)$ was widely used by Melnikov and Verdera (see, for example [MV]). The name is natural because if 
$\mu=\mathcal{H}^1|E$, 
$E$ being a compact set, then $\mu$ turns out to be a Besicovitch-Melnikov-Verdera rectifiable measure if and only if $E$ 
is a Besicovitch
rectifiable set. This is a difficult geometric result proved by G. David and J.-C. L\'eger. This result becomes especially 
difficult because
of the nonhomogeneity of the measure $\mu$, namely because $\liminf_{r\rightarrow 0}\frac{\mu(B(x,r))}{r} $ may apriori easily 
become $0$. 

In his paper [XT2] Xavier Tolsa gave the following application of Theorem 1a. We use the notations of Theorem 1.

\bigskip

{\bf Theorem (Xavier Tolsa):}
If $\mu$ satisfies non-uniform linear growth condition and if for any $M$ the assumption of a.e.  finiteness of the 
maximal singular  function 
$(K^*_{\Phi}1)(x)$ holds for $\Phi(x):= \dist(x,\C\setminus H_M)$, then $\mu$ is Besicovitch-Melnikov-Verdera rectifiable. 
If in addition
for $\mu$ a.e. $x$, $\limsup_{r\rightarrow 0}\frac{\mu(B(x,r))}{r}>0$, then $\supp\mu$ is Besicovitch rectifiable. If $E$ 
is a compact set
such that $\mathcal{H}^1(E) <\infty$, then $E$ is Besicovitch rectifiable if and only if the principal value of 
the Cauchy integral 
$C_{\mathcal{H}^1|E}(x)$
exists for $\mathcal{H}^1$ a.e. $x\in E$.

\bigskip

The last claim completely solves the conjecture of Mattila [Ma]. Mattila proved this result with the extra assumption of 
``non-uniform homogeneity":
$$
\liminf_{r\rightarrow 0}\frac{\mathcal{H}^1(B(x,r))}{r}>0\,\,\text{for}\,\,\mathcal{H}^1\,\,\text{a.e.}\,\, x\in E\,.
$$

{\bf Proof.}
We will prove the first assertion. The rest is not difficult to deduce. We choose $L,M$ so large that

$$
\mu(G_L(1)\cup H_M) <\frac12\,.
$$
  
Theorem 1a says that with $\Psi(x):=\dist(x,\C\setminus (G_L(1)\cup H_M))$ the operator $K_{\Psi}$ is bounded on $L^2(\mu)$ with 
norm at most $ALM$. Consider $G_{L,M}:= G_L(1)\cup H_M$, $f=1_{\C\setminus G_{L,M}}$. Then (*) implies 
$$
c^2(\mu|\C\setminus G_{L,M})) = \int_{\C\setminus G_{L,M})}|C1_{\C\setminus G_{L,M})}|^2\,d\,\mu= 
\int_{\C\setminus G_{L,M})}|K_{\Psi}f|^2\,d\,\mu \le ALM\,.
$$
The first equality is the famous formula of Melnikov-Verdera from [MV]. Notice that $\mu(\C\setminus G_{L,M})> \frac12$. 
But choosing larger $L,M$ we can make $\mu(\C\setminus G_{L,M})$ as close to $1$ as we wish (recall that our convention is that
$\|\mu\|=1$). So we can scoop the measure $\mu$ by pieces with finite curvature $c^2$. This proves the first claim of the Theorem. 

The  other
claims now follow easily. 
For example, the a.e. existence of the principal value of the Cauchy integral 
$C_{\mathcal{H}^1|E}(x)$ implies the a.e. finiteness of the maximal singular integral $C^*_{\mathcal{H}^1|E}(x)$. 
This and the non-uniform linear growth of 
$\mathcal{H}^1$ (it always has this property) imply that a.e. finiteness of the maximal singular integral $(K^*_{\Phi}1)(x)$ 
holds for $\Phi(x):=
\dist(x,\C\setminus H_M)$ and any $M$ (see Section VIII). Then $\mu=\mathcal{H}^1|E$ is a Besicovitch-Melnikov-Verdera rectifiable
measure  (by the first claim). The
result of David and L\'eger now shows that $E$ is Besicovitch rectifiable. 

\bigskip

To prove Theorems 1 (and, so, to prove the second claim of Theorem 1a) we need to use our decomposition into good 
and bad functions. Recall that we used the probability  space $(\Omega, P)$
of {\it pairs} of random dyadic lattices, $\omega=(\omega_1,\omega_2)$, here $\omega_i$ ``enumerates" the $i$-th ($i=1,2$) 
dyadic lattice $D_i$. 
These lattices $D_1,D_2$ are independent. We used also the notion of ``good" and ``bad" squares in 
$D_1$ and $D_2$. We also used the decomposition of sure functions $f,g\in L^2(\mu)$ to random functions

$$
  f=f_{good} +f_{bad},\,\,\,\, g=g_{good} + g_{bad}\,,
$$
$$
f_{bad} =\Sigma_{Q\in D_1, Q\,\,is\,\,bad} \Delta_Q f,\,\,\,\,g_{bad} =\Sigma_{R\in D_2, R\,\,is\,\,bad} \Delta_R g \,. 
$$
The proof of Theorem 1 (and, so, of 1a) is based on the following lemma.

\bigskip

{\bf Lemma 5:}
Let $\mu$ be a measure with non-uniform linear growth, let $H_M$ be the union of all $M$-non-Ahlfors discs,  let 
$$
\Phi(x) \ge \dist(x,\C\setminus H_M)\,,
$$
and let $\Phi$ is a $1$-Lipschitz function. Consider 
$(K_{\Phi,\e}f)(x) := \int_{|y-x|\ge\e}k_{\Phi}(x,y) f(y)\,d\,\mu(y)$.
Let $B$ be a finite constant such that 
$$
(K^*_{\Phi}1)(x) := \sup_{\e>0} |(K_{\Phi,\e}1)(x)| \le B \,\,\,\text{for} \,\,\mu\,\,\text{a. e.} \,\, x\,.
$$
Then
$$
|(K_{\Phi} f_{good}, g_{good})| \le ABM\|f\|\|g\|\,,
$$
$$
|(K_{\Phi} f, g)| \le ABM||f\|\|g\| + \|K_{\Phi}\| R(\omega,f,g),
$$
where the expectation of the remainder $R(\omega,f,g)$ has the following estimate: $\E R(\omega,f,g) \le \frac12 \|f\|\|g\|$.

\bigskip

The inequalities of the lemma imply immediately Theorem 1. In its turn, the last inequality follows from the 
first one and  the fact proved in the previous sections:

$$
 \E\|f_{bad}\|^2 \le 45^{-239}\|f\|^2\,,
$$
$$
\E\|g_{bad}\|^2 \le 45^{-239}\|g\|^2\,.
$$

The proof of the first inequality of the lemma takes a good part of previous sections. So, Theorem 1 and 1a 
are  proved.

\bigskip

What if we replace the function $1$ by a complex valued function $b$ (even, say, real valued but not always positive) in one of 
our main assumptions:
$$
(K^*_{\Phi}b)(x) := \sup_{\e>0} |(K_{\Phi,\e}b)(x)| \le B <\infty \,\,\,\text{for} \,\,\mu\,\,\text{a. e.} \,\,x \,?
$$
This is equivalent to still having the function $1$ but having {\it complex} measure $\mu$. We prefer to denote by $\mu$ only
positive measures, 
and to use the symbol $\nu$ for $b\,d\,\mu$. So now $b$ is an $L^{\infty}(\mu)$-function of norm $1$, and we assume that

$$
(C^*b)(x) := \sup_{\e>0} |(C^{\e}b)(x)| < \infty \,\,\,\text{for} \,\,\mu\,\,\text{a. e.} \,\,x\,. 
$$
We do not write the subscript $\mu$ because it will be always the same $\mu$.

We still assume everywhere below that $\mu$ has the non-uniform linear growth condition (unless it is stated otherwise).

Now we are in the framework of the $Tb$ theorem rather than the $T1$ theorem. The main problem we encounter is that our $b$ will 
{\it not be accretive}. 
The second problem (we always have it in this paper) is that $\mu$ has no doubling property.

We start again by considering the set $H_M$ of all $M$-non-Ahlfors discs for $\mu$. Again we can see that our assumption 
on $(C^*b)(x)$ 
implies (see Lemma 1)  {\it the a.e finiteness of the maximal singular operator}:

$$
(K^*_{\Phi} b)(x) <\infty \,\,\,\text{for} \,\,\mu\,\,\text{a. e.} \,\,x 
$$
for every $1$-Lipschitz  $\Phi$ such that $\Phi(x) \ge\dist(x,\C\setminus H_M))$. Exactly as before we can introduce 
the sets $G_L=G_L(b)=$ the union of $B(x,2\e_0(x))$, where $\e_0$ is the maximal radius for which 
$|(K_{\Phi,\e_0}b)(x)|\ge L$, 
and
$G_{L,M}= G_L\cup H_M$.

\bigskip

{\bf Lemma 6:}
Let 
$
(K^*_{\Phi} b)(x) <\infty \,\,\,\text{for} \,\,\mu\,\,\text{a. e.} \,\,x 
$
hold. Then $\mu(G_L)\rightarrow 0$ when $L\rightarrow \infty$.

\bigskip

Let $\Psi(x) =\dist(x,\C\setminus G_{L,M})$. The set $F_{\Psi} =\{x\in\supp\mu: \Psi(x) =0\}$ has  measure close to 
the full measure 
of $\mu$. Unfortunately, unlike 
in Theorem 1a, we cannot say that $K_{\Psi}$ is  bounded on $L^2(\mu)$. The place where the proof will break down is Lemma 5. 
The estimate
for good functions will not work. This is because $\Delta_Q f$ is now adapted to the function $b$. On squares where accretivity of 
$b$ becomes
very bad (or non-existent) the pieces 
$\Delta_Q f$ will blow up because the accretivity constant lives in their denominators. This was impossible for 
$b=1$---it is
accretive in any scale. To deal with this problem of non-accretivity of $b$ we need even more randomness: first let 
us assume that for a
certain positive $\eta$ the union of squares (``non-accretive squares") $Q\in D_1$ such that

$$
|\int_Q b\,d\,\mu|<\eta \,\mu(Q)
$$
has total measure less than $\delta$, and this is uniformly true for every random lattice $D_1$ ( so for $D_2$ also).

Let $T_i$ be the family of ``non-accretive" squares of $D_i$, $i=1,2$, in the above sense. Let $\omega \in \Omega$. 
Let $\mathcal{T}^{\omega}_i =\cup_{Q\in T_i} Q$.

We have the (strange) assumption that
$$
\qquad\qquad\qquad\qquad\qquad \mu(\mathcal{T}^{\omega}_1 \cup \mathcal{T}^{\omega}_2)\le \delta \,\, 
\text{for all}\,\,\omega\in \Omega\,.\qquad \qquad \qquad 
\qquad\qquad  (**)
$$ 

{\bf Lemma 7:}
Consider any $1$-Lipschitz function 
$\Phi_{\omega}$ such that $\Phi_{\omega}(x) \ge\dist(x, \C\setminus (G_{L,M}\cup \mathcal{T}^{\omega}_1 \cup \mathcal{T}^
{\omega}_2))$. Then
$$
|(K_{\Phi_{\omega}} f_{good}, g_{good})|\le ALM\,\eta^{-2} \|f\|\|g\|\,.
$$
 
\bigskip

\noindent This lemma is the result of our previous sections. Using the last inequality we can obviously write

$$
|(K_{\Phi_{\omega}} f, g)|\le ALM\,\eta^{-2} \|f\|\|g\| + \|K_{\Phi_{\omega}}\| R(\omega,f,g)
$$
with $R(\omega,f,g)$ having small average (exactly as in Lemma 5). But now it is not clear what to do with the random norm 
$\|K_{\Phi_{\omega}}\|$. We can consider a sure function $\Phi=\sup \Phi_{\omega}$. It is again $1$-Lipschitz and again

$$
|(K_{\Phi} f, g)|\le ALM\,\eta^{-2} \|f\|\|g\| + \|K_{\Phi}\| R(\omega,f,g)
$$
with small $\E R(\omega,f,g)$. So the bound for the norm of $\|K_{\Phi}\|$ will follow by averaging  
the previous inequality. 

But this is useless because our ``nice" set
$$
F_{\Phi} = \{x: \Phi(x) =0\} =\cap_{\omega} F_{\Phi_{\omega}} = \cap_{\omega} \{x: \Phi_{\omega}(x) =0\}
$$
and these random sets could easily have empty intersection.

We have, however, an extra ``strange" assumption  (**): 
$ \mu(\mathcal{T}^{\omega}_1 \cup \mathcal{T}^{\omega}_2)\le \delta \,\, \text{for all}\,\,\omega\in \Omega$. 
Then for sufficiently large $L,M$ we have

$$
\mu(G_{L,M} \cup  \mathcal{T}^{\omega}_1 \cup \mathcal{T}^{\omega}_2) \le 2\delta \,\, \text{for all}\,\,\omega\in \Omega\,.
$$

Notice that this means (by Fubini's theorem and Tchebyshov's inequality) that

$$
\mu\{x: P\{\omega : x \in G_{L,M} \cup  \mathcal{T}^{\omega}_1 \cup \mathcal{T}^{\omega}_2\} \le \sqrt{2\delta}\} 
\ge 1-\sqrt{2\delta}\,.
$$ 
We can now consider the expectation of $\Phi_{\omega}$ rather than maximum. Moreover, as we have done in Section IV, 
let us now consider 
{\it the truncated mathematical expectation}:

$$
\Psi(x) := \inf \{\E(\Phi_{\omega}(x) \, 1_S(\omega)): S\subset \Omega ,\, P(S) = 1-\sqrt{2\delta}\}\,.
$$

Now we have the good estimate for the zero set $F_{\Psi}$:

$$
\mu(F_{\Psi}) \ge 1-\sqrt{2\delta}\,.
$$

On the other hand, Lemma 7 can leads us to

\bigskip 

{\bf Theorem 2:}
Let $\mu$ have a non-uniform linear growth condition. Assume the a.e. finiteness of maximal singular operator, namely:
$$
(C^*b)(x) := \sup_{\e>0} |(C^{\e}b)(x)| < \infty \,\,\,\text{for} \,\,\mu\,\,\text{a. e.} \,\,x 
$$
We also assume  that $\mu$ has the non-uniform linear growth condition. Assume also (**). 
Let $\Phi_{\omega}(x) =\dist(x, \C\setminus (G_{L,M}\cup \mathcal{T}^{\omega}_1 \cup \mathcal{T}^{\omega}_2))$, and 
let $\Psi$ be the truncated mathematical
expectation of $\Phi_{\omega}$ defined above. Then

\noindent 1)
$$
|(K_{\Psi} f,g)| \le ALM\eta^{-2}\|f\|\|g\| + (\|K^*_{\Psi}\|+AM)\,R(\omega,f,g)
$$
where $\E R(\omega,f,g) \le A\delta \|f\|\|g\|$. 

\noindent 2) In particular, $\|K_{\Psi}\|_{L^2(\mu)\rightarrow L^2(\mu)} \le A LM\eta^{-2}$.

\noindent Automatically, for the set $F_{\Psi}$ (whose measure $\mu(F_{\Psi})\ge 1-\sqrt{2\delta}$) we have
$$
\|C\|_{L^2(F_{\Psi},\mu)\rightarrow L^2(F_{\Psi},\mu)} \le ALM\eta^{-2}\,.
$$

\bigskip

This theorem was proved by all the previous sections. However, the second part of the theorem requires the estimate of
$\|K^*_{\Psi}\|$ via
$\|K_{\Psi}\|$. This is done in [NTV2] for Ahlfors measures (i.e. measures having a uniform linear growth condition). Completely
similar reasoning for non-uniformly  Ahlfors measures (i.e. measures having a non-uniform linear growth condition) can be found in
Section XXV of the present paper. 

Theorem 2 gives the analytic part of Vitushkin's conjecture but  without the estimate of
how  large the rectifiable part of Vitushkin's compact is, and how rectifiable it is. This is because the assumption (**) 
does  not
happen very often. In fact, why should an arbitrary non-zero complex function $b$ (and in applications we usually do not know 
anything else  about $b$) 
be accretive except for a small set? In our previous sections we achieve (**) by localizing our considerations to 
an unspecified small disc
around a Lebesgue point $x_0$ of $b$, where $b(x_0) \neq 0$. Clearly, this way will not lead us to the proof of 
quantitative version of
Vitushkin's conjecture.

However, there is one piece of information which is usually available about $b$, and which has not been used so far. 
Namely, we have the 
accretivity of $b$ in one---the highest---scale:
$$
\|b\|_{\infty}=1,\,\,\, |\int_{\C} b\,d\,\mu |= \gamma > 0\,.\qquad\qquad\qquad\qquad\qquad\qquad\qquad 
(\gamma)
$$

This brings us to the quantitative version of $Tb$ theorem, where $b$ has accretivity only at the highest scale. 
We do not assume (**), 
but we assume $(\gamma)$. As always $\|\mu\|=1$.

\bigskip 

{\bf Theorem 3:}
Assume the a.e. finiteness of the maximal singular operator:
$$
(C^*b)(x) := \sup_{\e>0} |(C^{\e}b)(x)| <\infty \,\,\,\text{for} \,\,\mu\,\,\text{a. e.} \,\,x \,.
$$
Also assume  that $\mu$ has the non-uniform linear growth condition. Assume also $(\gamma)$. Then there exists a set 
$F$, $\mu(F) \ge \frac{3\gamma}{16}$, such that
$$
\|C\|_{L^2(F,d\mu)\rightarrow L^2(F,d\mu)} \le A L(\gamma) M(\gamma) \gamma^{-20}
$$
where $M(\gamma) =\inf\{M: \mu(H_M) <\frac{\gamma}{32}\}$ and $L(\gamma)=\inf\{L: \mu(G_L\setminus H_{M(\gamma)}) ) 
<\frac{\gamma}{32}\}$.

\bigskip 

The next theorem is the promised quantitative version of Vitushkin's conjecture. We will obtain it (easily) as a corollary of 
Theorem 3.

\bigskip

{\bf Theorem 4 (quantitative version of Vitushkin's conjecture):}
Let $E\subset \C$ be a compact set such that its length $\mathcal{H}^1(E) =M<\infty$ and its analytic capacity $\gamma(E)=
\gamma >0$. 
Then there exists a set $F$, $ \mathcal{H}^1(F) \ge \frac{\gamma}{16}$, such that 
$$
c^2(\mathcal{H}^1|F) \le A \, (\frac{\text{diam} E}{\gamma})(\frac{M}{\gamma})^{42} \mathcal{H}^1(F)\,.
$$

\bigskip

{\bf Proof of Theorem 3:} Consider $T_i=$ maximal squares from $D_i$ such that 
$$
|\int_Q b\,d\,\mu| \le \frac{\gamma}{2} \mu(Q)\,.
$$
Put $\mathcal{T}_i= \cup_{Q\in T_i}Q$, $i=1,2$. For brevity, let $E=\supp\mu$. Using $(\gamma)$ we have
$$
|\int_E b\,d\,\mu| = |\int_{\mathcal{T}_1} b\,d\,\mu| + |\int_{E\setminus \mathcal{T}_1} b\,d\,\mu|=
$$
$$
|\Sigma_{Q\in T_1}\int_Q b\,d\,\mu| + |\int_{E\setminus \mathcal{T}_1} b\,d\,\mu|  
\le \frac{\gamma}{2}\Sigma_{Q\in\mathcal{T}_1}\mu(Q) + 
\mu(E\setminus \mathcal{T}_1)   
\le
$$
$$
\frac{\gamma}{2} + \mu (E\setminus \mathcal{T}_1)\,.\qquad\qquad\qquad\qquad\qquad\qquad\qquad\qquad
$$
Therefore,
$$
\mu(E\setminus \mathcal{T}_1^{\omega}), \,\, \mu(E\setminus \mathcal{T}_2^{\omega})\ge \frac{\gamma}{2}\,.
$$
We wrote the superscript $\omega$  to emphasize that these are random sets. We want to show that for some detectable 
(=not very small) set 
of $x\in E$ the probability $p(x) := P\{\omega: x\in E\setminus (\mathcal{T}_1^{\omega}\cup\mathcal{T}_2^{\omega})\}$ 
is not too small. Denote $p_1(x) :=
P\{\omega: x\in E\setminus \mathcal{T}_1^{\omega}\}$. Notice that  the sets 
$E\setminus\mathcal{T}_1^{\omega},\,E\setminus \mathcal{T}_2^{\omega}$ are independent and 
that $E\setminus (\mathcal{T}_1^{\omega}\cup\mathcal{T}_2^{\omega})\} = 
(E\setminus\mathcal{T}_1^{\omega})\cap (E\setminus \mathcal{T}_2^{\omega})$. Therefore, $p(x)
=p_1(x)^2$. Also
$$
\int_E p_1(x)\,d\,\mu = \E\int 1_{E\setminus \mathcal{T}_1^{\omega}} \,d\,\mu =\E\mu(E\setminus \mathcal{T}_1^{\omega}) 
\ge \frac{\gamma}{2}\,.
$$

Now let us split $E= S\cup L$, where $S:=\{ x\in E: p_1(x) \le \frac{\gamma}{4}\}$ and  
$L:=\{ x\in E: p_1(x) > \frac{\gamma}{4}\}$. 
Then we have $\mu(L)\ge \frac{\gamma}{4}$. For $x\in L$, $p(x)=p_1^2(x) > \frac{\gamma^2}{16}$. 
For the sake of brevity we denote $\beta =
\frac{\gamma^2}{16}$. So
$$
\mu\{x\in E: P\{ \omega : x \in E \setminus (\mathcal{T}_1^{\omega}\cup\mathcal{T}^{\omega}_2)\}>\beta\} \ge \frac{\gamma}{4}\,.
$$

Now let us choose $M=M(\gamma),k=L(\gamma)$ to be smallest numbers such  that 
$$
\mu(H_M)\le \frac{\gamma}{32} ,\,\,\, \mu(G_L\setminus H_M) \le \frac{\gamma}{32}\,.
$$

Consider $O^{\omega}:= G_{L,M} \cup \mathcal{T}_1^{\omega}\cup \mathcal{T}_2^{\omega}$.  Put 
$\Phi_{\omega}(x) := \dist(x, \C\setminus O^{\omega})$. Thus,
$$
\mu\{x\in E: P\{\omega :\Phi_{\omega} (x)=0\}>\beta\} >\frac{3\gamma}{16}\,.
$$

Let us introduce sure $1$-Lipschitz function
$$
\Phi_0 (x) :=\inf_{S\subset\Omega,\,P(S) =\beta} \sup_{\omega\in S}\Phi_{\omega}(x)\,.
$$
Let us also fix a small positive number $\tau$ and put
$$
\Phi(x) := \Phi_0(x) +\tau\,.
$$
All estimates in the future will not depend on $\tau$.

We know that the zero set $F:=F_{\Phi_0}$ has detectable measure:
$$
\mu(F) >\frac{3\gamma}{16}\,.
$$

\bigskip

We will need a small modification of Lemma 7 of this section.

{\bf Lemma 7a:}
Consider any $1$-Lipschitz function 
$\Phi_{\omega}$ such that $\Phi_{\omega}(x) \ge\dist(x, \C\setminus (G_{L,M}\cup \mathcal{T}^{\omega}_1 \cup \mathcal{T}^
{\omega}_2))$. Fix a small positive number $\e$. Then we can decompose $f=f_{good}+f_{bad}, g=g_{good} +g_{bad}$ in such 
a way that

$$
\E\|f_{bad}\| \leq \e \|f\|,\,\, \E\|g_{bad}\|\leq \e\|g\|\,,
$$
and
$$
|(K_{\Phi_{\omega}} f_{good}, g_{good})|\le ALMC(\e)\,\eta^{-2} \|f\|\|g\|,\,\, \text{where}\,\, C(\e) \leq A\e^{-8}\,.
$$

\noindent All the previous sections were devoted to the proof of such a statement (called Lemma 7 in this section) with a fixed
small absolute constant
$\e$ (it has been chosen to be $45^{-239})$. But the same proof gives Lemma 7a because in our calculations in Section XXII we can
choose a very large $m$ and a very large $\wt{M}$ in accordance with the smallness of $\e$. They can be chosen to achieve our
first inequality of Lemma 7a. Then the second inequality of Lemma 7a follows from the bookkeeping of the estimate of the bilinear
form of the operator
$K_{\Phi_{\omega}}$ on good functions.
 
\bigskip

{\bf Main Lemma:}
Operator $C_{\Phi}$ is bounded on $L^2(\mu)$ by $AL(\gamma)M(\gamma)\gamma^{-20}$ 
(and the bound does not depend on $\tau$).

\bigskip

{\bf Proof.}
Fix $\e = a \gamma^2$.  Here $a$ is a small positive absolute constant. Recall that the splitting into good and 
bad functions can be made dependent on a number $\e$. Then 
$$
\E\|f_{bad}\|\le \e\|f\|,\,\,\, \E\|f_{bad}\|\le \e\|f\|\,.
$$
Lemma 7a (with $\eta =\gamma /2$) then  states the following:

$$
|(K_{\Phi\vee\Phi_{\omega}}f_{good}, g_{good})| \le ALMC(\e)\gamma^{-2},\,\,\,\text{with} \,\, C(\e) \le A \,\e^{-8}\,.
$$
We used the notations $\Phi\vee\Phi_{\omega}=\max(\Phi,\Phi_{\om})$. We use now Lemma 1.
$$
|(K_{\Phi\vee\Phi_{\omega}}f,g)|\le 
|(K_{\Phi\vee\Phi_{\omega}}f_{good},g_{good})| + |(C_{\Phi\vee\Phi_{\omega}} f_{bad}, g_{good})| + 
$$
$$
 |(C_{\Phi\vee\Phi_{\omega}} f_{good}, g_{bad})| +|(C_{\Phi\vee\Phi_{\omega}} f_{bad}, g_{bad}) + A\|M_{1,\Phi}f||\|g\|\,.
$$
Notice that $\Phi(x)\ge \dist (x, \C\setminus G_{L,M})$. Using Lemma 3 we make an estimate in the last term: 
$$
\|M_{1,\Phi}f\|\le AM\|f\|\,.
$$

The estimate of $ |(C_{\Phi\vee\Phi_{\omega}} f_{bad}, g_{good})| + ... $ involves an important lemma and several notations.
Let $k_{\om}(x,y)$ denote the kernel of $K_{\Phi\vee\Phi_{\omega}}$. Let $c_{\om}(x,y)$ denote the kernel of 
$C_{\Phi\vee\Phi_{\omega}}$.

Notice that 
$$
p_{\om}(x,y) :=|k_{\om}(x,y)-c_{\om}(x,y)|
$$ 
is a ``Poisson" type kernel. In particular,
$$
\int p_{\om}(x,y)|f(y)|\,d\,\mu(y) \le A \,(M_{1,\Phi}f)(x)
$$
Consider the averaging of the kernels: $k=\E k_{\om}$, $c=\E c_{\om}$, $p=\E p_{\om}$.
The same ``Poisson" property holds then for the average $p=\E p_{\om}$):

$$
\int p(x,y)|f(y)|\,d\,\mu(y) \le A \,(M_{1,\Phi}f)(x)\,.
$$
Let us also introduce operators $c^*, k^*$ as follows:
$$
(c^*f)(x):= \sup_{r>0}|\int_{|y-x|\ge r} c(x,y) f(y)\,d\,\mu(y)|,\,\,(k^*f)(x):= \sup_{r>0}|\int_{|y-x|\ge r} k(x,y) 
f(y)\,d\,\mu(y)|\,.
$$
The same ``Poisson" property holds then for the comparison of $k^*$ and $c^*$ 
(notice that $k$, $c$ are defined in such a way that $|k(x,y)|, |c(x,y)| \leq \frac1{\Phi(x)}$:
$$
(c^*f)(x)\le (k^*f)(x) + (M_{1,\Phi}f)(x)\,.
$$

We are ready to formulate the main inequalities:

$$
\qquad\qquad|(C_{\Phi}f)(x)|\le \frac{A}{\gamma^2} ((c^*f)(x)+ (M_{1,\Phi}f)(x))\,,\qquad\qquad\qquad\qquad\qquad\qquad (MI)
$$
$$
\qquad\qquad|(C_{\Phi\vee\Phi_{\om}}f)(x)|\le \frac{A}{\gamma^2} ((c^*f)(x)+(M_{1,\Phi}f)(x))\,.
\qquad\qquad\qquad\qquad\qquad (MI)
$$

Let us use (MI) to estimate 
$$ 
|(C_{\Phi\vee\Phi_{\omega}} f_{bad}, g_{good})| + |(C_{\Phi\vee\Phi_{\omega}} f_{good}, g_{bad})| +|(C_{\Phi\vee\Phi_{\omega}} 
f_{bad},
g_{bad})|\,. 
$$ 
After that we will prove (MI).
By (MI), Lemma 1 and the Poisson property  for the comparison of $k^*$ and $c^*$:
$|(C_{\Phi\vee\Phi_{\omega}} f_{bad}, g_{good})|\le \frac{A}{\gamma^2}(\|(c^*f_{bad})\|\|g\|+
\|M_{1,\Phi}f)\|\|g\|)\le \frac{A}{\gamma^2}\|(k^*f_{bad})\|\|g\| + \frac{A}{\gamma^2}\|M_{1,\Phi}f\|\|g\|$. We continue:
$$
|(C_{\Phi\vee\Phi_{\omega}} f_{bad}, g_{good})|\le \frac{A}{\gamma^2}\e\|k^*\|\|f\|\|g\| + \frac{A}{\gamma^2}M\|f\|\|g\|\,.
$$

Collecting our estimates for the good and bad function together, we get

$$
|(K_{\Phi\vee\Phi_{\omega}}f,g)|\le 
ALM\gamma^{-2} \e^{-8}\|f\|\|g\| + \frac{A}{\gamma^2}\e\|k^*\|\|f\|\|g\| + \frac{A}{\gamma^2}M\|f\|\|g\|\,.
$$

We already fixed $\e = a \gamma^2$. Thus  (with very small absolute $a$)

$$
|(K_{\Phi\vee\Phi_{\omega}}f,g)|\le 
ALM\gamma^{-18}\|f\|\|g\| + Aa\|k^*\|\|f\|\|g\| + AM\gamma^{-2}\|f\|\|g\|\,.
$$

\noindent Recall that $k$ denotes the average of the kernel of $K_{\Phi\vee\Phi_{\omega}}$. Averaging the previous inequality we
get

$$
\|k f\| \le ALM\gamma^{-18}\|f\| + Aa\|k^*\|\|f\| + \frac{A}{\gamma^2}M\|f\|\,.\qquad\qquad\qquad\qquad\qquad\qquad
\qquad\qquad\,(kI)
$$

In  Theorem 7.1 of [NTV2] it is proved that $ \|k^* f\| \le A_1 C\|f\| + A_2 C\|k\|\|f\|$, where $C$
stands for the Calder\'on-Zygmund constant of the kernel. Theorem 7.1 of [NTV2] is valid for operators  with Calder\'on-Zygmund
kernels.  This is
the case here because the averaging
$k$ of the Calder\'on-Zygmund kernels
$k_{\om}$  is still a
Calder\'on-Zygmund kernel. 

\noindent However, there is a difference between the sitation in [NTV2]
and the situation here. In [NTV2] one assumed that the measure $\mu$ has a uniform linear growth condition. 
Our $\mu$, however, has  only  the non-uniform linear growth condition (we call such measures non-uniformly Ahlfors measures). 
We are
going to formulate now an analog of Theorem 7.1 from [NTV2] that is valid for non-uniformly Ahlfors measures.
First, recall that given a Calder\'on-Zygmund kernel and a measure $\mu$ we say that the 
operator $T$ with kernel $k$ (see [NTV3])
is a Calder\'on-Zygmund operator if it is bounded on $L^2(\mu)$. Also recall that 
$$
\tilde{M}_{\beta}g(x) :=\sup_{r>0}\frac1{\mu(B(x,3r)}\left(\int_{B(x,r)}|g(y)|^{\beta}\,d\mu(y)\right)^{\frac1{\beta}}\,.
$$
When $\beta =1$ we write $\tilde{M}g(x)$ instead of $\tilde{M}_{1}g(x)$.

\bigskip

{\bf Theorem 5}. 
Let $\mu$ be a  non-uniformly Ahlfors measure. Fix a positive number $M$, and let 
$\mathcal{R}(x):=\sup\{r>0\,:\, \mu(B(x,r))>Mr\}$. Let $k(x,y)$ be a Calder\'on-Zygmund kernel having Calderr\'on-Zygmund
constant $C$ and such that
$$
|k(x,y)| \le \min \Bigl[\frac1{\mathcal{R}(x)},\frac1{\mathcal{R}(y)}\Bigr]\,.
$$
Let $T$ be a Calder\'on-Zygmund operator with kernel $k$. Fix $\beta\in (1,2)$. Then the
following Cotlar type inequality holds:
$$
(T^*f)(x) \le A_1 C [\tilde{M}Tf](x) + A_2 C M \tilde{M}_{\beta}f(x) + A_3 C \|T\|_{L^2(\mu)\rightarrow L^2(\mu)}
\tilde{M}_{\beta}f(x)\,.\qquad\qquad\qquad \,(CI) 
$$
 
\bigskip

The proof follows exactly the lines of the proof of Theorem 7.1 of [NTV2]. But for the sake of completness we give a full 
proof in Section XXV.

\bigskip

Combining this result with inequality (kI), we get

$$
\|k^* f\| \le ALM\gamma^{-18}\|f\| + A\,a\|k^*\|\|f\| + AM\|f\|\,.
$$
Finally, using the fact that $a$ is very small we get the estimate for the maximal singular operator:
$$
\|k^* f\| \le 2ALM\gamma^{-18}\|f\|\,.
$$
Now let us use again the ``Poisson" property  for the comparison of $k^*$ and $c^*$:
$
(c^*f)(x)\le (k^*f)(x) + (M_{1,\Phi}f)(x)
$ to get

$$
\|c^* f\| \le ALM\gamma^{-18}\|f\|\,.
$$

\noindent Let us use the first part of the main inequality (MI) to conclude now that

$$
\|C_{\Phi}f\| \le ALM\gamma^{-20}\|f\|\,.
$$

\noindent The main Lemma is proved. 

\bigskip

\noindent We are left to prove (MI).

\noindent The proof of (MI) is based on two ingredients---the calculation of the kernel $c(x,y)$ (average of 
$c_{\om}(x,y)$) and on the following lemma.

\noindent As usual, given $R\ge 0$, we denote by $(M_{1,R}f)(x)=\sup_{r>R}\frac1r\int_{B(x,r)}
|f(y)|\,d\,\mu(y)$.

\bigskip
 
{\bf Blanket Lemma:}
Let $b(x,y)$ be kernel such that $|b(x,y)|\le \frac{1}{|x-y|}$. Then we have a well-defined 
$(b^*f)(x):=\sup_{r>0}|\int_{|y-x|>r}b(x,y) f(y)\,d\,\mu(y)|$. 
Let $R > 0$ and let $\phi$ be a decreasing function on $[0,\infty)$, $0\le \phi\le 1$. 
Consider 
$$
(b^{\phi}_R f)(x):=|\int_{|y-x|>R}b(x,y)\phi(|x-y|) f(y)d\mu(y)|\,.
$$
Then
$$
(b^{\phi}_R f)(x) \le 2\,(b^*f)(x) + 2\,(M_{1,R}f)(x)\,.
$$

\bigskip

{\bf Proof.}
Consider annuli $A_k(x)=\{y: 2^{k-1}R\le |y-x|\le 2^{k}R\}$. 
Then
$$
(b^{\phi}_R f)(x) \approx \Sigma_{k\ge 1} \int_{A_k} b(x,y) \phi_k f(y)d\mu(y)
$$ 
where $\phi_k$ are some values (say, left end point values) of $\phi(t)$  for $t\in [2^{k-1}R, 2^{k}R]$, $k=1,2,...$.
More precisely ($\phi_0:=0$)

$$
(b^{\phi}_R f)(x) = \Sigma_{k\ge 1} (\phi_k-\phi_{k-1})\int_{|y-x|\ge 2^{k-1}R } b(x,y) f(y)\,d\,\mu(y) + \text{Discrepancy}\,.
$$

Thus, the monotonicity of $\phi$ implies

$$
|\text{The first term}| \le \phi_1 |\int_{|y-x|\ge R} b(x,y) f(y)\,d\,\mu(y)|+
$$
$$
\Sigma_{k\ge 2} (\phi_{k-1}-\phi_k)|\int_{|y-x|\ge 2^{k-1}R} b(x,y) f(y)\,d\,\mu(y)| \le 2\,(b^*f)(x) \sup \phi\,. 
$$
On the other hand, let us use the symbol $J_k$ to denote the jump (the oscillation) of the monotone function $\phi$ 
on the interval $[a_k, a_{k+1}]$. Then 
$$
|\text{Discrepancy}|\le \Sigma_{k\ge 1} J_k \frac1{2^{k-1}R}\int_{B(x, 2^{k}R)}|f(y)|\,d\,\mu(y)\,.
$$
We  continue the previous estimate 
as follows:
$$
|\text{Discrepancy}|\le 2\, (M_{1,R}f)(x)\Sigma_{k\ge 1} J_k\,.
$$
But $\phi$ was assumed to be monotone and $0\le\phi\le 1$, so the sum of the jumps is bounded by $1$.
The lemma is proved.

\bigskip

We continue the proof of (MI). Let $t\geq \Phi(x)$. Then 
$$
v(t):= P\{\om: \Phi\vee\Phi_{\om}(x) \le t\}\geq \gamma^2/16 \,.
$$ 
It is obvious that for $|x-y|< \Phi(x)$ we have
$v(|x-y|) =0$. 
Now let us compute the kernel $c(x,y)=\E c_{\om}(x,y)$. Clearly, 
$$
c(x,y) = \frac{v(|x-y|)}{x-y}=\frac{\chi_{\C\setminus B(x,\Phi(x))} v(|x-y|)}{x-y}\,.
$$

Put $\beta:=\gamma^2/16$.
To obtain (MI) we can apply the Blanket Lemma with $R(x) = \Phi(x)$ or $R(x) = \Phi\vee\Phi_{\om}(x) $, 
with $b(x,y)= \frac{c(x,y)}{\beta}$ 
and $\phi(t) =\frac{\beta}{v(t)}$. Theorem 3 is completely proved.

\bigskip

{\bf XXIV. The proof of Theorem 4. The quantitative version of Vitushkin's conjecture.}

\bigskip

Now let $\Gamma$ be a compact on $\C$ whose $\mathcal{H}^1$ measure is $L$ and whose analytic capacity is $\gamma$. 
We can think that $\Gamma$
consists of finitely many circle arcs. Consider $x\in \Gamma$ and $R(x)>0$ such that

$$
\frac{\mathcal{H}^1(B(x, R)\cap \Gamma)}{R} > \frac{160\pi L}{\gamma}\,.
$$

The union of such $B(x, R(x))$ is covered by $\cup B(x_j, 5R_j)$ and 
$$
\Sigma\, \mathcal{H}^1 (\partial B(x_j, 5R_j)) \le \frac{\gamma}{16}\,.
$$
Let $G$ be the boundary of the complement of $\cup_j B(x_j, 5R_j)\cup \Gamma$. Let $F=\Gamma\cap G$.
It is now clear that
$$
\mathcal{H}^1(G\setminus F)\le\frac{\gamma}{16}\,.
$$
It is easy to check that there is no $1000\,L/\gamma$-non-Ahlfors disc  for $G$. On the other hand, there exists
a function $b$ on $G$ such that its Cauchy integral is bounded by $1$ outside of $G$ (its Cauchy integral vanishes inside all 
$B(x_i,5R_i)$), such that
$\|b\|_{\infty}\le 1$, and such that
$$
|\int_G b\,d\,\mathcal{H}^1 |=\gamma\,.
$$
As $b$ we can take just the Ahlfors function of $\Gamma$ outside of $\cup B(x_i,5R_i)$ and zero inside.
In particular,

$$
(C^*\,b\,d\,\mathcal{H}^1)(x) \le A\frac{L}{\gamma}\,\,\,\text{for} \,\,\mathcal{H}^1\,\,\text{a.e} \,\, x\in G\,.
$$

Let us consider the normalized measure $\mu:=\mathcal{H}^1/L$ restricted on $\Gamma$. Then we are under the assumptions of Theorem
3, where  we can put 
$L:=\frac{1}{\gamma},\, M:=\frac{1}{\gamma},\, \gamma:=\frac{\gamma}{L}$ and get a set $F_0\subset E$ with $\mu(F_0)\ge
\frac{\gamma}{8L}$, that is with $\mathcal{H}^1(F_0)\ge
\frac{\gamma}{8}$, such that $\|C\|_{L^2(F_0,\mu)\rightarrow L^2(F_0,\mu)}\le A\gamma^{-2}(\gamma/L)^{-20}$.  That is 
$\|C\|_{L^2(F_0,\mathcal{H}^1)\rightarrow L^2(F_0,\mathcal{H}^1)}\le A\gamma^{-1}(\gamma/L)^{-21}$.

Consider $F^*:= F_0\cap F$. Then $\mathcal{H}^1(F_0)\ge\frac{\gamma}{8}$ and $\mathcal{H}^1(G\setminus F)\le\frac{\gamma}{16}$ 
imply that 
$$
\mathcal{H}^1(F^*)\ge\frac{\gamma}{16}\,.
$$ 
The advantage of $F^*$ is that it is contained in the original set $\Gamma$ because $F$ is. Also we have
$$
\|C\|_{L^2(F^*,\mathcal{H}^1)\rightarrow L^2(F^*,\mathcal{H}^1)}\le A\gamma^{-1}(\gamma/L)^{-21}
$$
just because $F^*\subset F_0$. The last relationship and the formula of Melnikov-Verdera shows

$$
c^2(\mathcal{H}^1|F^*)\le (A\gamma^{-1}(\gamma/L)^{-21})^2 \mathcal{H}^1(F^*)= A\gamma^{-2}(\gamma/L)^{-42}\mathcal{H}^1(F^*)\,.
$$
We tacitly assumed $\text{diam} \Gamma =1$. Thus, we have in general
$$
c^2(\mathcal{H}^1|F^*)\le  A(\frac{\text{diam} \Gamma}{\gamma})^2(\gamma/L)^{-42}\mathcal{H}^1(F^*)\,.
$$

Theorem 4 is proved.  

\bigskip

{\bf XXV. The proof of Theorem 5. Cotlar's inequality for non-uniformly Ahlfors measures.}

\bigskip

We start the proof by fixing $r>0, x\in \supp\mu$, and putting $\hat{r} =\max [r, \mathcal{R}(x)]$. Consider 
$(T^rf)(x):=\int_{y: |y-x|\ge r}k(x,y)f(y)\,d\mu(y)$. Put $r_j:= 3^j \hat{r}, \mu_j := \mu(B(x,r_j))$. 
Let $k$ be the smallest
index such that $\mu_{k+1} \leq 36 \mu_{k-1}$. It exists,  because otherwise, for every $k$, 
$\mu(B(x, \hat{r})) \leq 36^{-k}\mu_{2k}\leq 2M 36^{-k} r_{2k}$. This is because our radii are greater than 
$\mathcal{R}(x):=\sup\{r>0\,:\, \mu(B(x,r))>Mr\}$. We continue with $\mu(B(x, \hat{r})) \leq 2M 36^{-k} 3^{2k} \hat{r}=2M
2^{-2k}\hat{r}$. This contradicts the assumption $x\in\supp\mu$. 

Let $R:= r_{k-1}$. We estimate $|(T^r f)(x) - (T^{3R}f)(x)| \leq \int_{B(x,\hat{r})\setminus B(x,r)} |k(x,y)|\,|f(y)|\,d\mu(y)
+ \sum_{j=1}^k\int_{B(x,r_j)\setminus B(x,r_{j-1})}...$. The first term  vanishes if $\hat{r} >\mathcal{R}(x)$. Otherwise it is
bounded by
$$
\frac1{\mathcal{R}(x)}\int_{B(x,\hat{r})}|f(y)|\,d\mu(y)=\frac1{\hat{r}}\int_{B(x,\hat{r})}|f(y)|\,d\mu(y)\leq
$$
$$
\frac{\mu(B(x,3\hat{r})}{\hat{r}\mu(B(x,3\hat{r})}\int_{B(x,\hat{r})}|f(y)|\,d\mu(y)\,.
$$ 
And this is less than $AM\,\widetilde{M}f(x)$. Similarly 
$$
\int_{B(x,r_j)\setminus
B(x,r_{j-1})}|k(x,y)|\,|f(y)|\,d\mu(y)\leq\frac{\mu_{j+1}}{r_{j-1}\mu(B(x,r_{j+1})}\int_{B(x,r_j)}|f(y)|\,d\mu(y)\,.
$$
But we know that $r_{j-1} =3^{-k+j-1}r_k, \mu_{j+1}\leq 36(36)^{\frac{-k+j}{2}}\mu_k$. Hence $\frac{\mu_{j+1}}{r_{j-1}}\leq
36\cdot 3^{k-j+1}6^{-k+j}\frac{\mu_k}{r_k}\leq AM 2^{-k+j}$. Therefore,
$$
\sum_{j=1}^k\int_{B(x,r_j)\setminus
B(x,r_{j-1})}|k(x,y)|\,|f(y)|\,d\mu(y) \leq AM\sum_{j=1}^k 2^{-k+j}\frac1{\mu(B(x,r_{j+1})}\int_{B(x,r_j)}|f(y)|\,d\mu(y)\,.
$$
The last sum is obviously bounded by $AM\,\widetilde{M}f(x)$.
We finally get

$$
|(T^rf)(x) - (T^{3R}f)(x)| \leq AM\,\widetilde{M}f(x)\,.
$$

Now we need to estimate $(T^{3R}f)(x)$. Consider the average $V_R(x) :=\frac1{\mu(B(x,R)}\int_{B(x,R)} Tf\,d\mu$
First,
$$
|V_R(x)|\leq \frac{\mu(B(x,3R))}{\mu(B(x,R))}\widetilde{M}[Tf](x)\leq 36 \widetilde{M}[Tf](x)\,.
$$
Second,
$$
V_R(x) - (T^{3R}f)(x) = \int_{\C\setminus B(x,3R)}T^{'}[\delta_x - \frac1{\mu(B(x,R)}\chi_{B(x,R)}\,d\mu]f\,d\mu
- 
$$
$$
\frac1{\mu(B(x,R)}\int_{B(x,R)}T[f\chi_{B(x,3R)}]\,d\mu = I + II\,.
$$
Here $T^{'}$ denotes the operator with kernel $k(y,x)$.

\noindent {\bf Estimate of I.} Put $\eta = \delta_x - \frac1{\mu(B(x,R)}\chi_{B(x,R)}\,d\mu$. All radii greater
than $3R$ are $M$-Ahlfors for $\mu$. This and the fact  that $\eta(\C)=0$ allows us to use the Calder\'on-Zygmund property
of  $k(y,x)$ to prove as usual (see [NTV2] for example) that
$I \leq AM\,\|\eta\|\,\widetilde{M}f(x) \leq AM\,\widetilde{M}f(x)$.

\noindent {\bf Estimate of II.} Fix $\beta \in (1,2)$. Let $1/\al +\/\beta =1$.
$$
|II| \leq \frac1{\mu(B(x,R)}\|\chi_{B(x,R)}\|_{L^{\al}(\mu)}\|T(f\chi_{B(x,3R)})\|_{L^{\beta}(\mu)}\leq
\|T\|_{\beta}\frac{(\int_{B(x,3R)}|f|^{\beta}\,d\mu)^{\frac1{\beta}}}{\mu(B(x,R))^{\frac1{\beta}}}\,.
$$
Here we abbreviate $\|T\|_{\beta}:= \|T\|_{L^{\beta}(\mu)\rightarrow L^{\beta}(\mu)}$. We can continue
$$
|II| \leq \|T\|_{\beta}\frac{\mu(B(x,9R))^{\frac1{\beta}}\,(\widetilde{M}_{\beta}f)(x)}{\mu(B(x,R))^{\frac1{\beta}}}\leq
$$
$$
36^{\frac1{\beta}} \|T\|_{\beta}(\widetilde{M}_{\beta}f)(x)\leq A\|T\|_{\beta}(\widetilde{M}_{\beta}f)(x)\,.
$$

To estimate  $\|T\|_{\beta}$ via $\|T\|_{2}$ we need first

\bigskip

\noindent {\bf Estimate of weak type via $\|T\|_{2}$}.

\bigskip

\noindent {\bf Lemma (G. David).} For any measurable set $F$ and any point $x\in \supp\mu$,
$$
T^*\chi_F(x) \leq A_1 \widetilde{M}[T\chi_F](x) + A_2 M + A_3\|T\|_2\,.
$$
{\bf Proof.} Fix $x\in\supp\mu, r>0$. Put $\hat{r} =\max [r, \mathcal{R}(x)]$, where 
$\mathcal{R}(x):=\sup\{r>0\,:\, \mu(B(x,r))>Mr\}$. Consider $r_j =3^j\hat{r}$. Put $\mu_j := \mu(B(x,r_j))$. 
Let $k$ be the smallest
index such that $\mu_{k} \leq 6 \mu_{k-1}$. It exists. Otherwise, for every $k$, 
$\mu(B(x, \hat{r})) \leq 6^{-k}\mu_{k}\leq 2M 6^{-k} r_{k}$. This is because our radii are greater than 
$\mathcal{R}(x):=\sup\{r>0\,:\, \mu(B(x,r))>Mr\}$. We continue with $\mu(B(x, \hat{r})) \leq 2M 6^{-k} 3^{k} \hat{r}=2M
2^{-k}\hat{r}$. This contradicts the assumption $x\in\supp\mu$. Put $R=r_{k-1}$. We estimate $|(T^rf)(x) - (T^{3R}f)(x)|
 \leq \int_{B(x,\hat{r})\setminus B(x,r)} |k(x,y)|\,|\chi_F(y)|\,d\mu(y)
+ \sum_{j=1}^k\int_{B(x,r_j)\setminus B(x,r_{j-1})}...$. The first term  vanishes if $\hat{r} >\mathcal{R}(x)$. Otherwise it is
bounded by
$$
\frac1{\mathcal{R}(x)}\int_{B(x,\hat{r})}|\chi_F(y)|\,d\mu(y)=\frac1{\hat{r}}\int_{B(x,\hat{r})}|\chi_F(y)|\,d\mu(y)\leq 2M\,.
$$
Similarly 
$$
\int_{B(x,r_j)\setminus
B(x,r_{j-1})}|k(x,y)|\,|\chi_F(y)|\,d\mu(y)\leq\frac{\mu_{j}}{r_{j-1}}\,.
$$
But we know that $r_{j-1} =3^{-k+j-1}r_k, \mu_{j}\leq 6(6)^{-k+j}\mu_k$. Hence $\frac{\mu_{j}}{r_{j-1}}\leq
6 \cdot 3^{k-j+1}6^{-k+j}\frac{\mu_k}{r_k}\leq AM 2^{-k+j}$. Therefore,
$$
\sum_{j=1}^k\int_{B(x,r_j)\setminus
B(x,r_{j-1})}|k(x,y)|\,|\chi_F(y)|\,d\mu(y) \leq AM\sum_{j=1}^k 2^{-k+j}\leq AM\,.
$$

We finally get

$$
|(T^rf)(x) - (T^{3R}f)(x)| \leq AM\,.
$$

Now we need to estimate $(T^{3R}\chi_F)(x)$. Consider the average $V_R(x) :=\frac1{\mu(B(x,R)}\int_{B(x,R)} T\chi_F\,d\mu$
Firstly, by the choice of $R$, we have
$$
|V_R(x)|\leq \frac{\mu(B(x,3R))}{\mu(B(x,R))}\widetilde{M}[T\chi_F](x)\leq 6 \widetilde{M}[T\chi_F](x)\,.
$$
Second,
$$
V_R(x) - (T^{3R}f)(x) = \int_{\C\setminus B(x,3R)}T^{'}[\delta_x - \frac1{\mu(B(x,R)}\chi_{B(x,R)}\,d\mu]\chi_F\,d\mu
- 
$$
$$
\frac1{\mu(B(x,R)}\int_{B(x,R)}T[\chi_{F\cap B(x,3R)}]\,d\mu = I + II\,.
$$
Here $T^{'}$ denotes the operator with kernel $k(y,x)$.

\noindent {\bf Estimate of I.} Put $\eta = \delta_x - \frac1{\mu(B(x,R)}\chi_{B(x,R)}\,d\mu$. All radii greater
than $3R$ are $M$-Ahlfors for $\mu$. This and the fact  that $\eta(\C)=0$ allows us to use the Calder\'on-Zygmund property
of  $k(y,x)$ to prove as usual (see [NTV2] for example) that
$I \leq A\,\|\eta\|\,\sup_{\rho\geq R}\frac{\mu(B(x,\rho))}{\rho} \leq AM$.

\noindent {\bf Estimate of II.}
$$
|II| \leq \frac1{\mu(B(x,R)}\|\chi_{B(x,R)}\|_{L^{2}(\mu)}\|T(\chi_{F\cap B(x,3R)})\|_{L^{2}(\mu)}\leq
\|T\|_{2}\frac{(\int_{B(x,3R)}|\chi_F|^{2}\,d\mu)^{\frac1{2}}}{\mu(B(x,R))^{\frac1{2}}}\,.
$$
We can continue
$$
|II| \leq \|T\|_{\beta}\frac{\mu(B(x,3R))^{\frac1{2}}}{\mu(B(x,R))^{\frac1{2}}}\leq
$$
$$
6^{\frac1{2}} \|T\|_{2}\leq A\|T\|_{2}\,.
$$
The lemma is completely proved.

\bigskip

Now we are ready to repeat the considerations of Theorem 5.1 of [NTV2] (with small modifications due to the fact that $\mu$ is 
a non-uniformly Ahlfors measure).

We are going to prove now that
$$
\|T\|_{L^1(\mu)\rightarrow L^{1,\infty}} \leq A_1 C M + A_2 C \|T\|_2\,,\qquad\qquad\qquad\qquad\qquad\qquad\qquad
\qquad\qquad(W)
$$
where $C$ depend only on Calder\'on-Zygmund constants of the kernel of $T$.

Let $\nu\in M(\C)$ be a finite
linear combination of unit point masses with positive coefficients, 
i.e.,
$$
\nu=\sum_{i=1}^N \al_i\delta_{x_i}.
$$

\bigskip
{\bf Theorem 6.}
$$
\|T\nu\|\ci{L^{1,\infty}(\mu)}\le (A_1 C M + A_2 C \|T\|_2)\|\nu\|\,.
$$

\bigskip

Here there is no problem with the definition of $T\nu$: it is just the finite
sum \linebreak $\sum_{i=1}^N \al_i K(x,x_i)$, which makes sense everywhere except at
finitely many points.

\bigskip

{\bf Proof.}
In this proof $B(x,\rho)$ denotes closed ball, $B^{'}(x,\rho)$ denotes open ball.
Without loss of generality, we may assume
that $\|\nu\|=\sum_i\al_i=1$ (this is just a matter of normalization).
Thus we have to prove that

$\|T\nu\|\ci{L^{1,\infty}(\mu)}\le A_4$.
Fix some $t>0$ and
suppose first that $\mu(\C)>\frac1t$. Let
$B(x_1,\rho_1)$ be the smallest (closed) ball such that
$\mu(B(x_1,\rho_1))\ge\dfrac{\al_1}{t}$ (since the function
$\rho\to \mu(B(x_1,\rho))$ is increasing and continuous from the right,
tends to $0$ as $\rho\to 0$, and is greater than $\dfrac1t\ge
\dfrac{\al_1}{t}$ for
sufficiently large $\rho>0$,
such $\rho_1$ exists and is strictly positive).

Note that for the corresponding {\it open} ball $B'(x_1,\rho_1):=
\{y\in\C\,:\,\dist(x_1,y)<\rho_1\}$, we
have $\mu(B'(x_1,\rho_1))=\lim_{\rho\to \rho_1-0}\mu(B(x_1,\rho))\le
\dfrac{\al_1}{t}$.
Since the measure $\mu$ is $\sigma$-finite and
 non-atomic, one can choose a Borel set
$E_1$ satisfying
$$
B'(x_1,\rho_1)\subset E_1\subset B(x_1,\rho_1)\qquad\text{ and }\qquad
\mu(E_1)=\frac{\al_1}{t}.
$$
Let
$B(x_2,\rho_2)$ be the smallest ball such that $\mu(B(x_2,\rho_2)\setminus
E_1)\ge\dfrac{\al_2}{t}$ (since $\mu(\C)>\frac1t$, the measure of the remaining
part $\C\setminus E_1$ is still greater than
$\dfrac{1-\al_1}{t}\ge\dfrac{\al_2}{t}$).
Again for the corresponding open ball $B'(x_2,\rho_2)$,
we have $\mu(B'(x_2,\rho_2)\setminus
E_1)\le\dfrac{\al_2}{t}$, and therefore there exists a Borel set $E_2$
satisfying
$$
B'(x_2,\rho_2)\setminus E_1
\subset E_2
\subset B(x_2,\rho_2)\setminus E_1
\qquad\text{ and }\qquad
\mu(E_2)=\frac{\al_2}{t}.
$$
In general, for $i=3,4,\dots,N$, let
$B(x_i,\rho_i)$ be the smallest ball such that
$$
\mu\Bigl(B(x_i,\rho_i)\setminus
\bigcup_{\ell=1}^{i-1}E_\ell
\Bigr)\ge\frac{\al_i}{t},
$$
and let $E_i$ be a Borel set satisfying
$$
B'(x_i,\rho_i)\setminus
\bigcup_{\ell=1}^{i-1}E_\ell
\subset
E_i
\subset
B(x_i,\rho_i)\setminus
\bigcup_{\ell=1}^{i-1}E_\ell
\quad\text{ and }\quad
\mu(E_i)=\frac{\al_i}{t}.
$$
Put $E:=\bigcup_i E_i$.
Clearly
$$
\bigcup_i B'(x_i,\rho_i)
\subset
E
\subset
\bigcup_i B(x_i,\rho_i)
\qquad\text{ and }\qquad
\mu(E)=\frac1t.
$$
Now let us compare $T\nu$ to $t\,\sum_i\chi\ci{\C\setminus
B(x_i,2\rho_i)}
\cdot T\chi\ci{E_i}=:t\sigma$ outside $E$.
We have
$$
T\nu-t\s=\sum_i\f_i
$$
where
$$
\f_i=\al_i T\delta_{x_i}- t\,\chi\ci{\C\setminus B(x_i,2\rho_i)}\cdot
T\chi\ci{E_i}.
$$
Note now that
$$
\int_{\C\setminus E}|\f_i|d\mu\le \int_{\C\setminus
B(x_i,2\rho_i)}\bigl|T[\al_i
\delta_{x_i}-t\chi\ci{E_i}d\mu]\bigr|d\mu + \int_{B(x_i,2\rho_i)\setminus
B'(x_i,\rho_i)}
\al_i|T\delta_{x_i}|d\mu=:I+\al_i II\,.
$$
To estimate $I$, notice that it has the form $\int_{\C\setminus B(x,2\rho)}|T\eta|\,d\mu$ with the measure $\eta$ supported
by $B(x,\rho)$ and $\eta(\C) =0$. To estimate such an integral we put $\hat{r} :=\max [2\rho, R(x)]$ and split $\int_{\C\setminus
B(x,2\rho)}|T\eta|\,d\mu=
\int_{B(x,\hat{r})\setminus B(x,2\rho)}|T\eta|\,d\mu + \int_{\C\setminus B(x,\hat{r})}|T\eta|\,d\mu =:I_1+I_2$. The ntegral $I_2$
is estimated exactly as in Lemma 3.4 of [NTV2] because our measure is $2M$-Ahlfors for disks centered at $x$ with radii larger than
$\hat{r}$. Hence $I_2\leq ACM\|\eta\|\leq ACM\al_i$. On the other hand using the properties of the kernel of $T$ we see that
$$
I_1\leq  C\min [\frac{1}{2\rho},\frac{1}{R(x)}] \mu(B(x,\hat{r}))\|\eta\|\leq ACM\al_i\,.
$$
Hence $I\leq ACM\al_i$.

To estimate $II$ we notice that it has the form $\int_{B(x,2\rho)\setminus B(x,\rho)}|T\delta_x|\,d\mu$. This is almost the same
as $I_1$. Namely, $II \leq AC\min [\frac{1}{\rho},\frac{1}{R(x)}]\mu(B(x,2\rho))\leq \frac{AC\mu(B(x, 2 \max [R(x),\rho ]}{\max
[R(x),\rho ]}$. This is bounded by $ACM$ because our measure is $2M$-Ahlfors for disks centered at $x$ with radii larger than
$R(x)$. Finally $I +\al_i II \leq ACM\al_i$.

Finally we conclude that
$$
\int_{\C\setminus E}|T\nu-t\s|d\mu\le
ACM\sum_i\al_i=ACM,
$$
and thereby $|T\nu-t\s|\le ACMt$
everywhere on $\C\setminus E$, except, maybe, a set of
measure $\frac{1}{t}$. To accomplish the proof of the theorem,
we will show that for sufficiently large $B=B(C,M, \|T\|_2)$,
$$
\mu\{|\s|>B\}\le \frac2t.
$$
Then, combining all the above estimates, we shall get
$$
\mu\bigl\{x\in\C\,:\,|T\nu(x)|>(B+ACM)t\bigr\}\le\frac4t.
$$

We will apply the standard Stein-Weiss duality trick.
Assume that the inverse inequality
$\mu\{|\s|>B\}> \frac2t$
holds. Then either
$
\mu\{\s>B\}> \frac1t,
$
or
$
\mu\{\s<-B\}> \frac1t.
$
Assume for definiteness that the first case takes place
and choose some set $F\subset \C$ of measure exactly
$\frac1t$
such that $\s>B$ everywhere on $F$.
Then, clearly,
$$
\int_\C\s\chi\ci F d\mu>\frac{B}{t}.
$$
On the other hand, this integral can be computed as
$$
\sum_i \int_{\C} [T\chi\ci{E_i}]\cdot\chi\ci{F\setminus
B(x_i,2\rho_i)}\,d\mu
=\sum_i \int_{\C} \chi\ci{E_i}\cdot [T^{'}\chi\ci{F\setminus
B(x_i,2\rho_i)}]\,d\mu.
$$
Fix a point $x\in E_i\subset B(x_i,\rho_i)$. We will use again the property that $|K(x,y)|\leq \frac1{R(x)}$.
$$
|T^{'}\chi\ci{F\setminus B(x_i,2\rho_i)}(x)-T^{'}\chi\ci{F\setminus
B(x,\rho_i)}(x)|
\le
$$
$$
\int_{B(x_i,2\rho_i)\setminus B(x,\rho_i)}|K(y,x)|\,d\mu(y)\leq\frac{AC\mu(B(x, 3 \max [\rho_i, R(x) ]))}{\max [\rho_i, R(x) ]}
\leq ACM,
$$
because all disks centered at $x$ and of radii greater than $R(x)$ are $2M$-Ahlfors,
and therefore for every $x\in E_i\cap\supp\mu$,
$$
|T^{'}\chi\ci{F\setminus B(x_i,2\rho_i)}(x)|\le (T^{'})^\sharp \chi\ci
F(x)+ACM\le
2\cdot A\,\wt MT^{'}\chi\ci F(x)+ACM
$$
according to  Guy David's lemma. 
Hence
$$
\int_{\C} \s\chi\ci F d\mu\le ACM\mu(E)+2\cdot A\int_{\C}\chi\ci
E\cdot 
\wt MT^{'}\chi\ci F
d\mu.
$$
But the first term equals $\dfrac{ACM}{t}$ while the second one does
not exceed
$$
2\cdot 3^n\,
\|\chi\ci E\|\ci{L^2(\mu)}
\|\wt MT^{'}\chi\ci F\|\ci{L^2(\mu)}\le
\frac{2\cdot 3^n}{t}\|\wt M\|\ci{L^2(\mu)\to
L^2(\mu)}
\|T^{'}\|\ci{\!L^2(\mu){\to} L^2(\mu)\!}.
$$
Recalling that
$
\|T^{'}\|\ci{\!L^2(\mu){\to} L^2(\mu)\!}
=\|T\|\ci{\!L^2(\mu){\to} L^2(\mu)\!}
$, we see that one can take
$$
B=ACM+2\cdot 3^n\,\|\wt M\|\ci{L^2(\mu)\to
L^2(\mu)}\|T\|\ci{\!L^2(\mu){\to} L^2(\mu)\!}
$$
to get a contradiction. Since the norm $\|\wt M\|\ci{L^2(\mu)\to
L^2(\mu)}$ is bounded by some absolute constant (the constant in the
Marcinkiewicz interpolation theorem), we are done.

\bigskip

\centerline
{\bf References}

\bigskip

[CJS] R.~R.~Coifman, P.~W.~Jones and S.~Semmes,  Two elementary
proofs of the $L\sp 2$ boundedness of Cauchy integrals on Lipschitz
curves, J.~Amer.~Math.~Soc.~ 2, (1989), no.~3, 553--564.

[DM] G. David, P. Mattila, Removable sets for Lipschitz harmonic functions in the plane. 
Rev. Mat. Iberoamericana 16 (2000), no. 1, 137--215. 

[D1] G. David, Unrectifiable $1$-sets have vanishing analytic capacity. Rev. Mat. Iberoamericana 14 (1998), no. 2, 369--479.

[D2] G. David, Analytic capacity, Cauchy kernel, Menger curvature, and rectifiability. Harmonic analysis and partial differential  equations (Chicago, IL, 1996), 183--197, Chicago Lectures in Math., Univ. Chicago Press, Chicago, IL, 1999.

[D3] G. David, Analytic capacity, Calder\'on-Zygmund operators, and rectifiability. Publ. Mat. 43 (1999), no. 1, 3--25.

[Du] J. Doudziak, Vitushkin's conjecture for removal sets, Springer, 2010, 272 pp.

[GJ] J.~B~Garnett and P.~W.~Jones,  $BMO$ from dyadic $BMO$,
Pacific J.~Math.,  99, (1982), no.~2, 351--371.

[L] J-C. L\'eger, Menger curvature and rectifiability. Ann. of Math. (2) 149 (1999), no. 3, 831--869. 

[Ma] P. Mattila, Cauchy singular integrals and rectifiability in measures of the plane. Adv. Math. 115 (1995), no. 1, 1--34.

[MaP] P. Mattila, D. Preiss, David Rectifiable measures in ${R}\sp n$ and existence of principal values for singular integrals. 
J. London Math. Soc. (2) 52 (1995), no. 3, 482--496.

[MTV] J. Matheu, X. Tolsa, J. Verdera, The planar Cantor sets of zero analytic capacity and the local $T(b)$-Theorem. Preprint, Univ. Auton. Barcelona, 2001, pp. 1--12.

[M] M. Melnikov, Analytic capacity: a discrete approach and the curvature of measure. (Russian) Mat. Sb. 186 (1995), no. 6, 57--76.

[MP] P. Mattila, P.V. Paramonov, On geometric properties of harmonic $Lip_1$-capacity, Pac. J. math., 171 (1991), No. 2, 469--491.

[MV] M.S. Melnikov, J. Verdera, A geometric proof of $L^2$ boundedness of the Cauchy integral on lipschitz graph. 
Intern. Math. Res. Notices, 7 (1995), 325--331.

[MMV] P. Mattila, M. Melnikov, J. Verdera, The Cauchy integral, analytic capacity, and uniform rectifiability. 
Ann. of Math. (2) 144 (1996), no. 1, 127--136. 

[NTV1] F. Nazarov, S. Treil, A. Volberg, Cauchy integral and Calder\'on-Zygmund operators on nonhomogeneous spaces, 
Intern. Math. Res. Notices, 15 (1997), 703--726.

[NTV2] F. Nazarov, S. Treil, A. Volberg, Weak type estimates and Cotlar's inequality for Calder\`on-Zygmund operators on nonhomogeneous spaces, Intern. Math. Res. Notices, 9 (1998), 463--487.

[NTV3] F. Nazarov, S. Treil, A. Volberg, The $Tb$-theorem on non-homogeneous spaces. Acta Math. 190 (2003), no. 2, 151--239.

[To] X. Tolsa,   Analytic capacity, the Cauchy transform, and
non-homogeneous Calder\'on-Zygmund theory,
Manuscript, 2012, http://www.mat.uab.cat/$\sim$ xtolsa/llibreweb.pdf.

[XT1] X.~Tolsa, {\it Painlev\'e's problem and analytic capacity}, Proceedings of the 7th International Conference on Harmonic Analysis and Partial Differential Equations
El Escorial, Madrid (Spain), June 21--25, 2004, Collect. Math. (2006), 89--125.

[XT2] X.~Tolsa, {\it Bilipschitz maps, analytic capacity, and the Cauchy integral.} Ann. of Math. (2) 162 (2005), no. 3, 1243–1304.

[XT3] X. Tolsa, {\em Painlev\'{e}'s problem and the
semiadditivity of analytic capacity}, Acta Math. 190:1 (2003), 105--149.

[Vo] A. Volberg, Calder\'on--Zygmund capacities and operators on non-homogeneous spaces. CBMS lecture series, AMS, v. 100, 2003, pp. 165.

\end{document}